\documentclass[12pt]{article}
\usepackage{graphicx}
\makeatletter

\def\URA#1{\begin{picture}(3,2)\put(0,0){\vector(3,2){3}}%
        \put(1.5,.5){${\scriptstyle{#1}}$}\end{picture}}

\def\DRA#1{\begin{picture}(3,2)\put(0,2){\vector(3,-2){3}}%
        \put(1.5,1.2){${\scriptstyle{#1}}$}\end{picture}}
\def\DLA#1{\begin{picture}(3,2)\put(3,2){\vector(-3,-2){3}}%
        \put(1.5,.5){${\scriptstyle{#1}}$}\end{picture}}
\def\LRA#1#2{\@tempdimb=\c@enumiv\@tempdima%
   \vcenter{\offinterlineskip\halign{##\cr%
   \hfil${\scriptstyle{#1}}$\hfil\crcr%
   \hbox to \@tempdimb{\rightarrowfill}\cr%
   \noalign{\kern-1ex}%
   \hbox to \@tempdimb{\leftarrowfill}\cr%
   \hfil${\scriptstyle{#2}}$\hfil\crcr}}}
\def\RA#1{\@tempdimb=\c@enumiv\@tempdima\vbox{\offinterlineskip%
   \halign{##\cr\hfil${\scriptstyle {#1}}$\hfil\crcr%
   \hbox to \@tempdimb{\rightarrowfill}\cr}}}
\def\LA#1{\@tempdimb=\c@enumiv\@tempdima\vbox{\offinterlineskip%
   \halign{##\cr\hfil${\scriptstyle {#1}}$\hfil\crcr%
   \hbox to \@tempdimb{\leftarrowfill}\cr}}}

\def\DA#1{\strut\begin{picture}(1,2)\put(.5,2){\vector(0,-1){2}}%
        \put(.7,.8){${\scriptstyle{#1}}$}\end{picture}}
\def\diag{\leavevmode\bgroup\setcounter{enumiv}{1}%
   \unitlength1em \@tempdima3em \def\\{\crcr&}\vbox\bgroup%
   \def\multicolumn##1##2{\multispan##1\setcounter{enumiv}{##1}%
   \hfil{##2}\hfil\setcounter{enumiv}{1}}
   \offinterlineskip\halign\bgroup\vrule height.8em depth.7em %
   width0pt##&&\hfil${\displaystyle{##}}$\hfil\cr&}
\def\enddiag{\crcr\egroup\egroup\egroup}
\makeatother
\if@twoside \oddsidemargin 0pt \evensidemargin 0pt \marginparwidth 85pt
\else \oddsidemargin 0pt \evensidemargin 0pt
\fi
\advance\textheight by \topskip
\font\symbolfont=msbm10
\font\teneu=eufm10
\font\egteu=eufm8
\def\dn#1{\mathchoice{\hbox{\teneu #1}}{\hbox{\teneu #1}}%
   {\hbox{\egteu #1}}{\hbox{\egteu #1}}}

\def\ccc{\hbox{\symbolfont C}}

\def\fff{\hbox{\symbolfont F}}

\def\ggg{\hbox{\symbolfont G}}

\def\zzz{\hbox{\symbolfont Z}}
\def\ppp{\hbox{\symbolfont P}}
\def\qqq{\hbox{\symbolfont Q}}

\def\ttt{\hbox{\symbolfont T}}

\def\aa{\mathop{\cal A}\nolimits}
\def\bb{\mathop{\cal B}\nolimits}
\def\cc{\mathop{\cal C}\nolimits}
\def\dd{\mathop{\cal D}\nolimits}

\def\ii{\mathop{\cal I}\nolimits}
\def\kk{\mathop{\cal K}\nolimits}
\def\oo{\mathop{\cal O}\nolimits}

\def\qq{\mathop{\cal Q}\nolimits}

\def\pd{\mathop{\rm pd}\nolimits}
\def\height{\mathop{\rm ht}\nolimits}
\def\id{\mathop{\rm id}\nolimits}
\def\gl{\mathop{\rm gl.dim}\nolimits}

\def\height{\mathop{\rm ht}\nolimits}
\def\Spec{\mathop{\rm Spec}\nolimits}

\def\spep#1{\mathop{{}^{\bullet}\strut\kern-.1em{#1}}\nolimits}

\def\Ext{\mathop{\rm Ext}\nolimits}
\def\hom{\mathop{\rm Hom}\nolimits}
\def\Rhom{\mathop{\rm {\bf R}Hom}\nolimits}
\def\Lotimes{\mathop{\stackrel{\bf L}{\otimes}}\nolimits}
\def\endm{\mathop{\rm End}\nolimits}
\def\aut{\mathop{\rm Aut}\nolimits}

\def\out{\mathop{\rm Out}\nolimits}
\def\add{\mathop{\rm add}\nolimits}
\def\Mod{\mathop{\rm Mod}\nolimits}
\def\soc{\mathop{\rm soc}\nolimits}
\def\top{\mathop{\rm top}\nolimits}

\def\Flat{\mathop{\rm Flat}\nolimits}
\def\Pr{\mathop{\rm Pr}\nolimits}
\def\pr{\mathop{\rm pr}\nolimits}

\def\cm{\mathop{\rm CM}\nolimits}
\def\ann{\mathop{\rm ann}\nolimits}
\def\tr{\mathop{\rm Tr}\nolimits}

\def\Ker{\mathop{\rm Ker}\nolimits}

\def\Im{\mathop{\rm Im}\nolimits}

\def\mod{\mathop{\rm mod}\nolimits}
\def\grmod{\mathop{\rm grmod}\nolimits}
\def\Coh{\mathop{\rm Coh}\nolimits}
\def\ref{\mathop{\rm ref}\nolimits}
\def\flmod{\mathop{\rm fl}\nolimits}
\def\rank{\mathop{\rm rank}\nolimits}

\def\add{\mathop{\rm add}\nolimits}

\def\cc{\mathop{\cal C}\nolimits}

\def\tor{\mathop{\rm Tor}\nolimits}
\def\Mod{\mathop{\rm Mod}\nolimits}

\def\Ker{\mathop{\rm Ker}\nolimits}

\def\depth{\mathop{\rm depth}\nolimits}
\def\pd{\mathop{\rm pd}\nolimits}

\def\id{\mathop{\rm id}\nolimits}
\def\gl{\mathop{\rm gl.dim}\nolimits}

\def\height{\mathop{\rm ht}\nolimits}
\def\Spec{\mathop{\rm Spec}\nolimits}
\def\Max{\mathop{\rm Max}\nolimits}
\def\Supp{\mathop{\rm Supp}\nolimits}
\def\GL{\mathop{\rm GL}\nolimits}
\def\SL{\mathop{\rm SL}\nolimits}
\def\tri{\mathop{\rm T}\nolimits}
\def\ma{\mathop{\rm M}\nolimits}
\def\ext{\mathop{\rm Ext}\nolimits}

\def\op{\mathop{\rm op}\nolimits}
\def\b{\mathop{\rm b}\nolimits}

\def\qq{\mathop{\cal Q}\nolimits}

\def\length{\mathop{\rm length}\nolimits}

\def\dpic{\mathop{\rm DPic}\nolimits}
\def\auteq{\mathop{\rm Auteq}\nolimits}

\def\resdim#1#2{\mathop{{#1}\mbox{-{\rm dim}}\strut\kern2pt {#2}}\nolimits}

\def\timesl{\begin{picture}(16,8)\put(3,-1){\line(1,1){10}}\put(13,9){\line(0,-1){10}}\put(13,-1){\line(-1,1){10}}\end{picture}}
\def\ltimes{\begin{picture}(16,8)\put(3,-1){\line(1,1){10}}\put(3,9){\line(0,-1){10}}\put(13,-1){\line(-1,1){10}}\end{picture}}
\def\plim{\strut\kern-.5em\mathop{\raise-.2em\hbox{$\def\arraystretch{0}\begin{array}{c}\lim\\
\leftarrow\end{array}$}}\strut\kern-.5em}
\def\ilim{\strut\kern-.5em\mathop{\raise-.2em\hbox{$\def\arraystretch{0}\begin{array}{c}\lim\\
\rightarrow\end{array}$}}\strut\kern-.5em}
\def\XZ{0}
\def\XA{1}
\def\XB{2}
\def\XBA{2.1}
\def\XBB{2.2}
\def\XBC{2.3}
\def\XBD{2.4}
\def\XBE{2.5}
\def\XBF{2.6}
\def\XBG{2.7}
\def\XBH{2.8}
\def\XBI{2.9}
\def\XBJ{2.10}
\def\XC{3}
\def\XCA{3.1}
\def\XCB{3.2}
\def\XCC{3.3}
\def\XCD{3.4}
\def\XCE{3.5}
\def\XCF{3.6}
\def\XCG{3.7}
\def\XCH{3.8}
\def\XCI{3.9}
\def\XCJ{3.10}
\def\XCK{3.11}
\def\XCL{3.12}
\def\XCM{3.13}
\def\XCN{3.14}
\def\XD{4}
\def\XDA{4.1}
\def\XDB{4.2}
\def\XDC{4.3}
\def\XDD{4.4}
\def\XE{5}
\def\XEA{5.1}
\def\XEB{5.2}
\def\XEC{5.3}
\def\XED{5.4}
\def\XEE{5.5}
\def\XEF{5.6}
\def\XEG{5.7}
\def\XF{6}
\def\XFA{6.1}
\def\XFB{6.2}
\def\XFC{6.3}
\def\XFD{6.4}
\def\XFE{6.5}
\def\XFF{6.6}
\def\XFG{6.7}
\def\XFH{6.8}
\def\XFI{6.9}
\def\XFJ{6.10}
\def\XFK{6.11}
\def\XFL{6.12}
\def\XFM{6.13}
\def\XFN{6.14}
\def\XG{7}
\def\XGA{7.1}
\def\XGB{7.2}
\def\XGC{7.3}
\def\XH{8}
\def\XHA{8.1}
\def\XHB{8.2}
\def\XHC{8.3}
\def\XHD{8.4}
\def\XHE{8.5}
\def\XHF{8.6}
\def\XHG{8.7}
\def\XHH{8.8}
\def\XHI{8.9}
\def\XHJ{8.10}
\def\XHK{8.11}
\def\XHL{8.12}
\def\XHM{8.13}
\def\XHN{8.14}
\def\XHO{8.15}
\def\XHP{8.16}
\def\XHQ{8.17}
\def\XHR{8.18}
\def\XHS{8.19}
\def\tilt{\mathop{\rm tilt}\nolimits}

\begin{document}
\begin{center}
{\Large Fomin-Zelevinsky mutation and tilting modules over Calabi-Yau algebras}

\vskip1em
Osamu Iyama and Idun Reiten
\end{center}

\vskip1em{\small {\sc Abstract.} We say that an algebra $\Lambda$ over a commutative noetherian ring $R$ is Calabi-Yau of dimension $d$ ($d$-CY) if the shift functor $[d]$ gives a Serre functor on the bounded derived category of the finite length $\Lambda$-modules. We show that when $R$ is $d$-dimensional local Gorenstein the $d$-CY algebras are exactly the symmetric $R$-orders of global dimension $d$. We give a complete description of all tilting modules of projective dimension at most one for 2-CY algebras, and show that they are in bijection with elements of affine Weyl groups, preserving various natural partial orders. We show that there is a close connection between tilting theory for 3-CY algebras and the Fomin-Zelevinsky mutation of quivers (or matrices). We prove a conjecture of Van den Bergh on derived equivalence of non-commutative crepant resolutions.}

\vskip1em{\small
0. Introduction

1. Cluster algebras and Calabi-Yau conditions

2. Preliminaries on module-finite algebras

3. Calabi-Yau algebras and symmetric orders

4. Construction of tilting modules

5. Mutation on tilting modules

6. 2-Calabi-Yau algebras and affine Weyl groups

7. 3-Calabi-Yau algebras and cluster algebras

8. Non-commutative crepant resolutions}

\vskip1.5em{\bf\XZ. Introduction}

Let $\Lambda$ be an algebra over a commutative noetherian ring $R$, which is finitely generated as an $R$-module. Sometimes we assume that $R$ in addition satisfies one or more of the following conditions: local, complete, Gorenstein, normal. This paper deals mainly with algebras $\Lambda$ which are Calabi-Yau of dimension $d$, called $d$-CY algebras for short. This means that the shift functor $[d]$ gives a Serre functor on the bounded derived category of the finite length $\Lambda$-modules. The main aspects of $d$-CY algebras which we investigate are the following. We deduce interesting properties of $d$-CY algebras of a ring theoretic and module theoretic nature. In particular, we show that when $R$ is $d$-dimensional local Gorenstein, the $d$-CY algebras are exactly the symmetric $R$-orders of global dimension $d$.

A central part of our investigations deals with tilting modules, mainly those of projective dimension at most one. We give a complete description of all such tilting modules for 2-CY algebras, and show that they are in bijection with elements of affine Weyl groups, preserving various natural partial orders. We also investigate tilting theory for 3-CY algebras, where we show that there is a close connection with the Fomin-Zelevinsky mutation of quivers (or matrices), which they introduced in connection with their definition of cluster algebras. Tilting theory for 3-CY algebras also turns out to have a close connection with the theory of non-commutative crepant resolutions of Van den Bergh, and we prove some of our results in this more general setting (actually for a generalization of the definition of Van den Bergh). In particular, we prove a conjecture of Van den Bergh on derived equivalence of non-commutative crepant resolutions in \XHH.

Main examples of $d$-CY algebras are skew group algebras $S*G$, where $S$ is the formal power series ring in $d$ variables over a field $K$ of characteristic zero and $G$ is a finite subgroup of $\SL_d(K)$. For the 2-dimensional case, these algebras have been a major object in several branches of mathematics. The invariant ring $S^G$ has only finitely many indecomposable (maximal) CM (=Cohen-Macaulay) modules, and $S*G$ is the endomorphism ring of an additive generator $S$ of the category $\cm S^G$ of (maximal) CM $S^G$-modules satisfying $\gl S*G=2$ [He][A3]. Consequently, the Auslander-Reiten quiver of $S^G$ coincides with the quiver of $S*G$, and also with the McKay quiver of $G$, and they are the double of an extended Dynkin diagram [A3]. This is closely related to the minimal resolution $X$ of the quotient singularity $S^G$. A one-one correspondence between indecomposable non-free CM $S^G$-modules and irreducible exceptional curves of $X$ was constructed in [AV]. More strongly, by the McKay correspondence [Mc] explained geometrically by Gonzalez-Sprinberg and Verdier [GV], the McKay graph of $G$ coincides with the dual graph of irreducible exceptional curves of $X$. Later, Kapranov and Vasserot [KV] reformulated McKay correspondence in terms of derived equivalences between $S*G$ and $X$. The tilting modules we construct for 2-CY algebras, and hence for the algebras $S*G$ in dimension 2, give autoequivalences of $\dd^{\b}(\mod(S*G))$, and correspond to twist functors on $\dd^{\b}(\Coh X)$ constructed by Seidel-Thomas [ST] under McKay correspondence. It is natural to ask about the geometric meaning of the tilting modules for arbitrary $d$-CY algebras (including $S*G$) which we construct in section \XD. Also the category $\cm S^G$ for the case $d>2$ is studied in [I3,4][IY].

Although our study of Calabi-Yau algebras in this paper is done from a purely algebraic viewpoint, it is closely related to the study of the usual `geometric' Calabi-Yau varieties. 
%Bondal and Orlov conjectured in [BO] that 3-dimensional smooth varieties related by a flop are derived equivalent, especially, any birational 3-dimensional Calabi-Yau varieties are derived equivalent. Bridgeland [Bri1] proved these conjectures
By using the method of Fourier-Mukai transformations, Bridgeland [Bri1] proved a conjecture of Bondal and Orlov [BO], which states that any birational 3-dimensional Calabi-Yau varieties are derived equivalent. It was motivated by Bridgeland's proof and also by the 3-dimensional McKay correspondence due to Bridgeland-King-Reid [BKR] that Van den Bergh [Va1,2] introduced the concept of non-commutative crepant resolutions (NCCR for short), and gave a new proof of Bridgeland's theorem. A typical example of NCCR is given again by skew group rings. It was in this connection that Van den Bergh proposed an analogue of a conjecture of Bondal and Orlov, i.e. all crepant resolutions (including NCCR) of a normal Gorenstein domain are derived equivalent, and proved this conjecture for 3-dimensional terminal singularities. In section \XH, we prove the non-commutative part of this conjecture for the 3-dimensional case in a more general form, i.e. all NCCR of a module-finite algebra over a 3-dimensional normal Gorenstein domain are derived equivalent. Our method using tilting modules is surprisingly simple. There are here also interesting relationships with the maximal 1-orthogonal modules of [I3,4] and the maximal rigid modules of [GLS].

In the philosophy of non-commutative algebraic geometry, we regard algebras with finite global dimension as an algebraic/non-commutative analogue of smooth varieties. Although there is a classical Cohen's structure theorem for regular algebras in the commutative situation, finiteness of global dimension seems not to be sufficient in the non-commutative situation. Known successful theory of algebras with low global dimension is obtained by Reiten and Van den Bergh [RV1,2], Artin and Schelter [AS], and Rump [Rum1,2] by giving additional assumptions, e.g. order, Auslander condition, and so on. We hope that our CY algebras will provide a source of non-commutative regular algebras.

The theory of tilting modules originated in the representation theory of finite dimensional algebras, where a substantial theory has been developed around this concept. There have also been generalizations to other classes of algebras, but not many concrete examples of tilting modules beyond the finite dimensional case. When $T'$ is an almost complete tilting $\Lambda$-module, then there are up to isomorphism at most two complements [RS]. For a finite dimensional algebra it is never the case that there are always two complements. However, this holds for $d$-CY algebras for $d=2,3$, and constitutes an important part of our investigations. In particular, it is an essential property for trying to model the Fomin-Zelevinsky mutation of quivers of 3-CY algebras using tilting $\Lambda$-modules. A similar type of program was carried out for acyclic cluster algebras, leading to the notions of cluster categories and cluster tilted algebras [BMRRT][BMR1] (see [CCS] for the $A_n$ case), along with related work in [GLS]. Actually, possible connections with cluster algebras was one of the motivations for investigating tilting theory for 3-CY algebras. Another interesting property of tilting theory for 3-CY algebras is that an equivalence of derived categories induced by a sequence of equivalences determined by tilting modules is actually induced by a tilting module. This is usually not the case for finite dimensional algebras.

The tilting theory for 2-CY algebras is also quite interesting from a very different point of view. We give a description of the tilting modules (of projective dimension at most 1), which all turn out to be ideals. The 2-CY algebras have valued quivers which are given by generalized extended Dynkin diagrams. We establish a bijection between the tilting modules and the elements of the associated affine Weyl groups. There is in general a partial order on the tilting modules [RS], and in addition there is in this case the reverse inclusion order of ideals. Under our bijection they correspond to the right order and the Bruhat order on the affine Weyl group.

We now describe the content of each section. We start with some background material and motivation in section \XA, with particular emphasis on the importance of 2-dimensional and 3-dimensional Calabi-Yau properties for modelling quiver mutation and other essential ingredients in the definition of cluster algebras. In section \XB\ we collect basic results on some central concepts  from commutative algebra, like dimension and depth, along with properties of reflexive and derived equivalences. We characterize the $d$-CY algebras as symmetric orders of global dimension $d$ in section \XC, and investigate a related class which we call $d$-CY$^-$ algebras, which coincides with the class of symmetric orders. In section \XD\ we start with the tilting $\Lambda$-module $\Lambda =\bigoplus_{i=1}^nP_i$, where the $P_i$ are non-isomorphic indecomposable projective $\Lambda$-modules. We describe the indecomposable objects in the bounded derived category which can replace $P_i$ to give a tilting complex, and in particular a tilting module. Our emphasis is on CY algebras. In section \XE\ we give some basic properties of tilting theory for module-finite algebras, most of which are analogs of properties of finite dimensional algebras, and give stronger results for complements of almost complete tilting modules over CY algebras. In section \XF\ we investigate the structure of tilting modules for 2-CY algebras and their connection with affine Weyl groups. In section \XG\ we specialize to 3-CY algebras, and investigate the connection between Fomin-Zelevinsky mutation and tilting theory. The connection between tilting modules, reflexive modules and the non-commutative crepant resolutions of Van den Bergh is discussed in section \XH.

Parts of the results in this paper have been presented at conferences in Beijing, Banff, Tokushima, Hanamako, Trieste, Hannover and Oberwolfach.
There is related work on various aspects on CY algebras by for example Berger-Taillefer [BT], Bocklandt [Boc], Braun [Bra], Brown-Gordon-Stroppel [BGS], Chuang-Rouquier [CRo], Ginzburg [G], Rickard, Van den Bergh [Va3].

\vskip.5em{\bf Notation }
For a noetherian ring $\Lambda$, we denote by $J_\Lambda$ the Jacobson radical of $\Lambda$. By a module we mean a left module. We denote by $\Mod\Lambda$ the category of $\Lambda$-modules, by $\mod\Lambda$ the category of finitely generated $\Lambda$-modules, and by $\flmod\Lambda$ the category of $\Lambda$-modules of finite length. These categories are abelian. Moreover, we denote by $\Pr\Lambda$ the category of projective $\Lambda$-modules, and by $\pr\Lambda$ the category of finitely generated projective $\Lambda$-modules.
For an additive category $\cc$, we denote by $J_{\cc}$ the Jacobson radical of $\cc$.

For a commutative noetherian ring $R$ we denote by $\Spec R$ the set of prime ideals of $R$, and by $\Max R$ the set of maximal ideals of $R$. For $\dn{p}\in\Spec R$, we denote by $R_{\dn{p}}$ the localization of $R$ at $\dn{p}$, and by $\widehat{R}_{\dn{p}}$ the completion of $R$ at $\dn{p}$. For $X\in\Mod R$, we put $X_{\dn{p}}:=X\otimes_RR_{\dn{p}}$ and $\widehat{X}_{\dn{p}}:=X\otimes_R\widehat{R}_{\dn{p}}$. We denote by $\Supp{}_RX$ the support of $X$. When $R$ is a local ring the maximal ideal $\dn{p}$, we often write $\widehat{X}:=\widehat{X}_{\dn{p}}$. We denote by $(-)^*$ the functor $\hom_R(-,R):\mod R\to\mod R$.

\vskip.5em{\sc Acknowledgements }

In addition to the work of Iyama [I1,2] and Van den Bergh [Va1,2], some of this work was inspired by the 2004 Edinburgh meeting on ``derived categories, quivers and strings'' organized by A. King and M. Douglas. The second author would like to thank the participants of the conference, in particular D. Berenstein, T. Bridgeland, J. Chuang, B. Keller, A. King, J. Rickard, R. Rouquier and M. Van den Bergh for interesting conversations.

A major part of this work was done while the first author visited Trondheim in February 2005, with improvements during his visit in March 2006. He would like to thank the people at NTNU for hospitality and stimulating discussions. He also would like to thank R.-O. Buchweitz, A. Ishii, B. Keller, H. Krause, J. Miyachi, R. Takahashi, J. Rickard, R. Rouquier, W. Rump, T, Shoji, K. Yoshida and Y. Yoshino for valuable comments.

\vskip1.5em{\bf\XA. Cluster algebras and Calabi-Yau conditions}

In this section we recall some basic facts on cluster algebras, and
discuss the relevance of Calabi-Yau type conditions for being able
to model some of the concepts from the definition of a cluster
algebra in a categorical/module theoretical way.

We recall the basic definitions from [FZ1] in a generality
suitable for our purpose. We give the definitions in terms of
quivers rather than matrices. Let $\qq$ be a finite quiver with no
loops. If $\qq$ has oriented cycles of length two, we can associate
with $\qq$ a quiver $\bar{\qq}$ obtained by removing all pairs
%$\begin{smallmatrix} {\rightarrow} & \\ {\leftarrow} & \end{smallmatrix}$
$\def\arraystretch{.2}\left.\begin{array}{c}\longrightarrow\\ \longleftarrow\end{array}\right.$. Let
$F=\qqq(x_1,\cdots, x_n)$ where $\qqq$ denotes the
rational numbers and $x_1, \cdots, x_n$ are indeterminates, and let
$\qq$ be a finite quiver with $n$ vertices and no oriented cycles of
length at most two. The pair $(\underline{x},\qq)$ where $\underline{x}
=\{x_1, \cdots, x_n\}$ is called a seed. Let $b_{ij}$
denote the number of arrows from $i$ to $j$ in $\qq$, interpreted as minus
the number of arrows from $j$ to $i$ in $\qq$ if $b_{ij} <0$. For each
$k=1, \cdots, n$, the mutated quiver $\qq'=\mu_k(\qq)$ is defined as
follows, where now $b_{ij}'$ denotes the number of arrows from $i$
to $j$ in $\qq'$:
\begin{itemize}
\item[] $b_{ik}' = -b_{ik}$ and $b_{kj}' = -b_{kj}$.

\item[] If $b_{ik},\ b_{kj}>0$ (resp. $b_{ik},\ b_{kj}<0$), then $b_{ij}'=b_{ij}+b_{ik}b_{kj}$ (resp. $b_{ij}'=b_{ij}-b_{ik}b_{kj}$).

\item[] $b_{ij}' =b_{ij}$ otherwise.
\end{itemize}

There is also defined a mutation
$\mu_k(\underline{x},\qq)=(\underline{x}',\qq')$ of seeds, where\\
$\underline{x}'=\{x_1,\cdots,x_k',\cdots,x_n\}$ and
$x_kx_k'=\prod_{b_{ik}>0}x_i^{b_{ik}}
+\prod_{b_{ik}<0}x_i^{-b_{ik}}$

Starting with a seed $(\underline{x},\qq)$ and applying sequences of
mutations in all directions, we obtain a collection of seeds. The
\emph{clusters} are by definition the $n$ element subsets
$\underline{x},\ \underline{x}',\cdots$, and the \emph{cluster
  variables} are the union of the elements in the clusters. The
associated cluster algebra is the subring of $F$ generated by
the cluster variables.

Mutation of quivers is defined without reference to any clusters or
cluster variables. It is an interesting problem to identify classes
of quivers where we can `lift' mutation of quivers to `mutation' of
algebras associated with the quivers. Secondly, we would want a
more global model, by finding some category $\mathcal{C}$ with
a special type of objects being the anolog of clusters, and where
the indecomposable summands are the anolog of cluster variables.
The relevant objects $T$ in $\mathcal{C}$ should have a direct sum
decomposition $T=\bigoplus_{i=1}^nT_i$ into $n$
non-isomorphic indecomposable summands, where $n$ is the number of
vertices in the quiver. For each $k=1,\cdots,n$ there should be a
unique indecomposable object $\check{T}_k\ {\not\simeq}\ T_k$ such that
$T'=T/T_k \oplus\check{T}_k$ is an object of the relevant type.
Also there should be some nice way of connecting $T_k$ and
$\check{T}_k$, and the endomorphism algebras
$\endm_{\mathcal{C}}(T)$ and $\endm_{\mathcal{C}}(T')$ should be related via `mutation' of algebras.

Inspired by the connection with quiver representations given in
[MRZ], such a program was carried out for acyclic cluster
algebras in [BMRRT], [BMR1,2]. A cluster algebra
is said to be acyclic if the quiver in some seed has no oriented
cycles. Here the global approach was done first, and the crucial
category was the cluster category ${\mathcal C}_{\qq}$ associated with
a finite quiver $\qq$ without oriented cycles, via the path algebra
$K{\mathcal Q}$ for a field $K$. The central objects $T$ in ${\mathcal C}_{\qq}$,
called (cluster-)tilting objects, are induced by tilting modules
over finite dimensional hereditary algebras in the derived
equivalence class of $K{\mathcal Q}$. The associated endomorphism algebras
$\endm_{{\mathcal C}_{\qq}}(T)$  are called cluster-tilted
algebras [BMR1], and it was shown later that the passage from
$\endm_{{\mathcal C}_{\qq}}(T)$ to
$\endm_{{\mathcal C}_{\qq}}(T')$ provides a lifting of the
corresponding quiver mutation [BMR2]. In addition to providing
interesting connections with cluster algebras, including a framework
for obtaining results on cluster algebras, the class of
cluster-tilted algebras is also interesting in itself.

The cluster category ${\mathcal C}_{\qq}$ associated to $H=K{\mathcal Q}$ is by
definition the factor category
${\mathcal C}_{\qq}=\dd^{\b}(H)/ \tau^{-1}[1]$, where
$\dd^{\b}(H)$ denotes the bounded derived category of the
category of finitely generated $H$-modules, $\tau$ the AR-
translation on $\dd^{\b}(H)$ and $[1]$ the shift functor. The
category ${\mathcal C}_{\qq}$ is known to be triangulated [K3] and it
has Serre functor $[2]$, since $\tau_{\cc}=[1]$ is an AR translation on ${\mathcal C}_{\qq}$.
Hence it is Calabi-Yau of dimension $2$, that is, there is for $M$
and $N$ in ${\mathcal C}_{\qq}$ a functorial isomorphism $D(\hom(M,N))
\simeq \hom(N,M[2])$, where $D$ is the ordinary duality. That the (cluster-)tilting theory in
${\mathcal C}_{\qq}$ works  very nicely, including unique exchange of
an indecomposable summand of a tilting object in ${\mathcal C}_{\qq}$
to get another one giving rise to a new tilting object, is
intimately related to the 2-CY property. Actually it is a consequence of the 2-CY
property [IY] (see also [I5]).

In order to prove lifting of the quiver mutation to the class of
cluster-tilted algebras, an essential ingredient was proving
reduction to three simple modules [BMR2]. To be able to lift
quiver mutation in this case has to do with establishing a close
relationship between arrows and relation spaces, in other words,
between $\Ext^1$- and $\Ext^2$-groups (see [Bon2]. Here a 3-CY
type property of such cluster-tilted algebras would be useful, since
this would imply $D\Ext^2(S,S') \simeq \Ext^1(S',S)$ for simple
modules $S$ and $S'$. Actually, it has been shown in [KR] that
the cluster-tilted algebra $\Gamma=\endm_{\cc_{\qq}}(T)$
is Gorenstein of dimension at most 1, with the stable category
$\underline{\cm}\Gamma$ of CM modules being 3-CY.
There is a more general result starting with $\cc$ being 2-CY [KR].
So we see that having 2-CY and 3-CY type properties around is
useful for modelling some central ingredients of cluster algebra
theory.

The quivers we deal with in this paper are coming from algebras
satisfying some 3-CY conditions. They are certain algebras
$\Lambda$ which are finitely generated as modules over a commutative
noetherian ring of Krull dimension $3$, such that the
bounded derived category $\dd^{\b}(\flmod\Lambda)$ is 3-CY.
In this setting we are able to lift mutations of
quivers to such 3-CY algebras. Also, under some additional
assumptions, the relevant category $\cc$ is $\mod \Lambda$, or rather
the subcategory $\ref \Lambda$ of reflexive $\Lambda$-modules with
the relevant objects being the tilting $\Lambda$-modules of projective dimension at most one.
In this case there is some trace of a 2-CY property associated
with the projective resolution of simple modules.

So in this case we have the `opposite' of the case of the categories associated with acyclic cluster algebras; that we have 3-CY and a trace of 2-CY. 
When we have the 3-CY algebra $S*G$ where $S=K[[x,y,z]]$ with the field $K$ of characteristic zero and $G$ is a finite subgroup of $\SL_3(K)$ acting freely, there is here a close relationship to the stable category $\underline{\cm}R$ of (maximal) CM modules for $R=K[[x,y,z]]^{G}$. 
The category $\underline{\cm}R$ is
2-CY, as follows from work of Auslander [A2;III]. For if $R$ is a commutative complete local Gorenstein isolated singularity of
dimension $d$, then the AR-translation $\tau$ is $\Omega^{2-d}$,
which is $[d-2]$ in $\underline{\cm}R$, and hence $\underline{\cm}R$ is $(d-1)$-CY.
More generally, the stable category of lattices over symmetric orders have the same property (see [A2] for definitions and results).

Other classes of cluster algebras are investigated from the
modelling point of view in [GLS]. Here they deal with the case
`with coefficients'. But at the same time they obtain results on the
`no coefficients' case, which is here natural to compare with, from
the point of view of CY-conditions. Let $\Lambda$ be the
preprojective algebra of a Dynkin diagram and $\underline{\mod}
\Lambda$ the stable category of the finitely generated
$\Lambda$-modules. Here the AR-translation $\tau$ is $\Omega^{-1}$
(see [AR2;3.1,1.2][K3;8.5]), so that $\underline{\mod}\Lambda$ is 2-CY.
When $\Lambda$ is of finite representation type, $\underline{\mod}
\Lambda$ is equivalent to a cluster category [BMRRT]. The
special objects  in $\underline{\mod}\Lambda$ are, like in the
cluster category, objects $C$ maximal with the property
$\Ext^1_{\underline{\mod}\Lambda}(C,C)=0$. Then $C$ can be lifted to $\Lambda\oplus C$ in
$\mod \Lambda$, which has a similar property. In this case one also
has the algebras $\Gamma = \endm_{\Lambda}(\Lambda\oplus C)$,
in addition to the factor algebras $\underline{\endm}(C)$. They
have global dimension 3, and have some trace of being 3-CY, as shown
in [GLS]; see also [KR]. There are also other similarities
between our work and that of [GLS], with respect to the role of
tilting modules.

\vskip1.5em{\bf\XB. Preliminaries on module-finite algebras }

The main focus in this paper is on $d$-CY algebras and related algebras, especially for $d=2,3$. In this section we give some useful background material on the module theory for module-finite $R$-algebras, for a noetherian  commutative ring $R$. First we discuss dimension, depth, global dimension and their relationships. Second, assuming that $R$ is in addition a normal domain, we discuss reflexive modules and reflexive equivalence of algebras, in particular that being a symmetric algebra is preserved under reflexive equivalence. Finally, we give some results on derived categories needed in section 3.

Let $\Lambda$ be a module-finite $R$-algebra and $M\in\mod\Lambda$. Put
\[\dim{}_RM:=\dim(R/\ann{}_RM)\]
where $\ann{}_RM$ is an annihilator of the $R$-module $M$. Since the value of $\dim{}_RM$ is independent of the choice of central subring $R$ of $\Lambda$, we denote it by $\dim M$.

Now assume that $R$ is local. Put
\[\depth{}_RM:=\inf\{i\ge0\ |\ \ext^i_R(R/J_R,M)\neq0\}.\]
Then $\depth{}_RM$ coincides with the maximal length of $M$-regular sequences [Ma;16.7]. By a result of Goto-Nishida [GN2;3.2], we have an equality
\[\depth{}_RM=\inf\{i\ge0\ |\ \ext^i_\Lambda(\Lambda/J_\Lambda,M)\neq0\}.\]
In particular, the value of $\depth{}_RM$ is independent of the choice of central subring $R$ of $\Lambda$. Thus we denote it by $\depth M$. The following results will be quite useful.

\vskip.5em{\bf Proposition \XBA\ }[GN2;3.5] {\it
For any $M\in\mod\Lambda$, we have

$\depth M\le\dim M\le\id{}_\Lambda M=\sup\{i\in\zzz\ |\ \ext^i_\Lambda(\Lambda/J_\Lambda,M)\neq0\}$.

In particular, we have $\depth\Lambda\le\dim \Lambda\le\gl\Lambda$.}

\vskip.5em{\bf Propostion \XBB\ }{\it
For $\Lambda$ as above, we have $\gl\Lambda=\sup\{\pd{}_\Lambda M\ |\ M\in\flmod\Lambda\}$.}

\vskip.5em{\sc Proof }
We put $n:=\pd{}_\Lambda(\Lambda/J_\Lambda)$, and can clearly assume that $n<\infty$. We will show that $\pd{}_\Lambda M\le n$ holds for any $M\in\mod\Lambda$ by using induction on $\dim M$. If $\dim M=0$, then $M\in\flmod\Lambda$ holds, so we have $\pd{}_\Lambda M\le n$. Now assume that $\pd{}_\Lambda M\le n$ holds for any $M\in\mod\Lambda$ with $\dim M<m$, where $m>0$. Take any $M\in\mod\Lambda$ with $\dim M=m$. There is then an exact sequence
$0\to L\to M\to N\to 0$
with $L\in\flmod\Lambda$ and $\depth N\ge1$. Since we have $\pd{}_\Lambda L\le n$ by assumption, we only have to show $\pd{}_\Lambda N\le n$. Take an $N$-regular element $r\in R$ and consider the exact sequence
$0\to N\stackrel{r}{\to}N\to N/rN\to0$.
Since $\dim (N/rN)<m$ holds, we have $\pd{}_\Lambda(N/rN)\le n$ by the induction assumption. Applying Nakayama's lemma to the exact sequence
$\ext^n_\Lambda(N,-)\stackrel{r}{\to}\ext^n_\Lambda(N,-)\to0$,
we have $\pd{}_\Lambda N<n$. We see that $\pd{}_\Lambda M\le n$.\rule{5pt}{10pt}

\vskip.5em
We call $M\in\mod\Lambda$ a {\it Cohen-Macaulay} ({\it CM} for short) {\it $\Lambda$-module of dimension $n$} if $\depth M=\dim M=n$. We simply call a CM $\Lambda$-module of dimension $d(=\dim R)$ a (maximal) {\it CM $\Lambda$-module}. We denote by $\cm\Lambda$ the category of CM $\Lambda$-modules. We call $\Lambda$ an {\it $R$-order}, or just an {\it order} if $\Lambda\in\cm\Lambda$. If $R$ has a canonical module $\omega_R$ (e.g. $R$ is Gorenstein and $\omega_R=R$), then $\depth M=d-\sup\{i\ge0\ |\ \ext^i_R(M,\omega_R)\neq0\}$ and $\dim M=d-\inf\{i\ge0\ |\ \ext^i_R(M,\omega_R)\neq0\}$ [BH;3.5.11]. Thus $M$ is CM of dimension $n$ if and only if $\ext^{d-i}_R(M,\omega_R)=0$ for any $i\neq n$.

As \XBA\ suggests, orders with $\gl\Lambda=d$ are very special. For example, we have the following Auslander-Buchsbaum type equality.

\vskip.5em{\bf Proposition \XBC\ }{\it
Let $\Lambda$ be an order with $\gl\Lambda=d$. For any $M\in\mod\Lambda$, we have $\pd{}_\Lambda M+\depth M=d$.}
%Let $\Lambda$ be an order and $M\in\mod\Lambda$. If $\pd{}_\Lambda M\le d$, then $\pd{}_\Lambda M+\depth M=d$.}

\vskip.5em{\sc Proof }
Put $n:=\pd{}_\Lambda M$ and $t:=\depth M$.
We have a projective resolution $0\to P_n\stackrel{f_n}{\to}\cdots\stackrel{f_1}{\to}P_0\stackrel{f_0}{\to}M\to0$.
Put $M_i:=\Im f_i$ for $0\le i\le n$, then
we have an exact sequence $0\to M_{i+1}\stackrel{}{\to}P_i\to M_i\to0$ with $P_i\in\cm\Lambda$. Applying $\hom_{\Lambda}(\Lambda/J_\Lambda,-)$, we have $\depth M_i\ge\depth M_{i+1}-1$.
Thus we have $t=\depth M_0\ge \depth M_n-n=d-n$.
%an exact sequence $\ext^j_\Lambda(\Lambda/J_\Lambda,P_i)\stackrel{}{\to}\ext^j_\Lambda(\Lambda/J_\Lambda,M_i)\to\ext^{j+1}_\Lambda(\Lambda/J_\Lambda,M_{i+1})$ for any $j$.
%with $P_i\in\cm\Lambda$. Applying $\hom_{\Lambda}(\Lambda/J_\Lambda,-)$, we have $\ext^i_\Lambda(\Lambda/J_\Lambda,M)=0$ for any $0\le i\le d-n$.

On the other hand, take an $M$-regular sequence $(x_1,\cdots,x_t)$.
Put $N_i:=M/(x_1,\cdots,x_i)M$ for $0\le i\le t$, then
we have an exact sequence $0\to N_i\stackrel{x_{i+1}}{\to}
N_i\to N_{i+1}\to0$. Applying $\hom_{\Lambda}(\Lambda/J_\Lambda,-)$, we have an
exact sequence $\ext^j_\Lambda(\Lambda/J_\Lambda,N_i)\stackrel{x_{i+1}}{\to}
\ext^j_\Lambda(\Lambda/J_\Lambda,N_i)\to\ext^{j+1}_\Lambda(\Lambda/J_\Lambda,N_{i+1})
\to\ext^{j+1}_\Lambda(\Lambda/J_\Lambda,N_{i+1})$ for any $j$.
Using Nakayama's Lemma, we have $\pd{}_\Lambda N_{i+1}=\pd{}_\Lambda N_i+1$.
Consequently, we have $n=\pd{}_\Lambda N_0=\pd{}_\Lambda N_t-t\le d-t$.\rule{5pt}{10pt}

%We use induction on $\pd{}_\Lambda M$. If $M$ is projective, then $\depth M=d$ holds. Assume that the equation holds if $\pd{}_\Lambda M<i$. Take $N\in\mod\Lambda$ with $\pd{}_\Lambda N=i$. Let $0\to M\to P\to N\to0$ be exact, with $P$ projective. Then $\pd{}_\Lambda M=i-1$ holds. Since $\depth M=d-i+1$ holds by the induction assumption, we obtain $\depth N=d-i$.\rule{5pt}{10pt}

%(2) If $R$ is Gorenstein, then $\sup\{i\ge0\ |\ \ext^i_\Lambda(M,\Lambda^*)\neq0\}+\depth M=d$.}

%(2) Let $\cdots\to P_0\to M\to0$ be a projective resolution. Since $\Lambda$ is an order, $\ext^i_R(P_i,R)=0$ for any $i>0$. Applying $\hom_\Lambda(-,\Lambda^*)\simeq\hom_R(-,R)$, we see that $\ext^i_\Lambda(M,\Lambda^*)=\ext^i_R(M,R)$. Thus the assertion follows from the remark above.\rule{5pt}{10pt}

\vskip.5em
Let $R$ be an arbitrary commutative noetherian ring. We call $M\in\mod\Lambda$ a {\it CM $\Lambda$-module} if $M_{\dn{p}}\in\cm\Lambda_{\dn{p}}$ for any $\dn{p}\in\Spec R$, and we denote by $\cm\Lambda$ the category of CM $\Lambda$-modules. We call $\Lambda$ an {\it $R$-order}, or just an {\it order} if $\Lambda\in\cm\Lambda$.

We call $\Lambda$ a {\it symmetric} $R$-algebra if $\hom_R(\Lambda,R)$ is isomorphic to $\Lambda$ as a $(\Lambda,\Lambda)$-module.

\vskip.5em
Now assume that $R$ is a normal domain. We want to investigate reflexive equivalence, especially in connection with symmetric algebras. Recall that $(-)^*$ is the functor $\hom_R(-,R):\mod R\to\mod R$.
We call $M\in\mod\Lambda$ a {\it reflexive} $\Lambda$-module if the natural map $M\to M^{**}$ (not $\hom_{\Lambda^{\op}}(\hom_\Lambda(M,\Lambda),\Lambda)$!) is an isomorphism. It is well known that $M\in\mod\Lambda$ is reflexive if and only if $M$ satisfies Serre's S$_2$ condition $\depth{}_{R_{\dn{p}}}M_{\dn{p}}\ge\min\{2,\height\dn{p}\}$ for any $\dn{p}\in\Spec R$ [EG;0.B,3.6]. We denote by $\ref\Lambda$ the category of reflexive $\Lambda$-modules. Using the S$_2$ condition, one can easily check that $\ref\Lambda$ is closed under kernels and extensions. We have a functor $(-)^{**}=\hom_R(\hom_R(-,R),R):\mod\Lambda\to\ref\Lambda$ (e.g. \XBD(1) below). This gives a left adjoint of the inclusion functor $\ref\Lambda\to\mod\Lambda$. 

%Let $K$ be a quotient field of $R$ and let $X$ be in $\mod\Lambda$. It is well-known that $X^{**}=\bigcap_{\dn{p}}X_{\dn{p}}\ (\subset K\otimes_RX)$ holds for any $X\in\ref\Lambda$, where $\dn{p}$ runs over all height one prime ideals of $R$ and we denote by $X_{\dn{p}}$ the image of the induced map $X_{\dn{p}}=R_{\dn{p}}\otimes_RX\to K\otimes_RX$.

We say that two $R$-algebras $\Lambda$ and $\Gamma$ are {\it reflexive equivalent} if the additive categories $\ref\Lambda$ and $\ref\Gamma$ are equivalent. We call $M\in\ref\Lambda$ a {\it height one generator} (resp. {\it progenerator}, {\it projective}) if $M_{\dn{p}}$ is a generator (resp. a progenerator, projective) over $\Lambda_{\dn{p}}$ for any height one prime ideal $\dn{p}$ of $R$.
The proposition below shows that many algebras which are not Morita equivalent may be reflexive equivalent. Part (2) is used in [RV1].

\vskip.5em{\bf Proposition \XBD\ }{\it
(1) $\hom_\Lambda(X,Y)\in\ref R$ for any $X\in\mod\Lambda$ and $Y\in\ref\Lambda$.

(2) Let $\Lambda$ and $\Gamma$ be $R$-algebras which are reflexive $R$-modules.

\strut\kern.5em
(i) For any height one progenerator $M\in\ref\Lambda$, we have an equivalence $\fff:=\hom_\Lambda(M,-):\ref\Lambda\to\ref\endm_\Lambda(M)$.

\strut\kern.5em
(ii) Let $\fff:\ref\Lambda\to\ref\Gamma$ be a categorical equivalence. Then there exists a height one progenerator $M\in\ref\Lambda$ such that $\Gamma\simeq\endm_\Lambda(M)$ and $\fff\simeq\hom_\Lambda(M,-)$.

(3) If $\Lambda$ is a symmetric $R$-algebra, then so is $\endm_\Lambda(M)$ for any height one projective $M\in\ref\Lambda$. Thus symmetric algebras are closed under reflexive equivalences.}

\vskip.5em{\sc Proof }
(1) Take an exact sequence $\Lambda^n\to\Lambda^m\to X\to0$ in $\mod\Lambda$. Applying $\hom_\Lambda(-,Y)$, we obtain an exact sequence $0\to\hom_\Lambda(X,Y)\to Y^m\to Y^n$. Since $\ref R$ is closed under kernels, we have $\hom_\Lambda(X,Y)\in\ref R$.

(2)(i) Let $\Gamma:=\endm_\Lambda(M)$ and consider the functor $\ggg:=\hom_\Gamma(\fff(\Lambda),-):\ref\Gamma\to\ref\Lambda$. We have homomorphisms $f:(\fff(\Lambda)\otimes_\Lambda M)^{**}=(\hom_\Lambda(M,\Lambda)\otimes_\Lambda M)^{**}\to\endm_\Lambda(M)^{**}=\Gamma^{**}=\Gamma$ and $g:(M\otimes_\Gamma\fff(\Lambda))^{**}=(M\otimes_\Gamma\hom_\Lambda(M,\Lambda))^{**}\to\Lambda^{**}=\Lambda$. Since $M$ is a height one progenerator, $f_{\dn{p}}$ and $g_{\dn{p}}$ are isomorphisms for any height one prime ideal $\dn{p}$ of $R$. Since $f$ and $g$ are homomorphisms between reflexive $R$-modules, they are isomorphisms.
%, we obtain{\small\begin{eqnarray*}&(\fff(\Lambda)\otimes_\Lambda M)^{**}=\bigcap_{\dn{p}}(\fff(\Lambda)\otimes_\Lambda M)_{\dn{p}}=\bigcap_{\dn{p}}(\hom_{\Lambda_{\dn{p}}}(M_{\dn{p}},\Lambda_{\dn{p}})\otimes_{\Lambda_{\dn{p}}}M_{\dn{p}})=\bigcap_{\dn{p}}\endm_{\Lambda_{\dn{p}}}(M_{\dn{p}})=\bigcap_{\dn{p}}\Gamma_{\dn{p}}=\Gamma&\\ &(M\otimes_\Gamma\fff(\Lambda))^{**}=\bigcap_{\dn{p}}(M\otimes_\Gamma\fff(\Lambda))_{\dn{p}}=\bigcap_{\dn{p}}(M_{\dn{p}}\otimes_{\Gamma_{\dn{p}}}\hom_{\Lambda_{\dn{p}}}(M_{\dn{p}},\Lambda_{\dn{p}}))\stackrel{{\rm (a)}}{=}\bigcap_{\dn{p}}\Lambda_{\dn{p}}=\Lambda& \end{eqnarray*}} where $\dn{p}$ runs over all height one prime ideals of $R$ and (a) follows by Morita equivalence.
Using adjointness properties, we obtain
{\small\begin{eqnarray*}
&\fff\circ\ggg=\hom_\Lambda(M,\hom_\Gamma(\fff(\Lambda),-))=\hom_\Gamma(\fff(\Lambda)\otimes_\Lambda M,-)=\hom_\Gamma((\fff(\Lambda)\otimes_\Lambda M)^{**},-)=1&\\
&\ggg\circ\fff=\hom_\Gamma(\fff(\Lambda),\hom_\Lambda(M,-))=\hom_\Lambda(M\otimes_\Gamma\fff(\Lambda),-)=\hom_\Lambda((M\otimes_\Gamma\fff(\Lambda))^{**},-)=1.&
\end{eqnarray*}}
\vskip-1em
Hence $\fff:\ref\Lambda\to\ref\Gamma$ is an equivalence.

(ii) Left to the reader.

(3) See [I4;5.4.3(1)], for example. For completeness, we give a proof here. By the same argument as in the proof of (2), we have an isomorphism $f:(\hom_\Lambda(M,\Lambda)\otimes_\Lambda M)^{**}\to\Gamma$. Thus we have isomorphisms
{\small\[\Gamma^*\stackrel{f^*}{\simeq}(\hom_\Lambda(M,\Lambda)\otimes_\Lambda M)^*\stackrel{}{\simeq}(\hom_\Lambda(M,\Lambda^*)\otimes_\Lambda M)^*\simeq(M^*\otimes_\Lambda M)^*\simeq\hom_\Lambda(M,M^{**})\simeq\Gamma\]}
of $(\Gamma,\Gamma)$-modules.\rule{5pt}{10pt}

\vskip.5em
We want to recall some results on derived categories which will be useful in the next section. We start with basic notation and definitions. For the rest of the section $R$ is a commutative noetherian ring and $\Lambda$ is a module-finite $R$-algebra.

For an additive category $\aa$, we denote by $\cc(\aa)$ the category of complexes over $\aa$, by $\kk(\aa)$ the {\it homotopy category} of $\aa$, and by $\dd(\aa)$ the {\it derived category} of $\aa$ provided $\aa$ is abelian [Hap][Har1]. For $*=+,-$ or $\b$, we denote by $\cc^*(\aa)$ (resp. $\kk^*(\aa)$, $\dd^*(\aa)$) the full subcategory of $\cc(\aa)$ (resp. $\kk(\aa)$, $\dd(\aa)$) consisting of bounded below, bounded above or bounded complexes respectively. Moreover, for a full subcategory $\bb$ of $\aa$, we denote by $\cc^*_{\bb}(\aa)$ (resp. $\kk^*_{\bb}(\aa)$, $\dd^*_{\bb}(\aa)$) the full subcategory of $\cc^*(\aa)$ (resp. $\kk^*(\aa)$, $\dd^*(\aa)$) consisting of all objects $X$ such that the $i$-th homology $H^i(X)$ belongs to $\bb$ for any $i$.

%For $\aa=\Mod\Lambda,\ \mod\Lambda$ or $\flmod\Lambda$, we denote by $\dd^{\b}(\aa)$ (resp. $\dd^-(\aa)$) the bounded (resp. bounded above) derived category of $\aa$. For $\pp=\Pr\Lambda$ or $\pr\Lambda$, we denote by $\kk^{\b}(\pp)$ (resp. $\kk^-(\pp)$) the homotopy category of bounded (resp. bounded below) complexes on $\pp$. 
%We denote by $\dd^{\b}_{\flmod\Lambda}(\mod\Lambda)$ (resp. $\kk^{\b}_{\flmod\Lambda}(\mod\Lambda)$) the full subcategory of $\dd^{\b}(\mod\Lambda)$ (resp. $\kk^{\b}(\mod\Lambda)$) consisting of all objects $X$ such that $H^i(X)\in\flmod\Lambda$ for any $i$. 

%For $*=\b$ or $-$, we have the following commutative diagram of full and faithful functors {\small \[\begin{diag} \kk^*(\pr\Lambda)&\ \subset\ &\kk^*(\Pr\Lambda)\\ \cap&&\cap\\ \dd^*(\mod\Lambda)&\ \subset\ &\dd^*(\Mod\Lambda),\end{diag}\]}whose vertical functors are equivalences if $*=-$.

We have natural equivalences $\kk^-(\Pr\Lambda)\stackrel{\sim}{\to}\dd^-(\Mod\Lambda)$, $\kk^-(\pr\Lambda)\stackrel{\sim}{\to}\dd^-(\mod\Lambda)$ and $\dd^{\b}(\mod\Lambda)\stackrel{\sim}{\to}\dd^{\b}_{\mod\Lambda}(\Mod\Lambda)$. The next two results were pointed out to us by Rickard [Ri2].

\vskip.5em{\bf Lemma \XBE\ }{\it
For any $X\in\cc^{\b}_{\flmod\Lambda}(\mod\Lambda)$, there exists a quasi-isomorphism $X\to Y$ with $Y\in\cc^{\b}(\flmod\Lambda)$. Thus we have equivalences $\dd^{\b}(\flmod\Lambda)\stackrel{\sim}{\to}\dd^{\b}_{\flmod\Lambda}(\mod\Lambda)\stackrel{\sim}{\to}\dd^{\b}_{\flmod\Lambda}(\Mod\Lambda)$.}

\vskip.5em{\sc Proof }
It follows from [Ve;III.2] that we only have to check the condition:

($E_2$)$^{\op}$ Let $X\in\mod\Lambda$ and $Y\in\flmod\Lambda$ a submodule of $X$. Then there exists a submodule $Z$ of $X$ such that $Y\cap Z=0$ and $X/Z\in\flmod\Lambda$.

To see this, let $I$ be an arbitrary ideal of $R$. Since $R$ is noetherian, there exists $c>0$ such that $I^nX\cap Y=I^{n-c}(I^cX\cap Y)$ holds for any $n>c$ by the Artin-Rees Lemma [Ma;8.5]. Applying this to $I:=\ann{}_RY$, we have $I^{c+1}X\cap Y=I(I^cX\cap Y)=0$. Since $R/I$ is an artin ring, $Z:=I^{c+1}X$ satisfies the desired conditions.\rule{5pt}{10pt}

\vskip.5em
In \XCA\ we shall use the results below due to Rickard [Ri2]. We consider a descending chain $I_1\supset I_2\supset\cdots$ of ideals of $\Lambda$ such that $\Lambda^{(i)}:=\Lambda/I_i\in\flmod R$ for any $i\ge0$. Put $\widehat{\Lambda}:=\plim_{i\ge0}\ \Lambda^{(i)}$. Denote as usual by $\Lotimes$ (resp. $\Rhom$) the left (resp. right) derived functor of $\otimes$ (resp. $\hom$).

\vskip.5em{\bf Proposition \XBF\ }{\it
With the above notation, put $P^{(i)}:=\Lambda^{(i)}\Lotimes_\Lambda P$ and $\widehat{P}:=\widehat{\Lambda}\Lotimes_\Lambda P$ for $P\in \kk^{\b}(\pr\Lambda)$.

%(1) Assume $\Lambda=\plim_{i\ge0}\ \Lambda^{(i)}$. Then $P=\plim_{i\ge0}\ P^{(i)}$ in $\dd^-(\mod\Lambda)$, 

(1) $\plim_{i\ge0}\hom_{\dd(\Mod\Lambda)}(X,P^{(i)})=\hom_{\dd(\Mod\Lambda)}(X,\widehat{P})$ for any $X\in \dd^-(\mod\Lambda)$.

(2) If $I_i=J_\Lambda^i$, then $\ilim_{i\ge0}\hom_{\dd(\Mod\Lambda)}(P^{(i)},X)=\hom_{\dd(\Mod\Lambda)}(P,X)$ for any $X\in\dd^{\b}(\flmod\Lambda)$.}

\vskip.5em{\sc Proof }
Let $P$ be a bounded complex $\cdots\to P^{-1}\to P^0\to P^1\to\cdots$ in $\pr\Lambda$. Then $P^{(i)}$ is given by the bounded complex $\cdots\to(P^{-1})^{(i)}\to(P^{0})^{(i)}\to(P^{1})^{(i)}\to\cdots$ in $\pr\Lambda^{(i)}$.

(1) We can assume that $X$ is given by the bounded above complex $\cdots\to Q^{-1}\to Q^0\to Q^1\to\cdots$ in $\pr\Lambda$.
Since $\plim\hom_\Lambda(Q^s,(P^{t})^{(i)})=\hom_\Lambda(Q^s,\widehat{P}^t)$ holds, we have an isomorphism of complexes
{\small\[\begin{diag}
\cdots&\to&\plim\prod_{-s+t=-1}\hom_\Lambda(Q^s,(P^{t})^{(i)})&\stackrel{}{\to}&\plim\prod_{-s+t=0}\hom_\Lambda(Q^s,(P^{t})^{(i)})&\to&\plim\prod_{-s+t=1}\hom_\Lambda(Q^s,(P^{t})^{(i)})&\to&\cdots\\
&&|\wr&&|\wr&&|\wr\\
\cdots&\to&\prod_{-s+t=-1}\hom_\Lambda(Q^s,\widehat{P}^t)&\stackrel{}{\to}&\prod_{-s+t=0}\hom_\Lambda(Q^s,\widehat{P}^t)&\to&\prod_{-s+t=1}\hom_\Lambda(Q^s,\widehat{P}^t)&\to&\cdots.
\end{diag}\]}
Since each inverse system $(\prod_{-s+t=n}\hom_\Lambda(Q^s,(P^{t})^{(i)}))_i$ consists of finite length $R$-modules, it satisfies the Mittag-Leffler condition (e.g. [Har2;II.9]). Thus the $0$-th homology of the upper sequence in the diagram above is isomorphic to $\plim\hom_{\dd(\Mod\Lambda)}(X,P^{(i)})$, and we obtain the desired isomorphism.

(2) We will show that the natural map $\ilim\hom_{\dd(\Mod\Lambda)}(P^{(i)},X)\to\hom_{\dd(\Mod\Lambda)}(P,X)$ is bijective. Fix any $f\in\hom_{\dd(\Mod\Lambda)}(P,X)$. By \XBE, we can write $f=gs^{-1}$ for $g\in\hom_{\kk^{\b}(\mod\Lambda)}(P,X^\prime)$ and a quasi-isomorphism $s\in\hom_{\kk^{\b}(\flmod\Lambda)}(X,X^\prime)$. It follows from $X^\prime\in\kk^{\b}(\flmod\Lambda)$ that $g$ factors through some $P^{(i)}$, and hence $f$ also factors through $P^{(i)}$. Thus the above map is surjective.

We now show injectivity. Let $p:P\to P^{(i)}$ be the natural map. Assume that $f\in\hom_{\dd(\Mod\Lambda)}(P^{(i)},X)$ satisfies $pf=0$. Again we write $f=gs^{-1}$ for $g\in\hom_{\kk^{\b}(\mod\Lambda)}(P^{(i)},X^\prime)$ and a quasi-isomorphism $s\in\hom_{\kk^{\b}(\flmod\Lambda)}(X,X^\prime)$. Since $pg$ is null-homotopic, we can take a homotopy $a:P\to X^\prime[-1]$. It follows from $X^\prime\in\kk^{\b}(\flmod\Lambda)$ that $a$ factors through $P^{(j)}$ for some sufficiently large $j$. Then the composition of $P^{(j)}\to P^{(i)}$ and $g$ is null-homotopic.\rule{5pt}{10pt}

%Let $X$ be a bounded complex $\cdots\to X^{-1}\to X^0\to X^1\to\cdots$ in $\flmod\Lambda$. For sufficiently large $i$, $X$ is a bounded complex in $\flmod\Lambda^{(i)}$. For such $i$, $\hom_{\dd(\Mod\Lambda)}(P^{(i)},X)$ is the $0$-th homology of the complex $\Rhom_\Lambda(P^{(i)},X):$ \[\cdots\to\prod_{-s+t=-1}\hom(P^{s(i)},X^t)\to\prod_{-s+t=0}\hom(P^{s(i)},X^t)\to\prod_{-s+t=1}\hom(P^{s(i)},X^t)\to\cdots.\] Since any $X^t$ has a finite length, $\ilim\hom_{\dd(\Mod\Lambda)}(P^{s(i)},X^t)=\hom_{\dd(\Mod\Lambda)}(P^s,X^t)$ holds. Thus we have an isomorphism of complexes {\small\[\begin{diag} \cdots&\to&\ilim\prod_{-s+t=-1}\hom(P^{s(i)},X^t)&\to&\ilim\prod_{-s+t=0}\hom(P^{s(i)},X^t)&\to&\ilim\prod_{-s+t=1}\hom(P^{s(i)},X^t)&\to&\cdots\\ &&|\wr&&|\wr&&|\wr\\ \cdots&\to&\prod_{-s+t=-1}\hom(P^{s},X^t)&\to&\prod_{-s+t=0}\hom(P^{s},X^t)&\to&\prod_{-s+t=1}\hom(P^{s},X^t)&\to&\cdots.\end{diag}\]} Since $\ilim$ is an exact functor, we obtain the desired isomorphism by taking $0$-th homology of complexes above.\rule{5pt}{10pt}

\vskip.5em
The concept of {\it derived equivalence} is central for our work. Recall that $T\in\kk^{\b}(\pr\Lambda)$ is a {\it tilting complex} if $\hom_{\dd(\Mod\Lambda)}(T,T[i])=0$ for any $i\neq0$ and $T$ generates $\kk^{\b}(\pr\Lambda)$. If a $\Lambda$-module $T$ is a tilting complex, it is called a {\it tilting module}. We will mostly deal with tilting modules $T$ of projective dimension at most one. They satisfy the conditions (i) $\ext^1_\Lambda(T,T)=0$ and (ii) there exists an exact sequence $0\to\Lambda\to T_0\to T_1\to0$ with $T_i\in\add T$.

Rickard proved in [Ri1] that the following conditions (1)--(4) are equivalent.

(1) (resp. (2), (3)) $\kk^{\b}(\pr\Lambda)$ and $\kk^{\b}(\pr\Gamma)$ (resp. $\dd^{\b}(\mod\Lambda)$ and $\dd^{\b}(\mod\Gamma)$, $\kk^-(\Pr\Lambda)$ and $\kk^-(\Pr\Gamma)$) are triangle equivalent.

(4) There exists a tilting complex $T\in\kk^{\b}(\pr\Lambda)$ such that $\Gamma\simeq\endm_{\dd(\Mod\Lambda)}(T)$.

If these conditions are satisfied, we call $\Lambda$ and $\Gamma$ {\it derived equivalent}. We call $T$ in (4) above a {\it two-sided tilting complex} if $T\in\dd^{\b}(\mod\Lambda\otimes_{\zzz}\Gamma^{\op})$.

We have the following relationship with localizations.

\vskip.5em{\bf Lemma \XBG\ }{\it
Let $T\in\kk^{\b}(\pr\Lambda)$ be a tilting complex with $\Gamma:=\endm_{\dd(\Mod\Lambda)}(T)$. For any $\dn{p}\in\Spec R$, we have a tilting complex $T_{\dn{p}}\in\kk^{\b}(\pr\Lambda_{\dn{p}})$ with $\endm_{\dd(\Mod\Lambda_{\dn{p}})}(T_{\dn{p}})=\Gamma_{\dn{p}}$.}

\vskip.5em
Recall that a $\Lambda$-module $T$ is said to be a {\it partial tilting module} if $\pd{}_\Lambda T\le1$ and $\ext^1_\Lambda(T,T)=0$. The following proposition is a generalization of Bongartz's result [Bon1] for finite dimensional algebras. We sometimes call $X$ a {\it Bongartz complement} of $T$.

\vskip.5em{\bf Lemma \XBH\ }{\it
For any partial tilting $\Lambda$-module $T$, there exists $X\in\mod\Lambda$ such that $T\oplus X$ is a tilting $\Lambda$-module. Moreover, if $R$ is normal and $T\oplus\Lambda\in\ref\Lambda$, then $X\in\ref\Lambda$.}

\vskip.5em{\sc Proof }
Let $P\stackrel{f}{\to}\ext^1_\Lambda(T,\Lambda)\to0$ be exact, with $P$ projective in $\mod\endm_\Lambda(T)$. We can write $P=\hom_\Lambda(T,T^\prime)$ for $T'\in\add T$. It follows from Yoneda's lemma on $\add T$ that $f$ is given by $\sigma\in\ext^1_\Lambda(T',\Lambda)$. Take an exact sequence $0\to\Lambda\to X\to T^\prime\to0$ corresponding to $\sigma$. Then $\hom_\Lambda(T,T^\prime)\stackrel{f=(\bullet\sigma)}{\to}\ext^1_\Lambda(T,\Lambda)\to0$ is exact and $\pd{}_\Lambda(T\oplus X)\le1$ holds. Applying $\hom_\Lambda(T,-)$, we see that $\ext^1_\Lambda(T,X)=0$. Applying $\hom_\Lambda(-,T\oplus X)$, we get $\ext^1_\Lambda(T\oplus X,T\oplus X)=0$. Thus $T\oplus X$ is a tilting $\Lambda$-module. If $R$ is normal, then $\ref\Lambda$ is closed under extensions. Thus the second assertion follows.\rule{5pt}{10pt}

\vskip.5em
The following easy lemma is useful (e.g. [HU2;1.2], [Ye2;2.3]).

\vskip.5em{\bf Lemma \XBI\ }{\it
Let $T$ be a tilting $\Lambda$-module with a minimal projective resolution $0\to P_1\to P_0\to T\to0$. Then $\add(P_0\oplus P_1)=\add\Lambda$ and $\add P_0\cap\add P_1=0$. In particular, if $\Lambda$ is Morita equivalent to a local ring, then $T$ is projective.}

\vskip.5em
We will use the following canonical isomorphisms (e.g. [F;pp.153]), where we denote by $\Flat\Gamma$ the category of flat $\Gamma$-modules.

\vskip.5em{\bf Lemma \XBJ }{\it
(1) $\Rhom_\Lambda(Y\Lotimes_\Gamma X,Z)\simeq\Rhom_\Gamma(X,\Rhom_\Lambda(Y,Z))$ for any $X\in\dd^-(\Mod\Gamma)$, $Y\in\dd^-(\Mod\Lambda\otimes_R\Gamma^{\op})$ and $Z\in \dd^+(\Mod\Lambda)$.

(2) $\Rhom_\Lambda(X,Y)\Lotimes_\Gamma Z\simeq\Rhom_\Lambda(X,Y\Lotimes_\Gamma Z)$ for any $X\in\dd^-(\mod\Lambda)$, $Y\in \dd^{\b}(\Mod\Lambda\otimes_R\Gamma^{\op})$ and $Z\in\dd^{\b}(\Mod\Gamma)$, provided $X\in\kk^{\b}(\pr\Lambda)$ or $Z\in\kk^{\b}(\Flat\Gamma)$.}

\vskip1.5em{\bf\XC. Calabi-Yau algebras and symmetric orders}

Throughout this section, let $R$ be a commutative noetherian ring with $\dim R=d$ and $\Lambda$ a module-finite $R$-algebra. We denote by $E(X)$ the injective hull of $X\in\mod R$ and put $E:=\bigoplus_{\dn{p}\in\Max R}E(R/\dn{p})$.
%We denote by $0\to R\to I^0\to\cdots\to I^d\to0$ the minimal injective resolution of $R$. 
Then $E$ is an injective $R$-module, and we have a duality
{\small\[D:=\hom_R(-,E):\flmod R\to\flmod R\ \ \ (\mbox{resp. }\flmod\Lambda\to\flmod\Lambda^{\op},\ \dd^{\b}(\flmod R)\to\dd^{\b}(\flmod R),\ \dd^{\b}(\flmod\Lambda)\to\dd^{\b}(\flmod\Lambda^{\op}))\]}
called {\it Matlis duality}, such that $D\circ D$ is isomorphic to the identity functor [BH;3.1.3,3.2.13]. For example, if $R$ is a polynomial or power series ring over a field $K$, then $D$ is isomorphic to $\hom_K(-,K)$. Obviously, $X\in\Mod\Lambda$ belongs to $\flmod\Lambda$ if and only if it belongs to $\flmod R$ as an $R$-module. For $X,Y\in\dd^{\b}(\mod\Lambda)$, $\hom_{\dd(\Mod\Lambda)}(X,Y)$ belongs to $\flmod R$ if $X$ or $Y$ is in $\dd^{\b}(\flmod\Lambda)$.

For an integer $n$, we call $\Lambda$ {\it $n$-Calabi-Yau} ({\it $n$-CY} for short) if there exists a functorial isomorphism
\[\hom_{\dd(\Mod\Lambda)}(X,Y[n])\simeq D\hom_{\dd(\Mod\Lambda)}(Y,X)\ \ \ \ \ \ \ (*)\]
for any $X,Y\in \dd^{\b}(\flmod\Lambda)$. Similarly, we call $\Lambda$ {\it $n$-Calabi-Yau$^-$} ({\it $n$-CY$^-$} for short) if there exists a functorial isomorphism $(*)$ for any $X\in \dd^{\b}(\flmod\Lambda)$ and $Y\in\kk^{\b}(\pr\Lambda)$.

In this section we first give some basic results for $n$-CY and $n$-CY$^-$ algebras in \XCA. The main result is a characterization of these algebras in terms of symmetric orders in \XCB. Obviously, $n$-CY (resp. $n$-CY$^-$) algebras are closed under Morita equivalences. A finite product $\prod_i\Lambda_i$ of algebras is $n$-CY (resp. $n$-CY$^-$) if and only if so is each $\Lambda_i$. Part (1) and (7) in the following theorem are due to Rickard [Ri2].

\vskip.5em{\bf Theorem \XCA\ }{\it
(1) $n$-CY (resp. $n$-CY$^-$) algebras are closed under derived equivalences.

(2) $\Lambda$ is $n$-CY if and only if so is $\Lambda^{\op}$ (cf. \XCD(1)).

(3)(localization) $\Lambda$ is $n$-CY (resp. $n$-CY$^-$) if and only if so is $\Lambda_{\dn{p}}$ for any $\dn{p}\in\Max R\cap\Supp{}_R\Lambda$.

(4)(completion) If $R$ is local, then $\Lambda$ is $n$-CY if and only if so is $\widehat{\Lambda}$ (cf. \XCD(2)).

(5) Any $n$-CY algebra $\Lambda$ satisfies $\gl\Lambda=n$. 

(6) Any $n$-CY$^-$ algebra $\Lambda$ satisfies $\depth{}_{R_{\dn{p}}}\Lambda_{\dn{p}}=\dim \Lambda=\id{}_\Lambda\Lambda=n$ for any $\dn{p}\in\Max R\cap\Supp{}_R\Lambda$.

(7) $\Lambda$ is $n$-CY if and only if it is $n$-CY$^-$ and $\gl\Lambda<\infty$.

(8) If $\Lambda$ is $n$-CY, then there exists a functorial isomorphism
$\hom_{\dd(\Mod\Lambda)}(X,Y[n])\simeq D\hom_{\dd(\Mod\Lambda)}(Y,X)$ for any $X\in \dd^{\b}(\flmod\Lambda)$ and $Y\in \dd^{\b}(\mod\Lambda)$.}

\vskip.5em{\sc Proof }
(1) By [Re;6.3], $\kk^{\b}(\pr\Lambda)$ consists of the compact objects $X$ of $\kk^-(\Pr\Lambda)$, i.e. the functor $\hom_{\dd(\Mod\Lambda)}(X,-)$ on $\kk^-(\Pr\Lambda)$ commutes with arbitrary direct sums. On the other hand, $\dd^{\b}(\flmod\Lambda)$ consists of all objects $X\in\kk^-(\Pr\Lambda)$ such that $\bigoplus_{i\in\zzz}\hom_{\dd(\Mod\Lambda)}(Y,X[i])$ has finite length for any $Y\in\kk^{\b}(\pr\Lambda)$ (see \XBE). Thus any triangle equivalence $\kk^-(\Pr\Lambda)\to\kk^-(\Pr\Gamma)$ induces triangle equivalences $\kk^{\b}(\pr\Lambda)\to\kk^{\b}(\pr\Gamma)$ and $\dd^{\b}(\flmod\Lambda)\to \dd^{\b}(\flmod\Gamma)$. Thus the $n$-CY and $n$-CY$^-$ properties are preserved by derived equivalence.

(2) Immediate from Matlis duality.

(3)(4) Put $S:=R_{\dn{p}}$ for (3) and $S:=\widehat{R}$ for (4). Then $S$ is a flat $R$-module. We have a functor $(-\otimes_RS):\dd^{\b}(\Mod\Lambda)\to\dd^{\b}(\Mod\Lambda\otimes_RS)$ with isomorphisms
{\small\[\Rhom_{\Lambda\otimes_RS}(X\otimes_RS,Y\otimes_RS)\simeq\Rhom_\Lambda(X,Y\otimes_RS)\simeq\Rhom_\Lambda(X,Y)\otimes_RS\]}
for any $X,Y\in\dd^{\b}(\mod\Lambda)$ by \XBJ. Applying $H^0$, we have a functorial isomorphism
{\small\[\hom_{\dd(\Mod\Lambda\otimes_RS)}(X\otimes_RS,Y\otimes_RS)\simeq\hom_{\dd(\Mod\Lambda)}(X,Y)\otimes_RS.\]}
If one of $X,Y\in\dd^{\b}(\mod\Lambda)$ is contained in $\dd^{\b}(\flmod\Lambda)$, we have functorial isomorphisms
{\small\begin{eqnarray*}
\hom_{\dd(\Mod\Lambda)}(X,Y)\simeq\bigoplus_{\dn{p}}\hom_{\dd(\Mod\Lambda)}(X,Y)_{\dn{p}}\simeq\bigoplus_{\dn{p}}\hom_{\dd(\Mod\Lambda_{\dn{p}})}(X_{\dn{p}},Y_{\dn{p}})&&\mbox{for }(3),\\
\hom_{\dd(\Mod\Lambda)}(X,Y)\simeq\hom_{\dd(\Mod\Lambda)}(X,Y)^{\widehat{\ }}\simeq\hom_{\dd(\Mod\widehat{\Lambda})}(\widehat{X},\widehat{Y})&&\mbox{for }(4),
\end{eqnarray*}}
where the direct sum is finite. Thus the `if' part follows. The `only if' part also follows since the induced functors $\dd^{\b}(\flmod\Lambda)\to\dd^{\b}(\flmod\Lambda\otimes_RS)$ and $\kk^{\b}(\pr\Lambda)\to\kk^{\b}(\pr\Lambda_{\dn{p}})$ are dense.

(5) We can assume that $R$ is complete local by (3) and (4). For any $X,Y\in\flmod\Lambda$, we have
{\small\[\Ext^{i}_\Lambda(X,Y)\simeq D\ext^{n-i}_\Lambda(Y,X)=\left\{\begin{array}{cc}
0&(i>n),\\
D\hom_\Lambda(Y,X)&(i=n).\\
\end{array}\right.\]}
Considering a minimal projective resolution $\cdots\to P_1\to P_0\to X\to0$ of $X\in\flmod\Lambda$, we have $\hom_\Lambda(P_i,\Lambda/J_\Lambda)=\ext^i_\Lambda(X,\Lambda/J_{\Lambda})=0$ for $i>n$. Thus we see that $P_{n+1}=0$ and $\pd{}_\Lambda X\le n$. Putting $Y:=X$, we see that $\pd{}_\Lambda X=n$. By \XBB, we have $\gl\Lambda=n$.

(6) We can assume that $R$ is local by (3). For any $X\in\flmod\Lambda$, we have
{\small\[\ext^i_\Lambda(X,\Lambda)\simeq D\ext^{n-i}_\Lambda(\Lambda,X)=\left\{\begin{array}{cc}
0&(i\neq n),\\
DX&(i=n).\\
\end{array}\right.\]}
This implies $\depth\Lambda=\dim \Lambda=\id{}_\Lambda\Lambda=n $ by \XBA.

(7) The `if' part holds since $\gl\Lambda<\infty$ implies $\dd^{\b}(\flmod\Lambda)\subset\dd^{\b}(\mod\Lambda)=\kk^{\b}(\pr\Lambda)$. We now show the `only if' part. We can assume that $R$ is local by (3). Then $\Lambda^{(i)}:=\Lambda/J_\Lambda^i$ satisfies $\widehat{\Lambda}=\plim_{i\ge0}\Lambda^{(i)}$. Take $X\in \dd^{\b}(\flmod\Lambda)$ and $Y\in\kk^{\b}(\pr\Lambda)$. Then we have $Y^{(i)}\in \dd^{\b}(\flmod\Lambda)$. Since $\Lambda$ is $n$-CY, we have a functorial isomorphism
$D\hom_{\dd(\Mod\Lambda)}(Y^{(i)},X)\simeq\hom_{\dd(\Mod\Lambda)}(X,Y^{(i)}[n])$.
Taking $\plim$ on both sides and applying \XBF, we obtain functorial isomorphisms
{\small\[D\hom_{\dd(\Mod\Lambda)}(Y,X)\simeq\hom_{\dd(\Mod\Lambda)}(X,\widehat{Y}[n])=\hom_{\dd(\Mod\Lambda)}(X,Y[n])^{\widehat{\ }}=\hom_{\dd(\Mod\Lambda)}(X,Y[n]).\]}

\vskip-1em
(8) $\Lambda$ is $n$-CY$^-$ and $\dd^{\b}(\mod\Lambda)=\kk^{\b}(\pr\Lambda)$ by (7). Thus the assertion follows.\rule{5pt}{10pt}

\vskip.5em
We now state our main theorems in this section, which give a characterization of $n$-CY and $n$-CY$^-$ algebras.

\vskip.5em{\bf Theorem \XCB\ }{\it
Let $R$ be a local Gorenstein ring with $\dim R=d$ and $\Lambda$ a module-finite $R$-algebra. Assume that the structure morphism $R\to\Lambda$ is injective.

(1) If $\Lambda$ is $n$-CY or $n$-CY$^-$ for some integer $n$, then $n=d$.

(2) $\Lambda$ is $d$-CY$^-$ if and only if $\Lambda$ is a symmetric $R$-order (in the sense of section \XB).

%(ii) $\Lambda$ is $d$-CY$^-$ if and only if $\Lambda_{\dn{p}}$ is a symmetric $R_{\dn{p}}$-order for any $\dn{p}\in\Max R$ if and only if $\widehat{\Lambda}_{\dn{p}}$ is a symmetric $\widehat{R}_{\dn{p}}$-order for any $\dn{p}\in\Max R$.

(3) $\Lambda$ is $d$-CY if and only if $\Lambda$ is a symmetric $R$-order with $\gl\Lambda=d$.}
%(iii) $\Lambda$ is $d$-CY if and only if $\Lambda_{\dn{p}}$ is a symmetric $R_{\dn{p}}$-order with $\gl\Lambda_{\dn{p}}=d$ for any $\dn{p}\in\Max R$ if and only if $\widehat{\Lambda}_{\dn{p}}$ is a symmetric $\widehat{R}_{\dn{p}}$-order with $\gl\widehat{\Lambda}_{\dn{p}}=d$ for any $\dn{p}\in\Max R$.

\vskip.5em{\bf Theorem \XCC\ }{\it
Let $R$ be a Gorenstein ring with $\dim R=d$, $\Lambda$ a module-finite $R$-algebra and $n$ an integer. 

(1) The conditions (i)--(iii) are equivalent.

\strut\kern1em
(i) $\Lambda$ is $n$-CY$^-$.

\strut\kern1em
(ii) $\Lambda_{\dn{p}}$ is a CM $R_{\dn{p}}$-module of dimension $n$ and $\ext^{\height\dn{p}-n}_{R_{\dn{p}}}(\Lambda_{\dn{p}},R_{\dn{p}})\simeq\Lambda_{\dn{p}}$ as $(\Lambda_{\dn{p}},\Lambda_{\dn{p}})$-modules for any $\dn{p}\in\Max R\cap\Supp{}_R\Lambda$.

\strut\kern1em
(iii) $\widehat{\Lambda}_{\dn{p}}$ is a CM $\widehat{R}_{\dn{p}}$-module of dimension $n$ and $\ext^{\height\dn{p}-n}_{\widehat{R}_{\dn{p}}}(\widehat{\Lambda}_{\dn{p}},\widehat{R}_{\dn{p}})\simeq\widehat{\Lambda}_{\dn{p}}$ as $(\widehat{\Lambda}_{\dn{p}},\widehat{\Lambda}_{\dn{p}})$-modules for any $\dn{p}\in\Max R\cap\Supp{}_R\Lambda$.

(2) The conditions (i)--(iii) are equivalent.

\strut\kern1em
(i) $\Lambda$ is $n$-CY.

\strut\kern1em
(ii) (1)(ii) and $\gl\Lambda_{\dn{p}}=n$ for any $\dn{p}\in\Max R\cap\Supp{}_R\Lambda$.

\strut\kern1em
(iii) (1)(iii) and $\gl\widehat{\Lambda}_{\dn{p}}=n$ for any $\dn{p}\in\Max R\cap\Supp{}_R\Lambda$.}

\vskip.5em
We note that a $d$-CY$^-$ algebra $\Lambda$ over a non-local Gorenstein ring $R$ is not necessarily a symmetric $R$-algebra even if the structure morphism $R\to\Lambda$ is injective.\footnote{The following example was pointed to us by J. Miyachi: Let $R$ be a noetherian ring with a non-trivial Picard group and $I$ a non-free invertible ideal. Then $\Lambda:=R\ltimes I$ is not a symmetric $R$-algebra even though $\Lambda_{\dn{p}}$ is a symmetric $R_{\dn{p}}$-algebra for any $\dn{p}\in\Max R$. If in addition $R$ is Gorenstein and $\height\dn{p}=d$ for any $\dn{p}\in\Max R$, then $\Lambda$ is $d$-CY$^-$ by \XCB.}

Before proving our main theorems, we state some easy consequences.

\vskip.5em{\bf Corollary \XCD\ }{\it
Let $R$ be a Gorenstein ring with $\dim R=d$, $\Lambda$ a module-finite $R$-algebra and $n$ an integer. 

(1) $\Lambda$ is $n$-CY$^-$ if and only if so is $\Lambda^{\op}$ (cf. \XCA(2)).

(2) If $R$ is local, then $\Lambda$ is $n$-CY$^-$ if and only if so is $\widehat{\Lambda}$ (cf. \XCA(4)).

(3) If $\Lambda$ is $n$-CY (resp. $n$-CY$^-$), then $\Lambda_{\dn{p}}$ and $\widehat{\Lambda}_{\dn{p}}$ are $m$-CY (resp. $m$-CY$^-$) for any $\dn{p}\in\Supp{}_R\Lambda$ and $m:=\dim{}_{R_{\dn{p}}}\Lambda_{\dn{p}}$.

(4) If $R$ is a normal domain and $\Lambda$ is $d$-CY ($d\ge1$), then $\Lambda$ is a reflexive $R$-module and a maximal $R$-order.

(5) If $R$ is local and $\Lambda$ is $d$-CY$^-$, then the following assertions hold for any $i$.

\strut\kern1em
(i) There exists a functorial isomorphism $\ext^i_\Lambda(-,\Lambda)\simeq\ext^i_R(-,R)$ on $\mod\Lambda$.

\strut\kern1em
(ii) $\depth M=d-\sup\{i\ge0\ |\ \ext^i_\Lambda(M,\Lambda)\neq0\}$ and $\dim M=d-\inf\{i\ge0\ |\ \ext^i_\Lambda(M,\Lambda)\neq0\}$ for any $M\in\mod\Lambda$.

\strut\kern1em
(iii) $\ext^{d-i}_\Lambda(-,\Lambda)$ gives a duality between CM $\Lambda$-modules of dimension $i$ and CM $\Lambda^{\op}$-modules of dimension $i$.

(6) If $R$ is local, then symmetric $R$-orders (resp. symmetric $R$-orders of global dimension $d$) are closed under derived equivalences.}

\vskip.5em{\sc Proof }
(1) follows from the left-right symmetry of the condition \XCC(1)(ii), and (2) follows from the equivalence of \XCC(1)(ii) and (iii).
%(i) We can assume $R$ is complete local by \XCA(4). Since $\Lambda$ is symmetric (resp. $\Lambda\in\cm R$, $\gl\Lambda=d$) if and only if so is $\Lambda^{\op}$, the assertion follows. (Of course, the $d$-CY assertion also follows from Matlis dual.)

(3) By \XCA(7), we only have to show the assertion for CY$^-$. For any $\dn{p}\in\Supp{}_R\Lambda$, take $\dn{q}\in\Max R$ with $\dn{p}\subseteq\dn{q}$. Since $\ext^i_{R_{\dn{q}}}(\Lambda_{\dn{q}},R_{\dn{q}})=0$ ($i\neq\height\dn{q}-n$) and $\ext^{\height\dn{q}-n}_{R_{\dn{q}}}(\Lambda_{\dn{q}},R_{\dn{q}})\simeq\Lambda_{\dn{q}}$ as $(\Lambda_{\dn{q}},\Lambda_{\dn{q}})$-modules by \XCC(1), we have that $\ext^i_{R_{\dn{p}}}(\Lambda_{\dn{p}},R_{\dn{p}})=0$ ($i\neq\height\dn{q}-n$) and $\ext^{\height\dn{q}-n}_{R_{\dn{p}}}(\Lambda_{\dn{p}},R_{\dn{p}})\simeq\Lambda_{\dn{p}}$ as $(\Lambda_{\dn{p}},\Lambda_{\dn{p}})$-modules. Thus $\Lambda_{\dn{p}}$ and $\widehat{\Lambda}_{\dn{p}}$ are $m$-CY$^-$ by \XCC(1).

 %Since $\Lambda$ is a CM $R$-module of dimension $\dim M$, it follows from [Ma;17.3] that $\Lambda_{\dn{p}}$ is a CM $R_{\dn{p}}$-module of dimension $m$. assume that $\Lambda$ is $d$-CY. Since $\Lambda$ is a symmetric $R$-algebra by \XC, $\Lambda_{\dn{p}}$ is a symmetric $R_{\dn{p}}$-algebra.  Hence $\Lambda_{\dn{p}}$ is a symmetric $R_{\dn{p}}$-order, and so $\Lambda_{\dn{p}}$ is $n$-CY$^-$. For the $d$-CY case, we only need to observe that $\gl\Lambda_{\dn{p}}<\infty$. The assertion for $\widehat{\Lambda}_{\dn{p}}$ follows using \XCA(4).

(4) Since $\Lambda$ is a Cohen-Macaulay $R$-module, it is reflexive. By (3), $\widehat{\Lambda}_{\dn{p}}$ is $1$-CY for any $\dn{p}\in\Spec R$ with $\height\dn{p}=1$. Thus $\widehat{\Lambda}_{\dn{p}}$ is a maximal $\widehat{R}_{\dn{p}}$-order (see \XCM). These conditions imply that $\Lambda$ is a maximal order by a result of Auslander-Goldman [Re;11.4,11.5].

(5) $\Lambda^*\simeq\Lambda$ as $(\Lambda,\Lambda)$-modules by \XCB. Let $\cdots\to P_0\to M\to0$ be a projective resolution of a $\Lambda$-module $M$. Since $\Lambda\in\cm R$ by \XCB, we have $\ext^i_R(P_i,R)=0$ for any $i>0$. Applying $\hom_\Lambda(-,\Lambda)\simeq\hom_\Lambda(-,\Lambda^*)\simeq\hom_R(-,R)$, we see that $\ext^i_\Lambda(M,\Lambda)=\ext^i_R(M,R)$. The other assertions follow immediately.

(6) Since module-finite $R$-algebras are closed under derived equivalences [Ri1], the assertion follows by \XCA(1) and \XCB.\rule{5pt}{10pt}

\vskip.5em
For proving our main result in this section, it will be useful to investigate Nakayama functors in the context of derived categories. Now let $R$ be a local Gorenstein ring with $\dim R=d$ and $\Lambda$ a module-finite $R$-algebra. Recall that $R$ is a dualizing complex of $R$, i.e.
\[(-)^\dagger:=\Rhom_R(-,R):\dd(\mod R)\to \dd(\mod R)\]
gives a duality such that $(-)^{\dagger\dagger}$ is isomorphic to the identity functor [Har1;V2.1]. Obviously $(-)^{\dagger}$ induces a duality $(-)^\dagger:\dd^\pm(\mod\Lambda)\leftrightarrow \dd^\mp(\mod\Lambda^{\op})$. Define the {\it Nakayama functor} in the derived category by the composition
\[\nu:\dd^-(\mod\Lambda)\stackrel{\Rhom_\Lambda(-,\Lambda)}{\longrightarrow}\dd^+(\mod\Lambda^{\op})\stackrel{(-)^\dagger}{\longrightarrow}\dd^-(\mod\Lambda).\]
%Then $(-)^\dagger$ gives a duality \[(-)^\dagger:\dd^\pm(\mod\Lambda)\leftrightarrow \dd^\mp(\mod\Lambda^{\op}).\] We also have functors \begin{eqnarray*}\Rhom_\Lambda(-,\Lambda)&:&\dd^-(\mod\Lambda)\to \dd^+(\mod\Lambda^{\op}),\\ \Rhom_{\Lambda^{\op}}(-,\Lambda)&:&\dd^-(\mod\Lambda^{\op})\to \dd^+(\mod\Lambda). \end{eqnarray*}
%\begin{eqnarray*}\nu^-&:&\dd^+(\mod\Lambda)\stackrel{(-)^\dagger}{\longrightarrow}\dd^-(\mod\Lambda^{\op})\stackrel{\Rhom_{\Lambda^{\op}}(-,\Lambda)}{\longrightarrow}\dd^+(\mod\Lambda).\end{eqnarray*}

\vskip.5em{\bf Proposition \XCE\ }{\it
(1) With the above notation, we have the following isomorphisms of functors, $(-)^\dagger\simeq\Rhom_\Lambda(-,\Lambda^\dagger):\dd(\mod\Lambda)\to \dd(\mod\Lambda^{\op})$ and $(-)^\dagger\simeq\Rhom_{\Lambda^{\op}}(-,\Lambda^\dagger):\dd(\mod\Lambda^{\op})\to \dd(\mod\Lambda)$. Thus $\Lambda^\dagger\in \dd^{\b}(\mod\Lambda\otimes_R\Lambda^{\op})$ is a dualizing complex of $\Lambda$ in the sense of Yekutieli [Ye1].

(2) There exists an isomorphism $\nu\simeq\Lambda^\dagger\Lotimes_\Lambda(-)$ of functors on $\dd^-(\mod\Lambda)$.

(3) There exists a functorial isomorphism $\Rhom_\Lambda(X,\nu(Y))\simeq\Rhom_\Lambda(Y,X)^\dagger$ for any $X\in \dd^{\b}(\mod\Lambda)$ and $Y\in\kk^{\b}(\pr\Lambda)$.}

\vskip.5em{\sc Proof }
(1) $\Rhom_\Lambda(-,\Lambda^\dagger)\stackrel{{\rm\XBJ(1)}}{\simeq}(\Lambda\Lotimes_\Lambda-)^\dagger\simeq(-)^\dagger$.

(2) $\nu=\Rhom_\Lambda(-,\Lambda)^\dagger\stackrel{}{\simeq}\Rhom_\Lambda(-,\Lambda^{\dagger\dagger})^\dagger\stackrel{{\rm\XBJ(1)}}{\simeq}(\Lambda^\dagger\Lotimes_\Lambda-)^{\dagger\dagger}\stackrel{}{\simeq}\Lambda^\dagger\Lotimes_\Lambda-$.

(3) $\Rhom_\Lambda(X,\nu(Y))=\Rhom_\Lambda(X,\Rhom_\Lambda(Y,\Lambda)^\dagger)\stackrel{{\rm\XBJ(1)}}{\simeq}(\Rhom_\Lambda(Y,\Lambda)\Lotimes_\Lambda X)^\dagger$

$\stackrel{{\rm\XBJ(2)}}{\simeq}\Rhom_\Lambda(Y,X)^\dagger$.\rule{5pt}{10pt}

\vskip.5em
We need the following special case of the local duality theorem [F][Har1].

\vskip.5em{\bf Lemma \XCF\ }{\it
We have an isomorphism $(-)^\dagger\simeq [-d]\circ D$ of functors on $\dd^{\b}(\flmod R)$.}

\vskip.5em{\sc Proof }
We give a proof for completeness. Fix $X\in \dd^{\b}(\flmod R)$. Consider the following morphism in $\dd^{\b}(\Mod R)$
{\small\[\begin{diag}
I&&\cdots&\RA{}&I^0&\RA{}&I^1&\RA{}&\cdots&\RA{}&I^{d-1}&\RA{}&I^d&\RA{}&0&\cdots\\
\uparrow^a&&&&\uparrow&&\uparrow&&&&\uparrow&&\parallel\\
I^d[-d]&&\cdots&\RA{}&0&\RA{}&0&\RA{}&\cdots&\RA{}&0&\RA{}&I^d&\RA{}&0&\cdots,
\end{diag}\]}
where $I$ is a minimal injective resolution of the $R$-module $R$. Then $I^d=E$. Take a triangle
$E[-d]\stackrel{a}{\to}I\to I^\prime\to E[1-d]$.
Since $\hom_R(X,I^i)=0$ for any $X\in\flmod R$ and for any $i$ with $0\le i<d$, we have $\Rhom_R(X,I^\prime)=0$ for any $X\in \dd^{\b}(\flmod R)$. Thus we have an isomorphism $(DX)[-d]=\Rhom_R(X,E[-d])\simeq\Rhom_R(X,I)=X^\dagger$.\rule{5pt}{10pt}
%since $I$ is isomorphic to $R$ in $\dd^{\b}(\Mod R)$

\vskip.5em
We have the following `Serre duality theorem' for arbitrary module-finite $R$-algebras.

\vskip.5em{\bf Theorem \XCG\ }{\it
For any module-finite $R$-algebra $\Lambda$, we have a functorial isomorphism
\[\hom_{\dd(\Mod\Lambda)}(X,\nu(Y)[d])\simeq D\hom_{\dd(\Mod\Lambda)}(Y,X)\]
for any $X\in \dd^{\b}(\flmod\Lambda)$ and $Y\in\kk^{\b}(\pr\Lambda)$,}

\vskip.5em{\sc Proof }
$\hom_{\dd(\Mod\Lambda)}(X,\nu(Y)[d])=H^d(\Rhom_\Lambda(X,\nu(Y)))\stackrel{{\rm\XCE(3)}}{\simeq}H^d(\Rhom_\Lambda(Y,X)^\dagger)$

$\stackrel{{\rm\XCF}}{\simeq}H^0(D\Rhom_\Lambda(Y,X))=DH^0(\Rhom_\Lambda(Y,X))=D\hom_{\dd(\Mod\Lambda)}(Y,X)$.\rule{5pt}{10pt}

\vskip.5em
We now obtain the following crucial result to prove our main theorems \XCB\ and \XCC.
%We now obtain our desired characterization of $d$-CY$^-$ algebras in terms of symmetric orders, by at the same time showing that this is equivalent to the Nakayama functor $\nu$ being isomorphic to the identity functor on $\kk^{\b}(\pr\Lambda)$.

\vskip.5em{\bf Theorem \XCH\ }{\it
Let $R$ be a local Gorenstein ring with $\dim R=d$, $\Lambda$ a module-finite $R$-algebra and $n$ an integer. Then the following conditions are equivalent.

(1) $\Lambda$ is $n$-CY$^-$.

%(2) There exists a functorial isomorphism $\hom_{\dd(\Mod\Lambda)}(X,Y[n])\simeq D\hom_{\dd(\Mod\Lambda)}(Y,X)$ for any $X,Y\in \dd^{\b}(\flmod\Lambda)\cap\kk^{\b}(\pr\Lambda)$.

(2) There exist isomorphisms $\ext^n_\Lambda(-,\Lambda)\simeq D$ and $\ext^i_\Lambda(-,\Lambda)=0$ ($i\neq n$) of functors on $\flmod\Lambda$ which commute with the right action of $\Lambda$.

(3) $\Lambda^\dagger\simeq\Lambda[n-d]$ in $\dd^{\b}(\mod\Lambda\otimes_R\Lambda^{\op})$.

(4) There exists an isomorphism $\nu\simeq[n-d]$ of functors on $\kk^{\b}(\pr\Lambda)$.

(5) $\Lambda$ is a CM $R$-module of dimension $n$, and $\ext^{d-n}_R(\Lambda,R)\simeq\Lambda$ as $(\Lambda,\Lambda)$-modules.

($\widehat{i}$) The condition ($i$) replacing $(R,\Lambda)$ by $(\widehat{R},\widehat{\Lambda})$ ($1\le i\le 5$).}

\vskip.5em{\sc Proof }
(3)$\Leftrightarrow$(4) (resp. ($\widehat{3}$)$\Leftrightarrow$($\widehat{4}$)) follows by \XCE(2). (4)$\Rightarrow$(1) (resp. ($\widehat{4}$)$\Rightarrow$($\widehat{1}$)) follows by \XCG. (1)$\Rightarrow$(2) (resp. ($\widehat{1}$)$\Rightarrow$($\widehat{2}$)) and (2)$\Rightarrow$($\widehat{2}$) are obvious.

%(2)$\Rightarrow$(4) We have a functorial isomorphism \[\hom_{\dd^{\b}(\mod\Lambda)}(X,\nu(Y)[d])\stackrel{\XCG}{\simeq}D\hom_{\dd^{\b}(\mod\Lambda)}(Y,X)\stackrel{{\rm(2)}}{\simeq}\hom_{\dd^{\b}(\mod\Lambda)}(X,Y[n])\] for any $X\in \dd^{\b}(\flmod\Lambda)$ and $Y\in\kk^{\b}(\pr\Lambda)$. Thus we have an isomorphism $[d]\circ\nu\simeq[n]$ of functors on $\dd^{\b}(\flmod\Lambda)\cap\kk^{\b}(\pr\Lambda)$.
%(3)$\Rightarrow$(4) We have $\Rhom_\Lambda(-,\Lambda)\simeq[-n]\circ D\stackrel{\XCF}{\simeq}[d-n]\circ(-)^\dagger$ on $\flmod\Lambda$.

($\widehat{2}$)$\Rightarrow$($\widehat{3}$) We put $(R,\Lambda):=(\widehat{R},\widehat{\Lambda})$ for simplicity. It follows from \XBA\ that ${\Lambda}$ is a CM ${R}$-module of dimension $n$. Take a ${\Lambda}$-regular sequence $a_1,\cdots,a_n\in {R}$, and put $I_i^\prime:=\sum_{j=1}^n{R}a_j^i$ and $I_i:=I_i^\prime{\Lambda}$. Then ${\Lambda}/I_i\in\flmod{\Lambda}$ and ${\Lambda}=\plim_{i\ge0}\ {\Lambda}/I_i$. Since ${\Lambda}\Lotimes_{{R}}{R}/I_i^\prime={\Lambda}/I_i={R}/I_i^\prime\Lotimes_{{R}}{\Lambda}$ holds, we have ${\Lambda}^\dagger\Lotimes_{{\Lambda}}{\Lambda}/I_i={\Lambda}^\dagger\Lotimes_{{R}}{R}/I_i^\prime={R}/I_i^\prime\Lotimes_{{R}}{\Lambda}^\dagger={\Lambda}/I_i\Lotimes_{{\Lambda}}{\Lambda}^\dagger$. 
We have isomorphisms ${\Lambda}^\dagger\Lotimes_{{\Lambda}}(-)\stackrel{\XCE(2)}{\simeq}\nu=(-)^\dagger\circ\Rhom_{{\Lambda}}(-,{\Lambda})\stackrel{(\widehat{2})}{\simeq}(-)^\dagger\circ[-n]\circ D\stackrel{\XCF}{\simeq}[n-d]$ of functors on $\flmod{\Lambda}$. Thus we have isomorphisms
{\small\[{\Lambda}[n-d]\stackrel{{\rm\XBF(1)}}{=}\plim_{i\ge0}\ {\Lambda}/I_i[n-d]\stackrel{}{\simeq}\plim_{i\ge0}\ {\Lambda}^\dagger\Lotimes_{{\Lambda}}{\Lambda}/I_i\simeq\plim_{i\ge0}\ {\Lambda}/I_i\Lotimes_{{\Lambda}}{\Lambda}^\dagger\stackrel{{\rm\XBF(1)}}{=}{\Lambda}^\dagger.\]}
These isomorphisms commute with the right multiplication of ${\Lambda}$. Thus ($\widehat{3}$) holds.

(4)$\Leftrightarrow$(5) (resp. ($\widehat{4}$)$\Leftrightarrow$($\widehat{5}$)) We have $H^i(\Lambda^\dagger)=\ext^i_R(\Lambda,R)$. Thus $H^i(\Lambda^\dagger)=0$ holds for any $i\neq d-n$ if and only if $\Lambda$ is a CM $R$-module of dimension $n$. In this case, $\Lambda^\dagger\simeq\Lambda$ as $(\Lambda,\Lambda)$-modules if and only if $\ext^{d-n}_R(\Lambda,R)\simeq\Lambda$ as $(\Lambda,\Lambda)$-modules.

(5)$\Leftrightarrow$($\widehat{5}$) $\Lambda$ is a CM $R$-module of dimension $n$ if and only if $\widehat{\Lambda}$ is a CM $\widehat{R}$-module of dimension $n$. Since $(\Lambda\otimes_R\Lambda^{\op})^{\widehat{\ }}=\widehat{\Lambda}\otimes_{\widehat{R}}\widehat{\Lambda}^{\op}$ holds, it follows from \XCI\ below that $\Lambda\simeq\hom_R(\Lambda,R)$ as $(\Lambda\otimes_R\Lambda^{\op})$-modules if and only if $\widehat{\Lambda}\simeq\hom_{\widehat{R}}(\widehat{\Lambda},\widehat{R})$ as $(\widehat{\Lambda}\otimes_{\widehat{R}}\widehat{\Lambda}^{\op})$-modules.

Now one can easily check that all conditions are equivalent.\rule{5pt}{10pt}

\vskip.5em{\bf Lemma \XCI\ }{\it
Let $R$ be a local ring, $\Lambda$ a module-finite $R$-algebra and $M,N\in\mod\Lambda$. If $\widehat{M}\simeq\widehat{N}$ as $\widehat{\Lambda}$-modules, then $M\simeq N$ as $\Lambda$-modules.}

\vskip.5em{\sc Proof }
We modify the proof of [CRe;30.17] where the case $\dim R=1$ is treated. Let $f\in\hom_{\widehat{\Lambda}}(\widehat{M},\widehat{N})$ be an isomorphism with $g:=f^{-1}$. Since $\hom_{\widehat{\Lambda}}(\widehat{M},\widehat{N})=\hom_\Lambda(M,N)^{\widehat{\ }}$, we can take $f'\in\hom_\Lambda(M,N)$ and $g'\in\hom_\Lambda(N,M)$ with $f-f'\in\dn{p}\hom_\Lambda(M,N)$ and $g-g'\in\dn{p}\hom_\Lambda(N,M)$. Then $g'f'-1_N=g'f'-gf\in(\dn{p}\endm_{\widehat{\Lambda}}(\widehat{N}))\cap\endm_\Lambda(N)=\dn{p}\endm_\Lambda(N)$. Thus we have $N=g'f'(N)+\dn{p}N$. Using Nakayama's lemma, we have $g'f'\in\aut_\Lambda(N)$. Similarly, $f'g'\in\aut_\Lambda(M)$ holds, so $f'$ is an isomorphism.\rule{5pt}{10pt}

%\vskip.5em As a direct consequence we get our desired characterization of $d$-CY algebras. \vskip.5em{\bf Theorem }{\it Under the same assumptions in \XCH, $\Lambda$ is $n$-CY if and only if $\Lambda$ is a CM $R$-module of dimension $n$, $\ext^{d-n}_R(\Lambda,R)$ is isomorphic to $\Lambda$ as a $(\Lambda,\Lambda)$-module and $\gl\Lambda=n$.} \vskip.5em{\sc Proof } (1)$\Leftrightarrow$(2) follows from \XCA(7) and \XCH(1)$\Leftrightarrow$(2). (3)$\Rightarrow$(1) Since $\Lambda$ is $d$-CY$^-$ by \XCH\ and $\kk^{\b}(\pr\Lambda)=\dd^{\b}(\mod\Lambda)$ holds by $\gl\Lambda<\infty$, it is $d$-CY. (1)$\Rightarrow$(3) $\Lambda$ is a symmetric $R$-order by \XCH, and $\gl\Lambda=d$ holds by \XCA(5).\rule{5pt}{10pt}

\vskip.5em
Using \XCH\ together with \XCA(7), we have now completed the proof of \XCB\ and \XCC. 

In the rest of this section, we give some examples of $n$-CY algebras. Let us start with considering commutative CY and CY$^-$ algebras.

\vskip.5em{\bf Proposition \XCJ\ }{\it
Let $R$ be a commutative noetherian ring and $n$ an integer.

(1) $R$ is $n$-CY$^-$ if and only if $R$ is Gorenstein and $\dim R_{\dn{p}}=n$ for any $\dn{p}\in\Max R$.

(2) $R$ is $n$-CY if and only if $R$ is regular and $\dim R_{\dn{p}}=n$ for any $\dn{p}\in\Max R$.}

\vskip.5em{\sc Proof }
The `only if' part follows from \XCA(3)(5) and (6). The `if' part follows from \XCA(3) and \XCB.\rule{5pt}{10pt}

\vskip.5em
Next we consider $0$-CY algebras.

\vskip.5em{\bf Proposition \XCK\ }{\it
A finite dimensional algebra over a field is 0-CY if and only if it is a semisimple algebra.}

\vskip.5em{\sc Proof }
Since any 0-CY algebra has global dimension zero, it is semisimple. Conversely, it is well-known that any semisimple algebra over a field is symmetric [CRe;9.8].\rule{5pt}{10pt}

\vskip.5em
Now we consider 1-CY algebras over a complete discrete valuation ring $R$ with quotient field $K$. Recall that an $R$-order $\Lambda$ is called {\it hereditary} if $\gl\Lambda=1$. Let us recall briefly the structure theory of maximal and hereditary orders [Re][CRe]. For a ring $\Delta$ and $n>0$, we put $\tri_n(\Delta):=\{(x_{ij})_{1\le i,j\le n}\in\ma_n(\Delta)\ |\ x_{ij}\in J_{\Delta}\mbox{ if }i>j\}$. The following results are well-known. 

\vskip.5em{\bf Proposition \XCL\ }{\it
(1) Any finite dimensional division $K$-algebra $D$ contains a unique maximal $R$-order $\Delta_D$ [Re;12.8].

(2) An $R$-order is maximal if and only if it is Morita equivalent to $\Delta_{D_1}\times\cdots\times\Delta_{D_k}$ for some finite dimensional division $K$-algebras $D_i$ [Re;17,3].

(3) An $R$-order is hereditary if and only if it is Morita equivalent to $\tri_{n_1}(\Delta_{D_1})\times\cdots\times\tri_{n_k}(\Delta_{D_k})$ for some finite dimensional division $K$-algebras $D_i$ and $n_i>0$ [Re;39.14].}

%Let us recall the classical structure theorem of hereditary orders [CR??].  and {\it maximal} if there is no $R$-order $\Gamma$ satisfying $\Lambda\subset\Gamma\subset\Lambda\otimes_RK$ for a quotient field $K$ of $R$ except $\Lambda$.

\vskip.5em
We have the following relationship between 1-CY algebras and maximal orders.

\vskip.5em{\bf Proposition \XCM\ }{\it
Let $R$ be a complete discrete valuation ring and $\Lambda$ a module-finite $R$-algebra. If $\Lambda$ is 1-CY, then it is a maximal $R$-order.}

\vskip.5em{\sc Proof }
By \XCB, $\Lambda$ is a symmetric $R$-order with $\gl\Lambda=1$. Thus $\Lambda$ is Morita equivalent to $\tri_{n_1}(\Delta_{D_1})\times\cdots\times\tri_{n_k}(\Delta_{D_k})$ by \XCL(3). One can check that if $\tri_n(\Delta)$ is a symmetric $R$-algebra, then $n=1$ (e.g. [Hae;6.3]). Thus $\Lambda$ is maximal by \XCL(2).\rule{5pt}{10pt}

\vskip.5em
We note that a maximal order is not necessarily symmetric. Let $D$ be a central division $K$-algebra with $\dim_KD=n^2$. If the residue field of $R$ is finite, then $\hom_R(\Delta_D,R)$ is isomorphic to $J_{\Delta_D}^{1-n}$ as a $(\Delta_D,\Delta_D)$-module [Re;14.9]. If $n>1$, then $J_{\Delta_D}^{1-n}$ is never isomorphic to $\Delta_D$ by [Re;37.27]. We thank Wolfgang Rump for kindly explaining these results to us.

\vskip.5em
We now give other examples of $d$-CY algebras, where $d\ge2$. Let $K$ be a field of characteristic zero and $G$ a finite subgroup of $\SL_d(K)$ acting on $K^d$ naturally. The action of $G$ naturally extends to $S:=K[[x_1,\cdots,x_d]]$. We denote by $S^G$ the invariant subring, and by $S*G$ the skew group ring, i.e. a free $S$-module with a basis $G$, where the multiplication is given by $(s_1g_1)\cdot(s_2g_2)=(s_1g_1(s_2))(g_1g_2)$ for $s_i\in S$ and $g_i\in G$. We have the following result (c.f. [CRo]).

\vskip.5em{\bf Theorem \XCN\ }{\it 
$S*G$ is $d$-CY and a symmetric $S^G$-order with $\gl S*G=d$.}

\vskip.5em{\sc Proof }
Any finite subgroup $G$ of $\SL_d(K)$ is small in the sense that any $g\in G$ with $g\neq1$ satisfies $\rank(g-1)>1$. This implies $\endm_{S^G}(S)=S*G$ by a result of Auslander. See [A3][Yo] for $d=2$, and a similar argument works for arbitrary $d$.

Since $S*G$ is a free $S$-module, $S*G$ is an $S^G$-order. Since $S*G=\endm_{S^G}(S)$ holds, $S*G$ is a symmetric $S^G$-order by \XBD(3). Moreover, $\Ext^i_{S*G}(X,Y)=\Ext^i_S(X,Y)^G$ holds for any $X,Y\in\mod S*G$ and $i\in\zzz$ [A3][Yo]. Thus we obtain $\gl S*G=d$ by $\gl S=d$ and \XBA\ (see also [RR]).\rule{5pt}{10pt}

\vskip.5em
For a finite subgroup  $G$ of $\SL_d(K)$, we can draw the quiver of the algebra $S*G$ as the {\it McKay quiver} of $G$ [Mc] by using irreducible representations of $G$ and tensor products (see [A3][Yo][I4]). Now we give some examples. If $G=\langle{\rm diag}(\zeta,\zeta,\cdots,\zeta)\rangle\subset\SL_d(K)$ with $\zeta^d=1$, then $S*G$ is the completion of the path algebra of the following quiver with commutative relations $x_ix_j=x_jx_i$ for any $i$ and $j$.
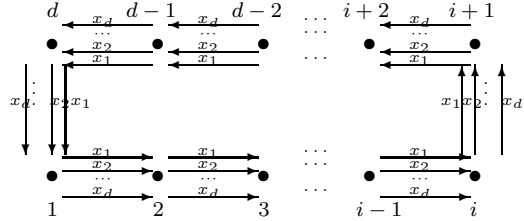
\[\begin{picture}(180,60)
\put(0,-2){\circle*{4}}
\put(-2,-17){\scriptsize 1}
\put(15,6){\tiny $x_1$}
\put(15,1){\tiny $x_2$}
\put(15,-9){\tiny $\stackrel{\cdots}{x_d}$}
\put(4,5){\vector(1,0){34}}
\put(4,0){\vector(1,0){34}}
\put(4,-10){\vector(1,0){34}}
\put(40,-2){\circle*{4}}
\put(38,-17){\scriptsize 2}
\put(55,6){\tiny $x_1$}
\put(55,1){\tiny $x_2$}
\put(55,-9){\tiny $\stackrel{\cdots}{x_d}$}
\put(44,5){\vector(1,0){34}}
\put(44,0){\vector(1,0){34}}
\put(44,-10){\vector(1,0){34}}
\put(80,-2){\circle*{4}}
\put(78,-17){\scriptsize 3}
\put(95,6){\tiny $\cdots$}
\put(95,1){\tiny $\cdots$}
\put(95,-9){\tiny $\stackrel{}{\cdots}$}
%\put(82,5){\vector(1,0){38}}
%\put(82,0){\vector(1,0){38}}
%\put(82,-10){\vector(1,0){38}}
\put(120,-2){\circle*{4}}
\put(115,-17){\scriptsize $i-1$}
\put(135,6){\tiny $x_1$}
\put(135,1){\tiny $x_2$}
\put(135,-9){\tiny $\stackrel{\cdots}{x_d}$}
\put(124,5){\vector(1,0){34}}
\put(124,0){\vector(1,0){34}}
\put(124,-10){\vector(1,0){34}}
\put(160,-2){\circle*{4}}
\put(158,-17){\scriptsize $i$}

\put(155,6){\vector(0,1){34}}
\put(160,6){\vector(0,1){34}}
\put(170,6){\vector(0,1){34}}
\put(147,25){\tiny $x_1$}
\put(155,25){\tiny $x_2$}
\put(163,25){\tiny $\vdots$}
\put(170,25){\tiny $x_d$}

\put(-16,25){\tiny $x_d$}
\put(-8,25){\tiny $\vdots$}
\put(-1,25){\tiny $x_2$}
\put(7,25){\tiny $x_1$}
\put(-10,40){\vector(0,-1){34}}
\put(0,40){\vector(0,-1){34}}
\put(5,40){\vector(0,-1){34}}

\put(0,48){\circle*{4}}
\put(-2,58){\scriptsize $d$}
\put(15,56){\tiny $x_d$}
\put(15,46){\tiny $\stackrel{\cdots}{x_2}$}
\put(15,41){\tiny $x_1$}
\put(38,55){\vector(-1,0){34}}
\put(38,45){\vector(-1,0){34}}
\put(38,40){\vector(-1,0){34}}
\put(40,48){\circle*{4}}
\put(28,58){\scriptsize $d-1$}
\put(55,56){\tiny $x_d$}
\put(55,46){\tiny $\stackrel{\cdots}{x_2}$}
\put(55,41){\tiny $x_1$}
\put(78,55){\vector(-1,0){34}}
\put(78,45){\vector(-1,0){34}}
\put(78,40){\vector(-1,0){34}}
\put(80,48){\circle*{4}}
\put(68,58){\scriptsize $d-2$}
\put(95,56){\tiny $\cdots$}
\put(95,51){\tiny $\cdots$}
\put(95,41){\tiny $\stackrel{}{\cdots}$}
%\put(82,55){\vector(1,0){38}}
%\put(82,50){\vector(1,0){38}}
%\put(82,40){\vector(1,0){38}}
\put(120,48){\circle*{4}}
\put(110,58){\scriptsize $i+2$}
\put(135,56){\tiny $x_d$}
\put(135,46){\tiny $\stackrel{\cdots}{x_2}$}
\put(135,41){\tiny $x_1$}
\put(158,55){\vector(-1,0){34}}
\put(158,45){\vector(-1,0){34}}
\put(158,40){\vector(-1,0){34}}
\put(160,48){\circle*{4}}
\put(150,58){\scriptsize $i+1$}
\end{picture}\]
\vskip.5em
If $d=3$ in the example above, then $S*G$ has the left quiver below. If $G=\langle{\rm diag}(\zeta,\zeta^2,\zeta^2)\rangle\subset\SL_5(K)$ with $\zeta^5=1$, then $S*G$ has the right quiver below.
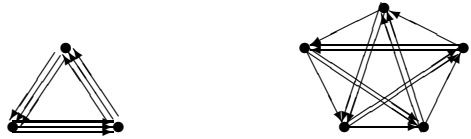
\[\begin{picture}(40,40)
\put(0,0){\circle*{4}}
\put(40,0){\circle*{4}}
\put(20,30){\circle*{4}}
%\put(-7,-7){\scriptsize 1}
%\put(42,-7){\scriptsize 2}
%\put(17,34){\scriptsize 3}
\put(0,2){\vector(1,0){38}}
\put(0,0){\vector(1,0){38}}
\put(0,-2){\vector(1,0){38}}
\put(40,4){\vector(-2,3){17}}
\put(38,3){\vector(-2,3){17}}
\put(36,2){\vector(-2,3){17}}
\put(16,28){\vector(-2,-3){17}}
\put(18,27){\vector(-2,-3){17}}
\put(20,26){\vector(-2,-3){17}}
\end{picture}
\ \ \ \ \ \ \ \ \ \ \ \ \ \ \ \ \ \ 
\begin{picture}(60,47)
\put(15,0){\circle*{4}}
\put(45,0){\circle*{4}}
\put(60,30){\circle*{4}}
\put(30,45){\circle*{4}}
\put(0,30){\circle*{4}}
\put(15,0){\vector(1,0){28}}
\put(45,0){\vector(1,2){14}}
\put(60,30){\vector(-2,1){28}}
\put(30,45){\vector(-2,-1){28}}
\put(0,30){\vector(1,-2){14}}
\put(15,-1){\vector(3,2){43}}
\put(15,1){\vector(3,2){43}}
\put(45,2){\vector(-1,3){14}}
\put(44,-2){\vector(-1,3){15}}
\put(60,29){\vector(-1,0){58}}
\put(60,31){\vector(-1,0){58}}
\put(30,43){\vector(-1,-3){14}}
\put(29,47){\vector(-1,-3){15}}
\put(0,29){\vector(3,-2){43}}
\put(0,31){\vector(3,-2){43}}
\end{picture}\]

%\[\begin{diag} \stackrel{1}{\bullet}&\def\arraystretch{-.2}\begin{array}{c} \stackrel{}{\longrightarrow}\\ \stackrel{}{\longrightarrow}\\ {\scriptstyle\cdots}\\ \stackrel{}{\longrightarrow}\end{array}& \stackrel{2}{\bullet}&\def\arraystretch{-.2}\begin{array}{c} \stackrel{}{\longrightarrow}\\ \stackrel{}{\longrightarrow}\\ {\scriptstyle\cdots}\\ \stackrel{}{\longrightarrow}\end{array}& \stackrel{3}{\bullet}&\cdots\cdots& \stackrel{i-1}{\bullet}&\def\arraystretch{-.2}\begin{array}{c} \stackrel{}{\longrightarrow}\\ \stackrel{}{\longrightarrow}\\ {\scriptstyle\cdots}\\ \stackrel{}{\longrightarrow}\end{array}& \stackrel{i}{\bullet}\\ \uparrow^{}\uparrow^{}{\scriptstyle}\vdots\uparrow^{}&&&&&&&&\downarrow^{}\downarrow^{}{\scriptstyle}\vdots\downarrow^{}\\ \stackrel{d}{\bullet}&\def\arraystretch{-.2}\begin{array}{c}\stackrel{}{\longleftarrow}\\ \stackrel{}{\longleftarrow}\\ {\scriptstyle\cdots}\\ \stackrel{}{\longleftarrow}\end{array}&\stackrel{d-1}{\bullet}&\def\arraystretch{-.2}\begin{array}{c}\stackrel{}{\longleftarrow}\\ \stackrel{}{\longleftarrow}\\ {\scriptstyle\cdots}\\ \stackrel{}{\longleftarrow}\end{array}&\stackrel{d-2}{\bullet}&\cdots\cdots&\stackrel{i+2}{\bullet}&\def\arraystretch{-.2}\begin{array}{c}\stackrel{}{\longleftarrow}\\ \stackrel{}{\longleftarrow}\\ {\scriptstyle\cdots}\\ \stackrel{}{\longleftarrow}\end{array}&\stackrel{i+1}{\bullet}\end{diag}\]

\vskip1.5em{\bf\XD. Construction of tilting modules }

Let $R$ be a complete local ring and $\Lambda$ a ring-indecomposable module-finite $R$-algebra. A central theme in this paper is the study of tilting modules for $d$-CY algebras for $d=2,3$, especially the tilting modules of projective dimension at most one. In particular, we are interested in the number of complements of almost complete tilting modules.
A basic partial tilting $\Lambda$-module $T$ is said to be an {\it almost complete tilting module} if $T$ has $(n-1)$ non-isomorphic indecomposable direct summands, where $n$ is the number of non-isomorphic simple $\Lambda$-modules. In this case, $X$ is called a {\it complement} of $T$ if $T\oplus X$ is a basic tilting $\Lambda$-module. 
For finite dimensional algebras it is known that there are at most two (and at least one) complements, in the case of projective dimension at most one, and it is never the case that all almost complete tilting modules have two complements. In the context of module-finite $R$-algebras $\Lambda$, we see in section \XE\ that the result on at most two complements still holds, but now there are algebras $\Lambda$ where all almost complete tilting modules have two complements, as we shall show in section \XE\ for 2-CY and 3-CY algebras. 

In this section we treat the special case of almost complete tilting modules which are projective. We show that for $d$-CY algebras with no loops in the quiver there are exactly two complements to tilting modules of projective dimension at most one, and give an explicit description of the non-projective one. Even though it will not be used later in this paper, we describe more generally all complements which give tilting complexes, in particular those which are tilting modules of projective dimension greater than one. In the first part of this section we work in the general context of module-finite $R$-algebras, and give here necessary (and sufficient) conditions on what the complements are. Then we use this to obtain a nice description for $d$-CY algebras.

Fix an indecomposable object $P\in\pr\Lambda$, and let $Q$ be a direct sum of all indecomposable projective $\Lambda$-modules which are not isomorphic to $P$. We want to find conditions for replacing $P$ with a complement $X$ of $Q$. For $n\ge0$, there exists a unique complex up to isomorphism
\[{\rm RA}_n:\cdots\to0\to A^0\stackrel{a^0}{\longrightarrow}\cdots\stackrel{a^{n-2}}{\longrightarrow}A^{n-1}\stackrel{a^{n-1}}{\longrightarrow}P\to0\to\cdots\]
%\stackrel{\def\arraystretch{.3}\begin{array}{c}{\scriptstyle P}\\ {\scriptstyle\parallel}\end{array}}{A^n}\to0\to\cdots\]
which gives the first $n$ terms of the minimal right $(\add Q)$-approximation sequence of $P$ [AS], i.e. $A^i\in\add Q$ and $a^i\in J_{\pr\Lambda}$ for any $i$ and $H^i(\Rhom_\Lambda(Q,{\rm RA}_n))=0$ for any $i\neq0$.
%\[\hom_\Lambda(Q,A^0)\stackrel{\cdot f^0}{\longrightarrow}\cdots\stackrel{\cdot f^{n-2}}{\longrightarrow}\hom_\Lambda(Q,A^{1-n})\stackrel{\cdot f^{n-1}}{\longrightarrow}\hom_\Lambda(Q,P)\to0\] is exact. Then ${\rm RA}_n$ is unique up to isomorphism.
Similarly, for $n\ge0$, there exists a unique complex up to isomorphism
\[{\rm LA}_n:\cdots\to0\to P
%\stackrel{\def\arraystretch{.3}\begin{array}{c}{\scriptstyle P}\\ {\scriptstyle\parallel}\end{array}}{B^{-n}}
\stackrel{b^{-n}}{\longrightarrow}B^{1-n}\stackrel{b^{1-n}}{\longrightarrow}\cdots\stackrel{b^{-1}}{\longrightarrow}B^0\to0\to\cdots\]
which gives the first $n$ terms of the minimal left $(\add Q)$-approximation sequence of $P$, i.e. $B^i\in\add Q$ and $b^i\in J_{\pr\Lambda}$ for any $i$ and $H^i(\Rhom_\Lambda({\rm LA}_n,Q))=0$ for any $i\neq0$.
%\[\hom_\Lambda(B^0,Q)\stackrel{g^{-1}\cdot}{\longrightarrow}\cdots\stackrel{g^{1-n}\cdot}{\longrightarrow}\hom_\Lambda(B^{n-1},Q)\stackrel{g^{-n}\cdot}{\longrightarrow}\hom_\Lambda(P,Q)\to0\] is exact. Then ${\rm LA}_n$ is unique up to isomorphism.

We then have the following necessary and sufficient conditions on complements of $Q$.

\vskip.5em{\bf Theorem \XDA\ }{\it
With the above notation and assumptions, we have the following.

(1) Let $X$ be an indecomposable object in $\kk^{\b}(\pr\Lambda)$. If $X\oplus Q$ is a tilting complex, then $X$ is isomorphic to ${\rm RA}_n$ or ${\rm LA}_n$ for some $n\ge0$.

(2) ${\rm RA}_n\oplus Q$ is a tilting complex if and only if $H^i(\Rhom_\Lambda({\rm RA}_n,Q))=0$ for any $i\neq0$.
%\[0\to\hom_\Lambda(P,Q)\stackrel{f^{n-1}\cdot}{\longrightarrow}\hom_\Lambda(A^{n-1},Q)\stackrel{f^{n-2}\cdot}{\longrightarrow}\cdots\stackrel{f^0\cdot}{\longrightarrow}\hom_\Lambda(A^0,Q)\]is exact.

(3) ${\rm LA}_n\oplus Q$ is a tilting complex if and only if $H^i(\Rhom_\Lambda(Q,{\rm LA}_n))=0$ for any $i\neq0$.}
%\[0\to\hom_\Lambda(Q,P)\stackrel{\cdot g^{-n}}{\longrightarrow}\hom_\Lambda(Q,B^{1-n})\stackrel{\cdot g^{1-n}}{\longrightarrow}\cdots\stackrel{\cdot g^{-1}}{\longrightarrow}\hom_\Lambda(Q,B^0)\]is exact.}

\vskip.5em{\sc Proof }
(1)(i) Let $X$ be a complex $\cdots\stackrel{c^{i-1}}{\longrightarrow}X^i\stackrel{c^{i}}{\longrightarrow}X^{i+1}\stackrel{c^{i+1}}{\longrightarrow}\cdots$ in $\kk^{\b}(\pr\Lambda)$ with $c^i\in J_{\pr\Lambda}$ for any $i$ such that $X\oplus Q$ is a tilting complex. We have $\hom_{\kk^{\b}(\pr\Lambda)}(Q,X[i])=H^i(\Rhom_\Lambda(Q,X))$ and $\hom_{\kk^{\b}(\pr\Lambda)}(X,Q[i])=H^i(\Rhom_\Lambda(X,Q))$.

(ii) Put $m:=\min\{i\ |\ X^i\neq0\}$. We will show that either $X^m\in\add P$ or $m=0$.

If $X^m\notin\add P$, then we can choose $f^m\in\hom_\Lambda(X^m,Q)$ not in $J_{\pr\Lambda}$. We extend $f^m$ to a chain morphism $f\in\hom_{\kk^{\b}(\pr\Lambda)}(X,Q[-m])$. If $m\neq0$, then we have $f=0$. Thus $f^m$ factors through $c^m\in J_{\pr\Lambda}$, a contradiction.

(iii) Put $n:=\max\{i\ |\ X^i\neq0\}$. Then the dual argument to (ii) shows that either $X^n\in\add P$ or $n=0$.

(iv) Since $\Lambda$ is ring-indecomposable, so is $\endm_{\kk^{\b}(\pr\Lambda)}(X\oplus Q)$ since $X\oplus Q$ is a tilting complex. Thus $m\le 0\le n$ holds. If $m=n=0$, then we have $X=P$. Otherwise, (ii) and (iii) imply that either ($m<0$ and $X^m\in\add P$) or ($0<n$ and $X^n\in\add P$).

(v) We will show that, if $0<n$, then $X$ is isomorphic to ${\rm RA}_n$.

Inductively, we will show that $X^i\in\add Q$ for any $i\neq n$. This is true for any $i$ with $i<m$. Assume that $X^i\in\add Q$ holds for any $i$ with $i<l$.
Since $H^i(\Rhom_\Lambda(Q,X))=0$ for any $i\neq0$ by (i), any $f^l\in\hom_\Lambda(X^l,X^n)$ can be extended to the following chain morphism $f\in\hom_{\kk^{\b}(\pr\Lambda)}(X,X[n-l])$.
{\small\[\begin{diag}
0&\longrightarrow&\cdots&\longrightarrow&0&\longrightarrow&X^m&\longrightarrow&\cdots&\longrightarrow&X^{l-1}&\longrightarrow&X^{l}&\stackrel{c^l}{\longrightarrow}&X^{l+1}&\longrightarrow&\cdots&\longrightarrow&X^n\\
\downarrow&&&&\downarrow&&\downarrow&&&&\downarrow&&\downarrow^{f^l}&&\downarrow&&&&\downarrow\\
X^m&\longrightarrow&\cdots&\longrightarrow&X^{n-l+m-1}&\longrightarrow&X^{n-l+m}&\longrightarrow&\cdots&\longrightarrow&X^{n-1}&\stackrel{c^{n-1}}{\longrightarrow}&X^n&\longrightarrow&0&\longrightarrow&\cdots&\longrightarrow&0
\end{diag}\]}
Since $\hom_{\kk^{\b}(\pr\Lambda)}(X,X[n-l])=0$, we obtain
{\small\[f^l\in c^l\hom_\Lambda(X^{l+1},X^n)+\hom_\Lambda(X^l,X^{n-1})c^{n-1}\subseteq J_{\pr\Lambda}.\]}
This implies $X^l\in\add Q$. 

Thus we have proved $X^i\in\add Q$ for any $i\neq n$. Hence (ii) implies $m=0$, and $X$ gives the first $n$ terms of the minimal right $(\add Q)$-approximation sequence of $X^n\in\add P$. Since $X$ is indecomposable, we have $X^n=P$ and $X={\rm RA}_{n}$.

(vi) A dual argument to (v) implies that, if $m<0$, then $X$ is isomorphic to ${\rm LA}_{-m}$.

(2) The `only if' part follows by (1)(i) above. We will show the `if' part. By (1)(i) again, $\hom_{\kk^{\b}(\pr\Lambda)}(Q,{\rm RA}_n[i])=0=\hom_{\kk^{\b}(\pr\Lambda)}({\rm RA}_n,Q[i])$ holds for any $i\neq0$.

Take any $f\in\hom_{\kk^{\b}(\pr\Lambda)}({\rm RA}_n,{\rm RA}_n[i])$. If $i>0$, then the conditions $A^i\in\add Q$ ($i\neq n$) and $H^i(\Rhom_\Lambda(Q,{\rm RA}_n))=0$ ($i\neq0$) imply that there exist $s^{n-i}, s^{n-i-1}, \cdots, s^0$ in the diagram below such that $f^j=a^js^{j+1}+s^ja^{j-1}$ for any $j$.
{\small\[\begin{diag}
0&\longrightarrow&\cdots&\longrightarrow&0&\RA{}&A^0&\RA{a^0}&\cdots\cdots&\RA{a^{n-i-2}}&A^{n-i-1}&\RA{a^{n-i-1}}&A^{n-i}&\longrightarrow&A^{n-i+1}&\longrightarrow&\cdots&\longrightarrow&A^n\\
\DA{}&&&&\DA{}&\DLA{s^0}&\DA{f^0}&\DLA{s^1}&&\DLA{s^{n-i-1}}&\DA{f^{n-i-1}}&\DLA{s^{n-i}}&\DA{f^{n-i}}&&\DA{}&&&&\DA{}\\
A^0&\longrightarrow&\cdots&\longrightarrow&A^{i-1}&\RA{a^{i-1}}&A^i&\RA{a^i}&\cdots\cdots&\RA{a^{n-2}}&A^{n-1}&\RA{a^{n-1}}&A^n&\longrightarrow&0&\longrightarrow&\cdots&\longrightarrow&0
\end{diag}\]}
Thus we have $f=0$. On the other hand, if $i<0$, then the conditions $A^i\in\add Q$ ($i\neq n$) and $H^i(\Rhom_\Lambda({\rm RA}_n,Q))=0$ ($i\neq0$) imply that there exist $s^{1-i}, s^{2-i}, \cdots, s^n$ in the diagram below such that $f^j=a^js^{j+1}+s^ja^{j-1}$ for any $j$.
{\small\[\begin{diag}
A^0&\longrightarrow&\cdots&\longrightarrow&A^{-i-1}&\longrightarrow&A^{-i}&\RA{a^{-i}}&A^{1-i}&\RA{a^{1-i}}&\cdots\cdots&\RA{a^{n-1}}&A^n&\longrightarrow&0&\longrightarrow&\cdots&\longrightarrow&0\\
\DA{}&&&&\DA{}&&\DA{f^{-i}}&\DLA{s^{1-i}}&\DA{f^{1-i}}&\DLA{s^{2-i}}&&\DLA{s^{n}}&\DA{f^{n}}&&\DA{}&&&&\DA{}\\
0&\longrightarrow&\cdots&\longrightarrow&0&\longrightarrow&A^0&\RA{a^0}&A^{1}&\RA{a^{1}}&\cdots\cdots&\RA{a^{n+i-2}}&A^{n+i}&\longrightarrow&A^{n+i+1}&\longrightarrow&\cdots&\longrightarrow&A^n
\end{diag}\]}
Thus we have $f=0$.

Consequently, $\hom_{\kk^{\b}(\pr\Lambda)}({\rm RA}_n\oplus Q,({\rm RA}_n\oplus Q)[i])=0$ holds for any $i\neq0$. Thus ${\rm RA}_n\oplus Q$ is a tilting complex since it clearly generates $\kk^{\b}(\pr\Lambda)$. We can show (3) dually.\rule{5pt}{10pt}

\vskip.5em
As a consequence, we obtain information on the number of possible complements.

\vskip.5em{\bf Corollary \XDB\ }{\it
For any $n>0$, $Q$ has at most $2n-1$ complements giving rise to tilting complexes with term length at most $n$ in $\kk^{\b}(\pr\Lambda)$. For any $n>0$, $Q$ has at most $n$ complements giving rise to tilting modules of projective dimension at most $n-1$.}

\vskip.5em{\sc Proof }
By \XDA, ${\rm LA}_i,{\rm RA}_i$ ($0\le i\le n$) are the possible complements with term length at most $n$ in $\kk^{\b}(\pr\Lambda)$. Note that ${\rm RA}_0=P={\rm LA}_0$. Thus the first assertion follows. Since $H^i({\rm RA}_i)$ never vanishes for any $i\neq0$, ${\rm RA}_i$ can never be isomorphic to a module. Thus the second assertion follows.\rule{5pt}{10pt}

\vskip.5em
We now give a basic result on $n$-CY algebras $\Lambda$, which we use to obtain more precise information on the number of complements. Following Seidel-Thomas [ST], we say that a simple $\Lambda$-module $S$ is {\it $n$-spherical} for $n>0$ if $\ext^i_\Lambda(S,S)=0$ for any $i$ with $i\neq0,n$. (See section \XF\ for a more general definition.)

\vskip.5em{\bf Proposition \XDC\ }{\it
%Let $\Lambda$ be a $d$-CY algebra, $S$ a simple $\Lambda$-module, $P$ a projective cover of $S$, and $Q$ a direct sum of all indecomposable projective $\Lambda$-modules which are not isomorphic to $P$.
Let $\Lambda$ be a basic $d$-CY algebra, $e$ a primitive idempotent of $\Lambda$, $P:=\Lambda e$ and $Q:=\Lambda(1-e)$. Take a minimal projective resolution
$0\to P_d\stackrel{f_d}{\to}P_{d-1}\stackrel{f_{d-1}}{\to}\cdots\stackrel{f_2}{\to}P_1\stackrel{f_1}{\to}P_0\to S\to0$ of $S:=P/J_\Lambda P$.

(1) We have an exact sequence
$0\to\hom_\Lambda(P_0,\Lambda)\stackrel{f_1\bullet}{\to}\cdots\stackrel{f_d\bullet}{\to}\hom_\Lambda(P_d,\Lambda)\to DS\to0$.

(2) $P_d\simeq P_0=P$.

(3) $\Omega^iS$ is indecomposable for any $i$ ($0\le i\le d$).

(4) $S$ is $d$-spherical if and only if $\bigoplus_{i=1}^{d-1}P_i\in\add Q$.

(5) If $S$ is $d$-spherical, then the induced morphism $\Omega^iS\to P_{i-1}$ by $f_i$ is a minimal left $(\add Q)$-approximation for any $i$ ($1<i\le d$).}

\vskip.5em{\sc Proof }
(1) We have an isomorphism
{\small\[\ext^i_\Lambda(S,\Lambda)\simeq D\ext^{d-i}_\Lambda(\Lambda,S)=\left\{\begin{array}{cc}0&(0\le i<d)\\
DS&(i=d)
\end{array}\right.\]}
Thus we have the desired exact sequence by applying $\hom_\Lambda(-,\Lambda)$.

(2) Since the projective cover of $DS$ is $\hom_\Lambda(P_0,\Lambda)$, the assertion follows by (1).

(3) This follows by using that $P=P_d$ is indecomposable.

(4) This is a direct consequence of the definition of $S$ being $d$-spherical.

(5) Applying $-\otimes_\Lambda Q$ to the exact sequence in (1), we get an exact sequence
$0\to\hom_\Lambda(P_0,Q)\stackrel{f_1\bullet}{\to}\cdots\stackrel{f_n\bullet}{\to}\hom_\Lambda(P_d,Q)\to0$ using that $(DS)\otimes_\Lambda Q=0$. Thus the assertion follows.\rule{5pt}{10pt}

\vskip.5em
We now obtain our sufficient conditions for an almost complete projective module to have exactly two completions giving rise to tilting modules of projective dimension at most one.

\vskip.5em{\bf Theorem \XDD\ }{\it
Let $\Lambda$ be a basic $d$-CY algebra, $e$ a primitive idempotent of $\Lambda$, $P:=\Lambda e$ and $Q:=\Lambda(1-e)$. Assume that $S:=P/J_\Lambda P$ is $d$-spherical.
%$P$ an indecomposable projective $\Lambda$-module and $Q$ a direct sum of the indecomposable projective $\Lambda$-modules which are not isomorphic to $P$. Assume that $S:=P/J_\Lambda P$ is $d$-spherical.

(1) Any ${\rm RA}_n$ and ${\rm LA}_n$ ($n\ge0$) are complements of $Q$.

(2) $Q$ has exactly $d$ complements $\Omega^{n}S$ ($1\le n\le d$) giving rise to tilting modules of finite projective dimension. They satisfy $\pd{}_\Lambda(\Omega^{n}S)=d-n$.

(3) $Q$ has exactly 2 complements $P$ and $\Omega^{d-1}S$ giving rise to tilting modules of projective dimension at most one. They are reflexive if $d\ge3$.}

\vskip.5em{\sc Proof }
(1) We use the notation in \XDC. Since $S$ is $d$-spherical, ${\rm RA}_n$ is obtained from a minimal projective resolution of $S$ as follows:
{\small\[\cdots\to P_1\to P_{d-1}\to\cdots\to P_1\to P_{d-1}\to\cdots\to P_1\to P_0\to0\to\cdots\]}
By \XDC, $H^i(\Rhom_\Lambda({\rm RA}_n,Q))=0$ holds for any $i\neq0$. Thus ${\rm RA}_n$ is a complement of $Q$. A similar argument works for ${\rm LA}_n$.

(2) $\Omega^nS$ is quasi-isomorphic to the complex $X=(\cdots\to0\to P_d\to P_{d-1}\to\cdots\to P_n\to0\to\cdots)$, which satisfies $\bigoplus_{i=n}^{d-1}P_i\in\add Q$ since $S$ is $d$-spherical. Since $\hom_\Lambda(P_n,Q)\to\cdots\to\hom_\Lambda(P_{d-1},Q)\to\hom_\Lambda(P_d,Q)\to0$ is exact by \XDC(5) and $0\to\hom_\Lambda(Q,P_d)\to\hom_\Lambda(Q,P_{d-1})\to\cdots\to\hom_\Lambda(Q,P_n)$ is also exact, we have that $X={\rm LA}_{d-n}$ is a complement of $Q$ by \XDA.

(3) Since $\Lambda$ is reflexive and $\ref\Lambda$ is closed under kernels, $\Omega^{d-1}S$ is reflexive if $d\ge3$.\rule{5pt}{10pt}

%\vskip1.5em
\newpage{\bf\XE. Mutation on tilting modules}

Let $R$ be a normal complete local Gorenstein domain and $\Lambda$ a module-finite $R$-algebra. Throughout this section, all (almost complete) tilting modules have projective dimension at most one. We denote by $\tilt_1\Lambda$ the set of isomorphism classes of basic tilting $\Lambda$-modules. We have seen in section \XD\ that for $d$-CY algebras with no loops in their quiver, the almost complete projective tilting modules have exactly two complements, and we have given an explicit description of the complements. We improve these results, by dropping the assumption that the almost complete tilting module is projective, and show that a more general class of $d$-CY algebras have the same property, including all 2-CY and 3-CY algebras. We also give a description of the complements. 

We start with some background material on tilting modules. This is taken from the theory of finite dimensional algebras [RS][HU1][U], but is stated in our more general context of module-finite $R$-algebras, where the results remain valid. Since the proofs are usually the same as for finite dimensional algebras, they are mostly omitted. Let us start with the following.
%The following proposition shows existence of at least one complements, which we call {\it Bongartz complement}.
%\vskip.5em{\bf Lemma \XEB\ }{\it Assume that $T\in\mod\Lambda$ satisfies $\ext^1_\Lambda(T,T)=0$ and $\pd{}_\Lambda T\le 1$. For any $C\in\mod\Lambda$, there exists an exact sequence $0\to C\to X\to T'\to0$ with $T'\in\add T$ and $\ext^1_\Lambda(T,X)=0$.} \vskip.5em{\sc Proof } Let $P\stackrel{f}{\to}\ext^1_\Lambda(T,C)\to0$ be exact, with $P$ projective in $\mod\endm_\Lambda(T)$. We can write $P=\hom_\Lambda(T,T^\prime)$ for $T'\in\add T$. It follows from Yoneda's lemma on $\add T$ that $f$ is given by $\sigma\in\ext^1_\Lambda(T',C)$. Take an exact sequence $0\to C\to X\to T^\prime\to0$ corresponding to $\sigma$. Then $\hom_\Lambda(T,T^\prime)\stackrel{\bullet\sigma}{\to}\ext^1_\Lambda(T,C)\to0$ is exact. Applying $\hom_\Lambda(T,-)$, we see that $\ext^1_\Lambda(T,X)=0$.\rule{5pt}{10pt}

\vskip.5em{\bf Proposition \XEA\ }{\it
Any almost complete tilting module $T$ has at least one complement (called a Bongartz complement constructed in \XBH) and at most two complements.}

\vskip.5em{\sc Proof }
This follows from \XBH\ and a similar argument as in [RS;1.3][U].\rule{5pt}{10pt}

\vskip.5em
To study the relationship between two complements of an almost complete tilting module, let us recall the following result [RS;1.3].

\vskip.5em{\bf Proposition \XEB\ }{\it
Let $T$ be an almost complete tilting $\Lambda$-module and $0\to Y\stackrel{g}{\to}T^\prime\stackrel{f}{\to}X\to0$ an exact sequence with $T'\in\add T$. Then the following conditions are equivalent.

(1) $X$ is a complement of $T$ and $f$ is a minimal right $(\add T)$-approximation.

(2) $Y$ is a complement of $T$, $g$ is a minimal left $(\add T)$-approximation, and $\pd{}_\Lambda X\le 1$.}

\vskip.5em{\sc Proof }
(1)$\Rightarrow$(2) Applying $\hom_\Lambda(T,-)$, we obtain $\ext^1_\Lambda(T,Y)=0$. Applying $\hom_\Lambda(-,T\oplus X)$, we obtain $\ext^1_\Lambda(Y,T\oplus X)=0$. We will show that $\hom_\Lambda(Y,T^\prime)\stackrel{\bullet f}{\to}\hom_\Lambda(Y,X)\to0$ is exact. Then we have $\ext^1_\Lambda(Y,Y)=0$ by applying $\hom_\Lambda(Y,-)$. Fix any $a\in\hom_\Lambda(Y,X)$. Since $\ext^1_\Lambda(X,X)=0$, there exists $b$ such that $a=gb$. Since $f$ is a right $(\add T)$-approximation of $X$, there exists $c$ such that $b=cf$. Thus $a=(gc)f$ holds.

Since $T\oplus X$ generates $\kk^{\b}(\pr\Lambda)$, it follows from the exact sequence $0\to Y\stackrel{}{\to}T^\prime\stackrel{}{\to}X\to0$ that $T\oplus Y$ also generates $\kk^{\b}(\pr\Lambda)$. Thus $T\oplus Y$ is a tilting $\Lambda$-module. It follows from $\ext^1_\Lambda(X,T)=0$ that $g$ is a left $(\add T)$-approximation. 

One can show (2)$\Rightarrow$(1) similarly.\rule{5pt}{10pt}

\vskip.5em
When the conditions of Proposition \XEB\ hold, put
\[\nu^-_X(T\oplus X):=T\oplus Y\ \ \ \mbox{ and }\ \ \ \nu^+_Y(T\oplus Y):=T\oplus X.\]
We call these operations {\it mutations}. For example, in \XDD(2), we have $\nu^-_{\Omega^{d-1}S}(Q\oplus\Omega^{d-1}S)=\Lambda$ and $\nu^+_{P}(\Lambda)=Q\oplus\Omega^{d-1}S$. For any basic tilting $\Lambda$-module $T$ and any indecomposable direct summand $X$ of $T$, at most one of $\nu^-_X(T)$ and $\nu^+_X(T)$ exists by \XEA, and we sometimes denote it by
\[\nu_X(T).\]
We put $T^\perp:=\{C\in\mod\Lambda\ |\ \ext^1_\Lambda(T,C)=0\}$. Following [RS] (see also [HU2]), we write
\[T\le U\]
if $T^\perp\supseteq U^\perp$. Then $\tilt_1\Lambda$ forms a partially ordered set with a unique minimal element $\Lambda$. One can easily check that, if $\nu^-_X(T)$ (resp. $\nu^+_X(T)$) exists, then $\nu^-_X(T)<T$ (resp. $T<\nu^+_X(T)$). Recall that the {\it Hasse quiver} of $\tilt_1\Lambda$ is the quiver with the set of vertices $\tilt_1\Lambda$, and we draw an arrow $T\to U$ ($T,U\in\tilt_1\Lambda$) if $T<U$ and there is no $V\in\tilt_1\Lambda$ such that $T<V<U$. The following proposition asserts that the arrows of the Hasse quiver of $\tilt_1\Lambda$ are given by mutation. 

\vskip.5em{\bf Proposition \XEC\ }{\it
%(5) If $T<U$, then there exists an exact sequence $0\to T_1\stackrel{}{\to}T_0\stackrel{}{\to}U\to0$ with $\add T=\add T_0\oplus T_1$ and $\add T_0\cap\add T_1=0$.
%(4) If $\nu^-_X(U)$ (resp. $\nu^+_X(U)$) exists, then $\nu^-_X(U)<U$ (resp. $U<\nu^+_X(U)$).
(1) For $T,U\in\tilt_1\Lambda$, the following conditions are equivalent.

\strut\kern1em
(i) $T<U$.

\strut\kern1em
(ii) There exists an indecomposable direct summand $X$ of $U$ such that $T\le\nu^-_X(U)$.

\strut\kern1em
(iii) There exists an indecomposable direct summand $Y$ of $T$ such that $\nu^+_Y(T)\le U$.

(2) For $T,U\in\tilt_1\Lambda$, the following conditions are equivalent.

\strut\kern1em
(i) There exists an arrow $T\to U$ in the Hasse quiver of $\tilt_1\Lambda$.

\strut\kern1em
(ii) There exists an indecomposable direct summand $X$ of $U$ such that $T=\nu^-_X(U)$.

\strut\kern1em
(iii) There exists an indecomposable direct summand $Y$ of $T$ such that $U=\nu^+_Y(T)$.

\strut\kern1em
(iv) There exists an almost complete tilting $\Lambda$-module which is a common direct summand of $T$ and $U$.}

\vskip.5em
The next result generalizes \XDC\ and \XDD\ in two directions. For one thing, we treat arbitrary tilting modules which are not necessarily projective. In addition, we drop the assumption in \XDD\ that $S$ is $d$-spherical, and we replace it by a weaker assumption on the depth and injective dimension of $\Gamma/I$. Notice that we can obtain \XDC\ and \XDD\ by putting $T=\Lambda$ in \XED.

%Now we consider the case when a simple module $S$ over an $d$-CY algebra $\Lambda$ is not $d$-spherical. We will show that, for the case $d=2$ or $3$, we can still construct a tilting $\Lambda$-module by taking an idempotent ideal $I$ and considering a minimal projective resolution of $\Lambda/I$. Our previous theorem \XDD is a special case of this construction. In the rest of this section, assume that $R$ is a normal domain. We start with a preliminary results valid for any dimension.

\vskip.5em{\bf Theorem \XED\ }{\it
Let $\Lambda$ be a $d$-CY algebra, $T$ a basic tilting $\Lambda$-module and $\Gamma:=\endm_\Lambda(T)$. For a primitive idempotent $e$ of $\Gamma$, put $P:=\Gamma e$, $Q:=\Gamma(1-e)$, $I:=\Gamma(1-e)\Gamma$ and $S:=\Gamma/(I+J_\Gamma)$. Assume that the equality $n:=\depth(\Gamma/I)=\id{}_{\Gamma/I}(\Gamma/I)$ holds.

(1) There exists a minimal projective resolution $0\to P_{d-n}\stackrel{f_{d-n}}{\to}\cdots\stackrel{f_1}{\to}P_0\stackrel{f_0}{\to}\Gamma/I\to0$ of the $\Gamma$-module $\Gamma/I$.

(2) We have a minimal projective resolution $0\to\hom_\Gamma(P_0,\Gamma)\stackrel{f_1\bullet}{\longrightarrow}\cdots\stackrel{f_{d-n}\bullet}{\longrightarrow}\hom_\Gamma(P_{d-n},\Gamma)\to\Gamma/I\to0$ of the $\Gamma^{\op}$-module $\Gamma/I$.

(3) $P_{d-n}\simeq P_0=P$ and $P_1,P_{d-n-1}\in\add Q$.

(4) Fix $i\ge0$. If $\tor^\Gamma_i(T,S)=0$, then $\tor^\Gamma_i(T,X)=0$ for any $X\in\mod(\Gamma/I)$.

(5) $\nu^+_{\Gamma e}(\Gamma)=\tr_{\Gamma^{\op}}(\Gamma/I)\oplus\Gamma(1-e)$ and $\nu^+_{e\Gamma}(\Gamma)=\tr_\Gamma(\Gamma/I)\oplus(1-e)\Gamma$.

(6) Precisely one of (i) or (ii) holds.

\strut\kern1em
(i) $T\otimes_\Gamma S=0$ and $\nu^-_{Te}(T)=\hom_{\Gamma^{\op}}(\nu^+_{e\Gamma}(\Gamma),T)$.

\strut\kern1em
(ii) $\tor^\Gamma_1(T,S)=0$ and $\nu^+_{Te}(T)=T\otimes_\Gamma\nu^+_{\Gamma e}(\Gamma)$.

(7) If $d-n\ge3$ and $T$ is reflexive, then $\nu_{Te}(T)$ is reflexive.}

\vskip.5em{\sc Proof }
(1) $\pd{}_\Gamma(\Gamma/I)=d-n$ holds by \XBC.

(2)(3) Since $\Gamma/I$ is a CM $\Gamma$-module of dimension $n$, it follows from \XCD(5) that $\ext^i_\Gamma(\Gamma/I,\Gamma)=0$ for any $i\neq d-n$ and $\ext^{d-n}_\Gamma(\Gamma/I,\Gamma)\simeq\ext^{d-n}_R(\Gamma/I,R)$ is a CM $\Gamma^{\op}$-module of dimension $n$. Since $\id{}_{\Gamma/I}(\Gamma/I)=n$, we have that $\ext^{d-n}_R(\Gamma/I,R)$ is a projective $(\Gamma/I)^{\op}$-module by [GN1;1.1(3)]. Since $\Gamma/I$ is local, $\ext^{d-n}_R(\Gamma/I,R)\simeq\Gamma/I$ as a $\Gamma^{\op}$-module. Now we can show (2) and (3) by a similar argument as in the proof of \XDC(1)(2).

(4) Use induction on $\dim X$ similarly as in the proof of \XBB.

(5)(6) Applying \XBI\ to the tilting $\Gamma^{\op}$-module $T$, precisely one of $T\otimes_\Gamma S=0$ or $\tor_1^\Gamma(T,S)=0$ holds. 

(i) Assume $T\otimes_\Gamma S=0$. Since $T\otimes_\Gamma(\Gamma/I)=0$ holds by (4), we have an exact sequence $0\to\hom_\Gamma(\tr_\Gamma(\Gamma/I),T)\to T\otimes_\Gamma P_1\stackrel{T\otimes f_1}{\longrightarrow}Te\to0$ by applying $\hom_\Gamma(-,T)$ to the exact sequence $\hom_\Gamma(P_0,\Gamma)\stackrel{f_1\bullet}{\to}\hom_\Gamma(P_1,\Gamma)\to\tr_\Gamma(\Gamma/I)\to0$. Since $f_1$ is a minimal right $(\add Q)$-approximation, $T\otimes f_1$ is a minimal right $(\add T(1-e))$-approximation. Thus we have $\nu^-_{Te}(T)=\hom_\Gamma(\tr_\Gamma(\Gamma/I),T)\oplus T(1-e)$.

(ii) Put $U:=T\otimes_\Gamma\tr_{\Gamma^{\op}}(\Gamma/I)$.
Since $\tor^\Gamma_i(T,\Gamma/I)=0$ holds for any $i>0$ by (4), we have an exact sequence $0\to T\otimes_\Gamma P_{d-n}\stackrel{T\otimes f_{d-n}}{\longrightarrow}\cdots\to T\otimes_\Gamma P_0\stackrel{T\otimes f_0}{\longrightarrow}T\otimes_\Gamma(\Gamma/I)\to0$. In particular, $0\to Te\stackrel{T\otimes f_{d-n}}{\longrightarrow}T\otimes_\Gamma P_{d-n-1}\to U\to0$ is exact by (2). Since $f_{d-n}$ is a minimal left $(\add Q)$-approximation, $T\otimes f_{d-n}$ is a minimal left $(\add T(1-e))$-approximation. Thus we only have to show $\pd{}_\Lambda U\le 1$, or equivalently, $\depth U\ge d-1$ by \XBC.

Take a $\Gamma/I$-regular sequence $(x_1,\cdots,x_n)$, and put $\overline{\Gamma}_i:=(\Gamma/I)/(x_1,\cdots,x_i)(\Gamma/I)$ for $i=1,\cdots,n$. Then for $i<n$ we have an exact sequence $0\to\overline{\Gamma}_i\stackrel{x_{i+1}}{\longrightarrow}\overline{\Gamma}_i\to\overline{\Gamma}_{i+1}\to0$. Applying $T\otimes_\Gamma-$, we have an exact sequence $0\to T\otimes_\Gamma\overline{\Gamma}_i\stackrel{x_{i+1}}{\longrightarrow}T\otimes_\Gamma\overline{\Gamma}_i\to T\otimes_\Gamma\overline{\Gamma}_{i+1}\to0$ since $\tor^\Gamma_1(T,\overline{\Gamma}_{i+1})=0$ by (4). This means that $(x_1,\cdots,x_n)$ is also a $(T\otimes_\Gamma(\Gamma/I))$-regular sequence. In particular, we have $\depth(T\otimes_\Gamma(\Gamma/I))\ge n$. Since $\depth(T\otimes_\Gamma P_i)\ge d-1$, the exact sequence $0\to U\to T\otimes_\Gamma P_{d-n-2}\to\cdots\to T\otimes_\Gamma P_0\stackrel{}{\to}T\otimes_\Gamma(\Gamma/I)\to0$ implies $\depth U\ge d-1$.

Putting $\Lambda=T=\Gamma$ in (ii), we have $\nu^+_{\Gamma e}(\Gamma)=\tr_{\Gamma^{\op}}(\Gamma/I)\oplus\Gamma(1-e)$ and $\nu^+_{e\Gamma}(\Gamma)=\tr_{\Gamma}(\Gamma/I)\oplus(1-e)\Gamma$. Thus the equalities in (5) and (6) follow.

(7) This is obvious for the case (6)(i). For the case (6)(ii), the assertion follows from the exact sequence $0\to U\to T\otimes_\Gamma P_{d-n-2}\to T\otimes_\Gamma P_{d-n-3}$.\rule{5pt}{10pt}

\vskip.5em
We notice here that we can regard the tilting modules constructed in (5) above as analogs of APR tilting modules [APR].

To apply \XED\ for the case $d=2$ and $3$, we need the following observation.

\vskip.5em{\bf Lemma \XEE\ }{\it
Let $\Gamma$ be a ring-indecomposable $d$-CY algebra. For an idempotent $e\neq1$ of $\Gamma$, put $I:=\Gamma(1-e)\Gamma$.

(1) $\Gamma\otimes_RK$ is a simple algebra for the quotient field $K$ of $R$.

(2) $\dim (\Gamma/I)\le\max\{0,d-2\}$.

(3) If $d\le 3$, then $\id_{\Gamma/I}(\Gamma/I)\le\max\{0,d-2\}$.

(4) If $d=3$ and $e$ is primitive, then $\depth(\Gamma/I)=\id_{\Gamma/I}(\Gamma/I)$.}

\vskip.5em{\sc Proof }
(1) Since $R$ is normal, $\Gamma=\bigcap_{\dn{p}}\Gamma_{\dn{p}}$ holds where $\dn{p}$ runs over all height one prime ideals of $R$. By \XCM\ and the structure theorem \XCL(2) of maximal orders over complete discrete valuation rings, $\Gamma_{\dn{p}}$ contains all central idempotents of $\Gamma\otimes_RK$. Thus $\Gamma$ contains all central idemptents of $\Gamma\otimes_RK$. Since $\Gamma$ is ring-indecomposable, $\Gamma\otimes_RK$ is simple.

(2) By the structure theorem of maximal orders, we have $(\Gamma/I)_{\dn{p}}=\Gamma_{\dn{p}}/\Gamma_{\dn{p}}e\Gamma_{\dn{p}}=0$ for any height one prime ideal $\dn{p}$ of $R$. Thus $\dim (\Gamma/I)\le d-2$ holds.

(3) If $d\le1$, then $\Gamma/I=0$ by \XCD(3), \XCL(2) and \XCM. Assume $d=2$ or $3$. Since $\mod(\Gamma/I)$ is extension closed in $\mod\Gamma$, we have $\ext^1_\Gamma(\Gamma/I,X)=0$ for any $X\in\mod(\Gamma/I)$. Since $\Gamma$ is $d$-CY, we have $\ext^{d-1}_\Gamma(X,\Gamma/I)=0$ for any $X\in\flmod(\Gamma/I)$. Assume $d=2$. Using again that $\mod(\Gamma/I)$ is extension closed in $\mod\Gamma$, we get $\ext^1_{\Gamma/I}(X,\Gamma/I)=0$. Thus $\id_{\Gamma/I}(\Gamma/I)=0$.

In the rest, assume $d=3$. For any $X\in\flmod(\Gamma/I)$, take an exact sequence $0\to Y\to(\Gamma/I)^n\to X\to0$. Applying $\hom_\Gamma(-,\Gamma/I)$, we get an exact sequence $0=\ext^1_\Gamma((\Gamma/I)^n,\Gamma/I)\to\ext^1_\Gamma(Y,\Gamma/I)\to\ext^2_\Gamma(X,\Gamma/I)=0$.
Thus we have $\ext^1_\Gamma(Y,\Gamma/I)=0$. Since $\mod(\Gamma/I)$ is extension closed in $\mod\Gamma$, we have $\ext^1_{\Gamma/I}(Y,\Gamma/I)=0$. Thus $\ext^2_{\Gamma/I}(X,\Gamma/I)=0$. It follows from \XBA\ that $\id_{\Gamma/I}(\Gamma/I)\le1$.

(4) Since $\Gamma/I$ is a local algebra with $\id_{\Gamma/I}(\Gamma/I)<\infty$ by (3), it follows from a result of Ramras [Ra;2.15] that $\depth(\Gamma/I)=\id_{\Gamma/I}(\Gamma/I)$ (see also [GN2;3.9]).\rule{5pt}{10pt}

\vskip.5em
We now get our desired result for 2-CY and 3-CY algebras.

\vskip.5em{\bf Theorem \XEF\ }{\it
Let $\Lambda$ be a ring-indecomposable $d$-CY algebra with $d=2$ or $3$, $T$ a basic tilting $\Lambda$-module and $\Gamma:=\endm_\Lambda(T)$. For a primitive idempotent $e\neq1$ of $\Gamma$, put $I:=\Gamma(1-e)\Gamma$ and $S:=\Gamma/(I+J_\Gamma)$. 

(1) $n:=\depth(\Gamma/I)=\id_{\Gamma/I}(\Gamma/I)$, and $n=d-2$ or $d-3$.

(2) $\nu^+_{\Gamma e}(\Gamma)=\tr_{\Gamma^{\op}}(\Gamma/I)\oplus\Gamma(1-e)$ and $\nu^+_{e\Gamma}(\Gamma)=\tr_\Gamma(\Gamma/I)\oplus(1-e)\Gamma$.

(3) $\nu^+_{\Gamma e}(\Gamma)=I=\nu^+_{e\Gamma}(\Gamma)$ if $n=d-2$, and $\nu^+_{\Gamma e}(\Gamma)=\hom_{\Gamma^{\op}}(\nu^+_{e\Gamma}(\Gamma),\Gamma)$ if $n=d-3$.

(4) Precisely one of (i) or (ii) holds.

\strut\kern1em
(i) $T\otimes_\Gamma S=0$ and $\nu^-_{Te}(T)=\hom_{\Gamma^{\op}}(\nu^+_{e\Gamma}(\Gamma),T)$.

\strut\kern1em
(ii) $\tor^\Gamma_1(T,S)=0$ and $\nu^+_{Te}(T)=T\otimes_\Gamma\nu^+_{\Gamma e}(\Gamma)$.

(5) If $n=d-3$ and $T$ is reflexive, then $\nu_{Te}(T)$ is reflexive and $\nu_{Te}(T)=(T\otimes_\Gamma\nu^+_{\Gamma e}(\Gamma))^{**}$.

(6) We have isomorphisms $\endm_\Gamma(\nu^+_{\Gamma e}(\Gamma))\simeq\endm_\Lambda(\nu_{Te}(T))\simeq\endm_{\Gamma^{\op}}(\nu^+_{e\Gamma}(\Gamma))^{\op}$. This is isomorphic to $\Gamma$ if $n=d-2$.}

\vskip.5em{\sc Proof }
(1) is shown in \XEE, and (2) and (4) are shown in \XED. One can check (3) easily by using the exact sequences in \XED(1)(2).

(5) $\nu_{Te}(T)$ is reflexive by \XED(7). The assertion for the case (4)(ii) is obvious. For the case (4)(i), the assertion follows from $\nu_{Te}(T)=(T\otimes_\Gamma\hom_{\Gamma^{\op}}(\nu^+_{e\Gamma}(\Gamma),\Gamma))^{**}$ and (3).

(6) We have ring morphisms $a:=\hom_{\Gamma^{\op}}(-,T)_{\nu^+_{e\Gamma}(\Gamma),\nu^+_{e\Gamma}(\Gamma)}:\endm_{\Gamma^{\op}}(\nu^+_{e\Gamma}(\Gamma))^{\op}\to\endm_\Lambda(\nu^-_{Te}(T))$ for (4)(i) and $b:=(T\otimes_\Gamma-)_{\nu^+_{\Gamma e}(\Gamma),\nu^+_{\Gamma e}(\Gamma)}:\endm_\Gamma(\nu^+_{\Gamma e}(\Gamma))\to\endm_\Lambda(\nu^+_{Te}(T))$ for (4)(ii) between $d$-CY algebras. For any $\dn{p}\in\Spec R$ with $\height\dn{p}=1$, then $\nu^+_{\Gamma e}(\Gamma)_{\dn{p}}$, $\nu^+_{e\Gamma}(\Gamma)_{\dn{p}}$ and $T_{\dn{p}}$ are progenerators. Thus $a_{\dn{p}}$ and $b_{\dn{p}}$ are isomorphisms. Since $a$ and $b$ are morphisms between reflexive $R$-modules, they are isomorphisms.

If $n=d-2$, then we have ring morphisms $\Gamma\to\endm_\Gamma(I)=\endm_\Gamma(\nu^+_{\Gamma e}(\Gamma))$ and $\Gamma\to\endm_{\Gamma^{\op}}(I)^{\op}=\endm_{\Gamma^{\op}}(\nu^+_{e\Gamma}(\Gamma))^{\op}$ by (3), and $I_{\dn{p}}=\Gamma_{\dn{p}}$ for any $\dn{p}\in\Spec R$ with $\height\dn{p}=1$. If $n=d-3$, then we have a ring morphism $\hom_{\Gamma^{\op}}(-,\Gamma)_{\nu^+_{e\Gamma}(\Gamma),\nu^+_{e\Gamma}(\Gamma)}:\endm_{\Gamma^{\op}}(\nu^+_{e\Gamma}(\Gamma))^{\op}\to\endm_\Gamma(\nu^+_{\Gamma e}(\Gamma))$ by (3). We can show that these are isomorphisms by a similar argument as above.\rule{5pt}{10pt}

%Let $I:=\Lambda e\Lambda$ and $S$ a simple $\Lambda/I$-module. Since $\mod\Lambda/I$ is extension closed in $\mod\Lambda$, we have $\ext^1_\Lambda(\Lambda/I,X)=0$ for any $X\in\mod\Lambda/I$. Since $\Lambda$ is 3-CY, we have $\ext^2_\Lambda(X,\Lambda/I)=0$. Take an exact sequence \[0\to Y\to(\Lambda/I)^n\to X\to0.\] Taking $\hom_\Lambda(-,\Lambda/I)$, we have an exact sequence \[0=\ext^1_\Lambda((\Lambda/I)^n,\Lambda/I)\to\ext^1_\Lambda(Y,\Lambda/I)\to\ext^2_\Lambda(X,\Lambda/I,\Lambda/I)=0.\] Thus we have $\ext^1_\Lambda(Y,\Lambda/I)=0$. Since $\mod\Lambda/I$ is extension closed in $\mod\Lambda$ again, we have $\ext^1_{\Lambda/I}(Y,\Lambda/I)=0$. Thus $\ext^2_{\Lambda/I}(X,\Lambda/I)=0$ holds.\rule{5pt}{10pt}
\vskip.5em
The following generalization of \XDD(3) follows immediately from \XEF(4)(5).

\vskip.5em{\bf Corollary \XEG\ }{\it
Let $\Lambda$ be a basic ring-indecomposable non-local $d$-CY algebra with $d=2$ or $3$, and let $n$ be the number of simple $\Lambda$-modules.

(1) Any almost complete tilting $\Lambda$-module has exactly two complements. Thus any vertex in the Hasse quiver of $\tilt_1\Lambda$ has precisely $n$ neighbours.

(2) Assume $d=3$. If $T$ is a reflexive tilting $\Lambda$-module and the quiver of $\endm_\Lambda(T)$ has no loops, then all $n$ neighbours of $T$ are again reflexive.}

\vskip.5em
We have the following natural questions.

\vskip.5em{\bf Questions }
(1) Is the Hasse quiver of $\tilt_1\Lambda$ connected? In other words, can any tilting $\Lambda$-module be obtained by applying successive mutations to $\Lambda$? We will show in the next section that this is the case for $d=2$.

(2) Is \XEG\ valid for arbitrary $d$? By \XED, it is enough to show $\depth(\Gamma/I)=\id{}_{\Gamma/I}(\Gamma/I)$ for any ring-indecomposable $d$-CY algebra $\Gamma$ and any primitive idempotent $e$ of $\Gamma$ with $I:=\Gamma(1-e)\Gamma$.

\vskip1.5em
{\bf\XF. 2-Calabi-Yau algebras and affine Weyl groups}

Let $R$ be a normal complete local Gorenstein domain with $\dim R=d$ and $\Lambda$ a basic module-finite $R$-algebra which is $d$-CY. Thus $\Lambda$ is a symmetric $R$-order with $\gl\Lambda=d$. If $d=0$ or $1$, then  $\Lambda$ is Morita equivalent to a finite product of local rings by \XCK\ and \XCM, and any tilting $\Lambda$-module is projective by \XBI. So the next question is to determine tilting $\Lambda$-modules for the case $d=2$. Assume for the rest of this section that $d=2$, and we only deal with tilting modules of projective dimension at most one.

In this section we show that each almost complete tilting $\Lambda$-module has exactly two complements, and also give an explicit description of them. The tilting modules are all ideals, and we describe the set $\tilt_1\Lambda$ of tilting modules as a monoid of ideals generated by a finite set of idempotent ideals. We show that this set $\tilt_1\Lambda$ is in bijection with the elements of the affine Weyl group $W$ associated with the quiver of $\Lambda$, which is given by a generalized extended Dynkin diagram. The group $W$ has two natural partial orders, and we show that they coincide with the two orders on $\tilt_1\Lambda$, where one is recalled in section \XE\ and the other one comes from inclusion of ideals.

When $\Lambda$ is local, then any tilting module over $\Lambda$ is projective. In the rest of this section, we assume that $\Lambda$ is non-local.
Let $e_1,\cdots,e_n$ be a complete set of orthogonal primitive idempotents of $\Lambda$. Put $I_i:=\Lambda(1-e_i)\Lambda$. Then $S_i:=\Lambda/(I_i+J_\Lambda)$ is the simple $\Lambda$-module (resp. simple $\Lambda^{\op}$-module) corresponding to $e_i$. Then $I_i$ is maximal amongst left (resp. right) ideals $I$ of $\Lambda$ such that any composition factor of $\Lambda/I$ is $S_i$.
We have shown in \XEF(3) that $I_i$ is a tilting $\Lambda$-module, and moreover if $T$ is a tilitng $\Lambda$-module, then either $\hom_\Lambda(I_i,T)$ or $I_iT$ is another tilting $\Lambda$-module. For later applications, we investigate these modules more carefully.

%For any tilting $\Lambda$-module $T$, we have tilting complexes $I_i\Lotimes_\Lambda T$ and $\Rhom_\Lambda(I_i,T)$. Let us investigate when they become $\Lambda$-modules.

\vskip.5em{\bf Proposition \XFA\ }{\it
Let $T$ be a tilting $\Lambda$-module and $1\le i\le n$. 

(1) We have $\depth T\ge1$ and natural inclusions $I_iT\subseteq T\subseteq\hom_\Lambda(I_i,T)$ such that $\hom_\Lambda(I_i,T)/I_iT\in\flmod\Lambda$.

(2) Presisely one of (i) or (ii) holds.

\strut\kern1em
(i) $S_i\otimes_\Lambda T=0$, $I_iT=T$ and $\nu^-_{Te_i}(T)=\Rhom_\Lambda(I_i,T)=\hom_\Lambda(I_i,T)$.

\strut\kern1em
(ii) $\tor_1^\Lambda(S_i,T)=0$, $\hom_\Lambda(I_i,T)=T$ and $\nu^+_{Te_i}(T)=I_i\Lotimes_\Lambda T=I_i\otimes_\Lambda T=I_iT$.}

\vskip.5em{\sc Proof }
(1) We have $\depth T\ge1$ by \XBC. It follows from \XEE(2) that $\Lambda/I_i$ is artin. Applying $\hom_\Lambda(-,T)$ to the exact sequence $0\to I_i\to\Lambda\to\Lambda/I_i\to0$, we get an exact sequence $0\to\hom_\Lambda(\Lambda/I_i,T)\to T\to\hom_\Lambda(I_i,T)\to\ext^1_\Lambda(\Lambda/I_i,T)\to0$. Since $\ext^i_\Lambda(\Lambda/I_i,T)$ has finite length for any $i$, we have $\hom_\Lambda(\Lambda/I_i,T)=0$ by $\depth T\ge1$. Thus we have inclusions $I_iT\subseteq T\subseteq\hom_\Lambda(I_i,T)$ with $\hom_\Lambda(I_i,T)/I_iT\in\flmod\Lambda$.

(2) By \XEF(3) and (4), we only have to show equalities $I_iT=T$ and $\Rhom_\Lambda(I_i,T)=\hom_\Lambda(I_i,T)$ in (i) and $\hom_\Lambda(I_i,T)=T$ and $I_i\Lotimes_\Lambda T=I_i\otimes_\Lambda T=I_iT$ in (ii).

(i) Since any composition factor of $T/I_iT$ is $S_i$, we have $I_iT=T$ since $S_i\otimes_\Lambda T=0$. Since
{\small\[\ext^l_\Lambda(I_i,T)=\ext^{l+1}_\Lambda(\Lambda/I_i,T)=\left\{\begin{array}{cc}0&(l>1)\\ D\hom_\Lambda(T,\Lambda/I_i)=0&(l=1)
\end{array}\right.\]}
holds, we have $\Rhom_\Lambda(I_i,T)=\hom_\Lambda(I_i,T)$.

(ii) Any composition factor of $\hom_\Lambda(I_i,T)/T$ is $S_i$, so we have $\hom_\Lambda(I_i,T)=T$ since $\tor^\Lambda_1(S_i,T)=0$. Since $\tor^l_\Lambda(I_i,T)=\tor^{l+1}_\Lambda(\Lambda/I_i,T)=0$ holds for any $l\neq0$, we have $I_i\Lotimes_\Lambda T=I_i\otimes_\Lambda T$. Applying $-\otimes_\Lambda T$ to the exact sequence $0\to I_i\stackrel{}{\to}\Lambda\to\Lambda/I_i\to0$, we have an exact sequence $0=\tor^1_\Lambda(\Lambda/I_i,T)\to I_i\otimes_\Lambda T\stackrel{f}{\to}T$. Thus we have $I_i\otimes_\Lambda T\simeq\Im f=I_iT$.\rule{5pt}{10pt}

\vskip.5em
The set of 2-sided ideals of $\Lambda$ forms a monoid by multiplication of ideals. We denote by $\ii(\Lambda)$ the submonoid generated by the ideals $I_1,\cdots,I_n$. Our first main result in this section is that all tilting $\Lambda$-modules are obtained in this way.

\vskip.5em{\bf Theorem \XFB\ }{\it
$\ii(\Lambda)=\tilt_1\Lambda$ and $\ii(\Lambda)=\tilt_1\Lambda^{\op}$.}

\vskip.5em{\sc Proof }
If $T$ is a tilting $\Lambda$-module, then so is $I_iT$ by \XFA(2). Thus $\ii(\Lambda)$ consists of tilting $\Lambda$-modules. We only have to show that any basic tilting $\Lambda$-module is isomorphic to some element of $\ii(\Lambda)$. We will use the functor $(-)^*\simeq\hom_\Lambda(-,\Lambda)$ (\XCD(5)).

(i) For any $T\in\tilt_1\Lambda$, we will show that $T^{**}$ is a projective $\Lambda$-module such that $T$ is a submodule of $T^{**}$ and $T^{**}/T$ has finite length.

By [AB], we have an exact sequence
\[0\to\ext^1_\Lambda(\tr T,\Lambda)\to T\to T^{**}\to\ext^2_\Lambda(\tr T,\Lambda)\to0.\]
For any $\dn{p}\in\Spec R\backslash\Max R$, it follows from \XBG\ that $\widehat{T}_{\dn{p}}$ is a tilting module over $\widehat{\Lambda}_{\dn{p}}$, which is $0$ or $1$-CY by \XCD(3). Thus $\widehat{T}_{\dn{p}}$ (and hence $(\tr T)^{\widehat{\ }}_{\dn{p}}$) is a projective $\widehat{\Lambda}_{\dn{p}}$-module as we remarked previously. Hence $\ext^i_\Lambda(\tr T,\Lambda)^{\widehat{\ }}_{\dn{p}}=0$ holds for $i=1,2$. This implies that $\ext^i_\Lambda(\tr T,\Lambda)$ has finite length for $i=1,2$. Since $\depth T\ge1$ holds by \XBC, we have $\ext^1_\Lambda(\tr T,\Lambda)=0$. Since $T^{**}\in\ref\Lambda$, we have $\depth T^{**}\ge2$. Thus $T^{**}$ is a projective $\Lambda$-module by \XBC.

(ii) Take a simple submodule $S_i$ of $\soc(T^{**}/T)$. Then $\tor^\Lambda_1(S_i,T)=\tor^\Lambda_1(DS_i,T)=D\ext^1_\Lambda(T,S_i)=\ext^1_\Lambda(S_i,T)\neq0$ holds. By \XFA, $T_1:=\hom_\Lambda(I_i,T)$ is again a tilting $\Lambda$-module with $T\subset T_1$ and $T_1/T\in\flmod\Lambda$. Applying $(-)^*$ to the exact sequence $0\to T\to T_1\to T_1/T\to0$, we get $T^*=T_1^*$ since $\depth\Lambda=2$. Consequently, we have inclusions $T\subset T_1\subseteq T^{**}=T_1^{**}$. Repeating this process, we obtain an increasing sequence
\[T=T_0\subset T_1\subset\cdots\subset T_m=T^{**}\]
of tilting $\Lambda$-modules with $T_{k-1}=I_{i_k}T_{k}$ for any $k$.
Since $T^{**}$ is a projective tilting $\Lambda$-module, we have $\add{}_\Lambda T^{**}=\add{}_\Lambda\Lambda$. Thus $\add{}_\Lambda T=\add{}_\Lambda(I_{i_1}\cdots I_{i_m})$.\rule{5pt}{10pt}

\vskip.5em
For any $T\in\tilt_1\Lambda$ and $1\le i\le n$, there exists an arrow $T\stackrel{}{\to}\nu^+_{Te_i}(T)$ or $\nu^-_{Te_i}(T)\stackrel{}{\to}T$ in the Hasse quiver of the poset $\tilt_1\Lambda$ by \XEC(2). We denote these arrows by $\stackrel{i}{\to}$ in the rest of this section. We have the following consequence.

\vskip.5em{\bf Corollary \XFC\ }{\it
The Hasse quiver of $\tilt_1\Lambda$ is connected, and any vertex $T$ has precisely $n$ neighbours $\nu_{Te_i}(T)$ ($1\le i\le n$). Moreover, $\endm_\Lambda(T)=\Lambda$ holds for any $T\in\tilt_1\Lambda$.}

\vskip.5em{\sc Proof }
The second assertion holds by \XEC(2) and \XFA. We will show that $\tilt_1\Lambda$ is connected. Any $T_0\in\tilt_1\Lambda$ can be written as $T_0=I_{a_1}\cdots I_{a_m}$ by \XFB. We can assume that $m$ is minimal. Put $T_i:=I_{a_1}\cdots I_{a_{m-i}}$. Then we have a strictly increasing sequence
$T_0\subset T_1\subset\cdots\subset T_{m-1}\subset T_m=\Lambda$.
By \XFA, there exists a path $\Lambda\stackrel{a_1}{\to}T_{m-1}\stackrel{a_2}{\to}\cdots\stackrel{a_m}{\to}T_0$ in $\tilt_1\Lambda$.
Now the third assertion follows by \XEF(6).\rule{5pt}{10pt}
%(2) We can assume $T\in\ii(\Lambda)$ by \XFB. Let $a:T\to\Lambda$ be the inclusion. Since $T$ is an ideal of $\Lambda$ such that $\Lambda/T$ has finite length, then $a^*:\Lambda\to T^*$ and $a^{**}:T^{**}\to\Lambda$ are bijective. The functor $(-)^{**}$ gives a ring morphism $f:\endm_\Lambda(T)\to\endm_\Lambda(\Lambda)=\Lambda$. Since $T$ is a two-sided ideal of $\Lambda$, $f$ is an isomorphism.\rule{5pt}{10pt}

\vskip.5em
To give an explicit description of $\ii(\Lambda)$, we determine the quiver of $\Lambda$. Following Happel-Preiser-Ringel [HPR1,2], we call the valued graphs below {\it generalized extended Dynkin diagrams}.

(i) An extended Dynkin diagram,

(ii)\ \ \ $\begin{picture}(10,10)\put(7,4){\circle{15}}\end{picture}\bullet\mbox{-----}\bullet\mbox{-----}\bullet\mbox{-----}\bullet\cdots\cdots\bullet\mbox{-----}\bullet\mbox{-----}\bullet\mbox{-----}\bullet\begin{picture}(10,10)\put(3,4){\circle{15}}\end{picture}$

(iii)\ \ \ $\begin{picture}(10,10)\put(7,4){\circle{15}}\end{picture}\bullet\mbox{-----}\bullet\mbox{-----}\bullet\mbox{-----}\bullet\cdots\cdots\bullet\mbox{-----}\bullet\mbox{-----}\bullet\stackrel{(a\ b)}{\mbox{-----}}\bullet$\ \ \ \ \ \ \ \ \ \ $(a,b)=(2,1)$ or $(1,2)$

\vskip.4em
(iv)\ \ \ $\begin{picture}(10,10)\put(7,4){\circle{15}}\end{picture}\bullet\mbox{-----}\bullet\mbox{-----}\bullet\mbox{-----}\bullet\cdots\cdots\bullet\mbox{-----}\bullet\mbox{-----}\bullet\begin{picture}(10,10)\put(-3,3){\line(3,1){22}}\put(-3,3){\line(3,-1){22}}\put(20,7){$\bullet$}\put(20,-7){$\bullet$}\end{picture}$

\vskip.5em
For a generalized extended Dynkin diagram $\Delta$, define a valued quiver called the {\it double} of $\Delta$ as follows: We replace a valued edge $\bullet\stackrel{(a\ b)}{\mbox{-----}}\bullet$ by two valued arrows $\bullet\def\arraystretch{.3}\begin{array}{cc}\stackrel{(a\ b)}{\longrightarrow}\\
\stackrel{(b\ a)}{\longleftarrow}\end{array}\bullet$ of opposite direction. We replace a loop $\begin{picture}(10,10)\put(7,4){\circle{15}}\end{picture}\bullet$ by an arrow from a vertex to itself.

We can describe valued quivers of 2-CY algebras (c.f. [Boc]). It is an interesting question whether all double of generalized extended Dynkin diagrams occur in this way.

\vskip.5em{\bf Proposition \XFD\ }{\it
The valued quiver of any basic ring-indecomposable 2-CY algebra is a double of a generalized extended Dynkin diagram.}

\vskip.5em{\sc Proof }
By \XDC, the quiver of $\Lambda$ is a double of some graph $\Delta$. For each vertex $e_i$ of $\Lambda$, put $d_i:=\rank_R\Lambda e_i$. Then $d_i$ gives a positive additive function on $\Delta$. By [HPR1], $\Delta$ is a generalized extended Dynkin diagram.\rule{5pt}{10pt}

\vskip.5em
We will give an explicit description of $\ii(\Lambda)$ in terms of affine Weyl groups. By \XFD, the quiver of $\Lambda$ is a double of a generalized extended Dynkin diagram $\Delta$. We denote by $W$ the {\it affine Weyl group} associated with $\Delta$ defined as follows [Hu][BB]: Put $m(i,i):=1$. For $i\neq j$, put
\[m(i,j):=\left\{\begin{array}{ll}
2\ \ \ \ \ &\mbox{no edge between $i$ and $j$,}\\
3&\stackrel{i}{\bullet}\stackrel{}{\mbox{------}}\stackrel{j}{\bullet},\\
4&\begin{picture}(10,10)\put(7,4){\circle{15}}\end{picture}\stackrel{i}{\bullet}\stackrel{}{\mbox{------}}\stackrel{j}{\bullet},\ \ 
\stackrel{i}{\bullet}\stackrel{}{\mbox{------}}\stackrel{j}{\bullet}\begin{picture}(10,10)\put(3,4){\circle{15}}\end{picture}\ \mbox{ or }
\stackrel{i}{\bullet}\stackrel{(a\ b)}{\mbox{------}}\stackrel{j}{\bullet}\ ((a\ b)=(1\ 2)\mbox{ or }(2\ 1))\\
6&\stackrel{i}{\bullet}\stackrel{(a\ b)}{\mbox{------}}\stackrel{j}{\bullet}\ ((a\ b)=(1\ 3)\mbox{ or }(3\ 1)),\\
\infty&\stackrel{i}{\bullet}\stackrel{(2\ 2)}{\mbox{------}}\stackrel{j}{\bullet},\ \ 
\begin{picture}(10,10)\put(7,4){\circle{15}}\end{picture}\stackrel{i}{\bullet}\stackrel{}{\mbox{------}}\stackrel{j}{\bullet}\begin{picture}(10,10)\put(3,4){\circle{15}}\end{picture},\ \ \begin{picture}(10,10)\put(7,4){\circle{15}}\end{picture}\stackrel{i}{\bullet}\stackrel{(a\ b)}{\mbox{------}}\stackrel{j}{\bullet}\ \mbox{ or }
\stackrel{i}{\bullet}\stackrel{(a\ b)}{\mbox{------}}\stackrel{j}{\bullet}\begin{picture}(10,10)\put(3,4){\circle{15}}\end{picture}\ ((a\ b)=(1\ 2)\mbox{ or }(2\ 1)).
\end{array}\right.\]
Then $W$ is presented by generators $s_1,\cdots,s_n$ and relations $(s_is_j)^{m(i,j)}=1$. We shall also deal with the affine braid group associated to $\Delta$ to study autoequivalences of the derived category $\dd^{\b}(\mod\Lambda)$. This group $B$ is presented by generators $t_1,\cdots,t_n$ and relations $t_it_jt_i\cdots=t_jt_it_j\cdots$, where both sides are product of $m(i,j)$ generators.

We have seen in section \XE\ that $\tilt_1\Lambda$ has a natural order. In view of \XFB\ we have in addition the order given by inclusion of ideals. The affine Weyl group $W$ also has two partial orders. We want to show that there is a bijection between the elements of $\tilt_1\Lambda$ and $W$, respecting partial orders.

So let us recall the {\it Bruhat order} $\le$, {\it right order} $\le_R$ and {\it left order} $\le_L$ on $W$ [BB] ($\le_R$ is called a {\it weak order} in [Hu]). The {\it length} $l(w)$ of $w\in W$ is the minimal value of $k$ for any expression $w=s_{a_1}\cdots s_{a_k}$ of $w$, and we call an expression with $k=l(w)$ {\it reduced}. Fix $w,w^\prime\in W$. We draw an arrow $w^\prime\to w$ if both $w=w^\prime s$ and $l(w^\prime)<l(w)$ hold for $s=xs_ix^{-1}$ for some $i$ and $x\in W$. Similarly, we draw an arrow $w^\prime\to_Rw$ (resp. $w^\prime\to_Lw$) if both $w=w^\prime s_i$ (resp. $w=s_iw^\prime$) and $l(w^\prime)<l(w)$ hold for some $i$.
%The (weak) Bruhat order is defined so that this quiver gives the Hasse quiver of $W$:
For $w,w^\prime\in W$, we define $w^\prime\le w$ (resp. $w^\prime\le_Rw$, $w^\prime\le_Lw$) if and only if there is a path from $w^\prime$ to $w$ consisting of arrows $\to$ (resp. $\to_R$, $\to_L$).
For any reduced expression $w=s_{a_1}\cdots s_{a_k}$, it is well-known that $w^\prime\le w$ holds if and only if $w^\prime=s_{a_{i_1}}\cdots s_{a_{i_q}}$ for some $1\le i_1<\cdots<i_q\le k$ [Hu;5.10]. Thus the Bruhat order on $W$ is left-right symmetric (i.e. $w^\prime\le w$ if and only if $w^\prime{}^{-1}\le w^{-1}$), but the right order and the left order are not.

We are now in the position to state our main result on the connection between $\tilt_1\Lambda$ and the affine Weyl group $W$. The crucial role is played by the mutation of tilting modules given in \XFA.

\vskip.5em{\bf Theorem \XFE\ }{\it
(1) $W$ acts transitively and freely on $\tilt_1\Lambda$ by
\[T^{s_i}:=\nu_{Te_i}(T)=\left\{\begin{array}{ll}
\nu^-_{Te_i}(T)=\hom_{\Lambda^{\op}}(I_i,T)&\mbox{if $T\otimes_\Lambda S_i=0$}\\
\nu^+_{Te_i}(T)=TI_i&\mbox{if $\tor_1^\Lambda(T,S_i)=0$}
\end{array}\right.\]
for any $T\in\tilt_1\Lambda$ and $1\le i\le n$.

(2) Under the induced bijection $W\ni w\mapsto\Lambda^w\in\ii(\Lambda)$, 

\strut\kern1em(i) the Bruhat order on $W$ coincides with the reverse inclusion relation on $\ii(\Lambda)$,

\strut\kern1em(ii) the right order on $W$ coincides with the order on $\tilt_1\Lambda=\ii(\Lambda)$ (section \XE),

\strut\kern1em(iii) the left order on $W$ coincides with the order on $\tilt_1\Lambda^{\op}=\ii(\Lambda)$,

\strut\kern1em(iv) $\Lambda^w=I_{a_1}\cdots I_{a_k}$ holds for any reduced expression $w=s_{a_1}\cdots s_{a_k}$.}

\vskip.5em
In order to prove this result we consider the action of tilting complexes on the derived category. Let $K_0(\Lambda)$ be the Grothendieck group of $\Lambda$ and $K_0(\Lambda)_{\ccc}:=K_0(\Lambda)\otimes_{\zzz}\ccc$, which has the basis $\{[P_i]\ |\ 1\le i\le n\}$ because $\Lambda$ has finite global dimension. For an arbitrary two-sided tilting complex $T$ of $\Lambda$, we have an autoequivalence $T\Lotimes_\Lambda-$ of $\dd^{\b}(\mod\Lambda)$. Thus we get a map $\tilt_1\Lambda\to\GL(K_0(\Lambda)_{\ccc})$ defined by $T\mapsto[T\Lotimes_\Lambda-]$. We will compare with the contragredient of the {\it geometric representation} of $W$ [BB;4.1,4.2] defined as follows: Put $k_{i,i}:=-2$. For $i\neq j$, put $(k_{i,j},k_{j,i}):=(a,b)$ for $\stackrel{i}{\bullet}\stackrel{(a\ b)}{\mbox{------}}\stackrel{j}{\bullet}$, $(k_{i,j},k_{j,i}):=(a,2b)$ for $\begin{picture}(10,10)\put(7,4){\circle{15}}\end{picture}\stackrel{i}{\bullet}\stackrel{(a,b)}{\mbox{------}}\stackrel{j}{\bullet}$, and $(k_{i,j},k_{j,i}):=(2,2)$ for $\begin{picture}(10,10)\put(7,4){\circle{15}}\end{picture}\stackrel{i}{\bullet}\stackrel{}{\mbox{------}}\stackrel{j}{\bullet}\begin{picture}(10,10)\put(3,4){\circle{15}}\end{picture}$. Let $V^*$ be a vector space with basis $\alpha^*_1,\cdots,\alpha^*_n$. Define $\sigma^*_i\in\GL(V^*)$ by $\sigma^*_i(p):=p+p_i\sum_{j=1}^nk_{i,j}\alpha^*_j$ for $p=\sum_{j=1}^np_j\alpha^*_j$. It is well-known that the map $s_i\mapsto\sigma^*_i$ extends uniquely to an injective homomorphism $\sigma^*:W\to\GL(V^*)$, $w\mapsto\sigma^*_w$ [BB;4.2.7]. The following result, which is also interesting itself, shows that the autoequivalence induced by a tilting module has similar properties as $\sigma^*$. As we shall explain later, it is closely related to a result of Seidel-Thomas [ST].

\vskip.5em{\bf Theorem \XFF\ }{\it
(1) Let $V^*\to K_0(\Lambda)_{\ccc}$ be an isomorphism defined by $\alpha^*_i\mapsto[P_i]$. Then the induced isomorphism $\GL(V^*)\to\GL(K_0(\Lambda)_{\ccc})$ satisfies $\sigma^*_w\mapsto[\Lambda^w\Lotimes_\Lambda-]$ for any $w\in W$.
%For any $w\in W$, the automorphism $\Lambda^w\Lotimes_\Lambda-$ of $K_0(\Lambda)$ coincides with the geometric representation of $W$, when identifying $K_0(\Lambda)$ and $V$.

(2) We have an action $t_i\mapsto(I_i\Lotimes_\Lambda-)$ of the braid group $B$ on $\dd^{\b}(\mod\Lambda)$.}

\vskip.5em
We will give a proof of \XFE\ and \XFF\ after giving a series of preliminary results. Let us start with the following observation on tilting modules associated with a set of simple modules.

\vskip.5em{\bf Lemma \XFG\ }{\it
Let ${\bf S}$ be a set of simple $\Lambda^{\op}$-modules and $T\in\tilt_1\Lambda=\ii(\Lambda)$. Put
\[
e:=\sum_{1\le i\le n,\ S_i\in{\bf S}}e_i,\ \ \ I_{\bf S}:=\Lambda(1-e)\Lambda,\ \ \ U:=TI_{\bf S}\ \mbox{ and }\ V:=\hom_{\Lambda^{\op}}(I_{\bf S},T).
\]

\vskip-.5em
(1) $U$ is minimal amongst sub $\Lambda^{\op}$-modules of $T$ such that any composition factor of $T/U$ is in ${\bf S}$, and $V$ is maximal amongst sub $\Lambda^{\op}$-modules of $\Lambda$ such that $T\subset V$ and any composition factor of $V/T$ is in ${\bf S}$.

(2) $\top T_\Lambda$ and $\soc(\Lambda/T)_\Lambda$ do not have any common composition factor.

(3) $U,V\in\tilt_1\Lambda$.

(4) There is a path $V\to\cdots\to T\to\cdots\to U$ in $\tilt_1\Lambda$ with arrows indexed by ${\bf S}$.

(5) $V/U$ is a projective-injective $(\Lambda/I_{\bf S})^{\op}$-module which is a generator-cogenerator.}

%\[V/U\in\left\{\begin{array}{cc}\add(S_i\oplus S_j)_\Lambda&\mbox{if there is no edge between $i$ and $j$ in $\Delta$,}\\ \add(S_i\oplus S_j\oplus{S_i\choose S_j}\oplus{S_j\choose S_i})_\Lambda&\mbox{otherwise,}\end{array}\right.\] where we denote by ${S_i\choose S_j}$ (resp. ${S_j\choose S_i}$) the middle term of the non-zero element in $\ext^1_{\Lambda^{\op}}(S_i,S_j)$ (resp. $\ext^1_{\Lambda^{\op}}(S_j,S_i)$) which is uniquely determined up to isomorphism.}

\vskip.5em{\sc Proof }
(1) is obvious. We obtain (2) by \XFA(2). Since we can obtain $U$ (resp. $V$) by applying $\nu^+$ (resp. $\nu^-$) to $T$ repeatedly, we have (3) and (4). We will now show (5). Since $\Lambda/I_{\bf S}$ is selfinjective by \XEE(3), we only have to show that $V/U$ is progenerator. Let
$0\to P_1\to P_0\stackrel{f}{\to}V\to0$
be a minimal projective resolution of the $\Lambda^{\op}$-module $V$. By the choice of $V$, $\soc(\Lambda/V)_\Lambda$ does not contain any module in ${\bf S}$. We have a minimal projective resolution
$P_0^*\to P_1^*\to\ext^2_{\Lambda^{\op}}(\Lambda/V,\Lambda)\to0$.
Since $\ext^2_{\Lambda^{\op}}(\Lambda/V,\Lambda)=D(\Lambda/V)$ holds, $\top P_1$ does not contain any module in ${\bf S}$. Thus $P_1\subseteq P_0I_{\bf S}$ holds. Thus we have $P_0I_{\bf S}/P_1=U$ and $f^{-1}(U)=P_0I_{\bf S}$. Then $V/U\simeq P_0/P_0I_{\bf S}=P_0\otimes_\Lambda(\Lambda/I_{\bf S})$ implies that $V/U$ is a projective $(\Lambda/I_{\bf S})^{\op}$-module. Moreover, \XBI\ implies that any module in ${\bf S}$ appears in $\top P_0=\top V$. Thus $V/U$ is a generator.\rule{5pt}{10pt}

%Any composition factor of the $\Lambda^{\op}$-module $V/U$ is $S_i$ or $S_j$. Since $\Lambda$ is 2-spherical, we have $\ext^1_{\Lambda^{\op}}(S_i,S_i)=0=\ext^1_{\Lambda^{\op}}(S_j,S_j)$. If there is no edge between $i$ and $j$, then we have $\ext^1_{\Lambda^{\op}}(S_i\oplus S_j,S_i\oplus S_j)=0$. Thus we have $V/U\in\add(S_i\oplus S_j)_\Lambda$. Othersiwe, we have $\dim_k\ext^1_{\Lambda^{\op}}(S_i,S_j)=1=\dim_k\ext^1_{\Lambda^{\op}}(S_j,S_i)$. Hence ${S_i\choose S_j}$ and ${S_j\choose S_i}$ are uniquely determined. One can easily check that $\ext^1_{\Lambda^{\op}}(S_i\oplus S_j,{S_i\choose S_j}\oplus{S_j\choose S_i})=0=\ext^1_{\Lambda^{\op}}({S_i\choose S_j}\oplus{S_j\choose S_i},S_i\oplus S_j)$. This implies that there is no indecomposable $\Lambda^{\op}$-module $X$ with length at least 3 whose composition factor is $S_i$ or $S_j$.\rule{5pt}{10pt}

\vskip.5em
The next result is needed for showing that the action of $W$ described in \XFE\ is well-defined.

\vskip.5em{\bf Proposition \XFH\ }{\it $T^{s_i^2}=T$ holds for any $T\in\tilt_1\Lambda$ and $1\le i\le n$.}

\vskip.5em{\sc Proof }
Put ${\bf S}:=\{S_i\}$ and define $U$ and $V$ by \XFG. Since any composition factor of the $(\Lambda/I_{\bf S})^{\op}$-module $V/U$ is $S_i$ by \XFG(1), either $U=T$ or $V=T$ holds by \XFG(2). If $U=T$, then we have $T^{s_i}=V$ and $V^{s_i}=T$. If $V=T$, then we have $T^{s_i}=U$ and $U^{s_i}=T$.\rule{5pt}{10pt}

\vskip.5em
Now we prove theorem \XFE\ for the case when $\Lambda$ has only two simple modules. Then the quiver of $\Lambda$ is the double of one of the following valued graphs:

$\stackrel{}{\bullet}\stackrel{(2\ 2)}{\mbox{------}}\stackrel{}{\bullet},\ \ 
\begin{picture}(10,10)\put(7,4){\circle{15}}\end{picture}\stackrel{}{\bullet}\stackrel{}{\mbox{------}}\stackrel{}{\bullet}\begin{picture}(10,10)\put(3,4){\circle{15}}\end{picture}\ \mbox{ or }\ 
\begin{picture}(10,10)\put(7,4){\circle{15}}\end{picture}\stackrel{}{\bullet}\stackrel{(a\ b)}{\mbox{------}}\stackrel{}{\bullet}\ ((a\ b)=(1\ 2)\mbox{ or }(2\ 1)).$

The corresponding affine Weyl group $W$ is presented by generators $s_1,s_2$ and relations $s_1^2=s_2^2=1$. Thus we obtain \XFE(1) in these cases by \XFH. Since any element in $W$ can be written as $s_1s_2s_1s_2\cdots$ or $s_2s_1s_2s_1\cdots$, we obtain \XFE(2) in these cases by the following proposition.

\vskip.5em{\bf Proposition \XFI\ }{\it
Assume that $\Lambda$ has only two simple modules.

(1) The Hasse quiver of $\tilt_1\Lambda$ is as follows:
%{\small\[\begin{diag}\Lambda&\RA{2}&I_2&\RA{1}&I_2I_1&\RA{2}&I_2I_1I_2&\RA{1}&I_2I_1I_2I_1&\RA{2}&\cdots\\ &\DRA{1}&I_1&\RA{2}&I_1I_2&\RA{1}&I_1I_2I_1&\RA{2}&I_1I_2I_1I_2&\RA{1}&\cdots\end{diag}\]}

\ \ \ \ \ \ \ \ \ \ {\small$\begin{array}{ccccccccccc}
\Lambda&\stackrel{2}{\longrightarrow}&I_2&\stackrel{1}{\longrightarrow}&I_2I_1&\stackrel{2}{\longrightarrow}&I_2I_1I_2&\stackrel{1}{\longrightarrow}&I_2I_1I_2I_1&\stackrel{2}{\longrightarrow}&\cdots\\
&\begin{picture}(20,0)\put(2,-2){$\scriptstyle 1$}\put(2,10){\vector(1,-1){15}}\end{picture}&&&&&&&&&\\
&&I_1&\stackrel{2}{\longrightarrow}&I_1I_2&\stackrel{1}{\longrightarrow}&I_1I_2I_1&\stackrel{2}{\longrightarrow}&I_1I_2I_1I_2&\stackrel{1}{\longrightarrow}&\cdots
\end{array}$}

(2) The Hasse quiver of $\tilt_1\Lambda^{\op}$ is as follows:

\ \ \ \ \ \ \ \ \ \ {\small$\begin{array}{ccccccccccc}
\Lambda&\stackrel{2}{\longrightarrow}&I_2&&I_2I_1&&I_2I_1I_2&&I_2I_1I_2I_1&&\cdots\\
&\begin{picture}(20,0)\put(2,-1){$\scriptstyle 1$}\put(2,10){\vector(1,-1){15}}\end{picture}&&\begin{picture}(20,0)\put(2,10){$\scriptstyle 1$}\put(12,10){$\scriptstyle 2$}\put(2,-5){\vector(1,1){15}}\put(2,10){\vector(1,-1){15}}\end{picture}&&\begin{picture}(20,0)\put(2,10){$\scriptstyle 1$}\put(12,10){$\scriptstyle 2$}\put(2,-5){\vector(1,1){15}}\put(2,10){\vector(1,-1){15}}\end{picture}&&\begin{picture}(20,0)\put(2,10){$\scriptstyle 1$}\put(12,10){$\scriptstyle 2$}\put(2,-5){\vector(1,1){15}}\put(2,10){\vector(1,-1){15}}\end{picture}&&\begin{picture}(20,0)\put(2,10){$\scriptstyle 1$}\put(12,10){$\scriptstyle 2$}\put(2,-5){\vector(1,1){15}}\put(2,10){\vector(1,-1){15}}\end{picture}&\\
&&I_1&&I_1I_2&&I_1I_2I_1&&I_1I_2I_1I_2&&\cdots
\end{array}$}

%{\small\[\begin{diag}\Lambda&\RA{2}&I_2&\RA{1}&I_1I_2&\RA{2}&I_2I_1I_2&\RA{1}&I_1I_2I_1I_2&\RA{2}&\cdots\\ &\DRA{1}&I_1&\RA{2}&I_2I_1&\RA{1}&I_1I_2I_1&\RA{2}&I_2I_1I_2I_1&\RA{2}&\cdots\end{diag}\]}

(3) The Hasse quiver of the reverse inclusion order in $\ii(\Lambda)$ is as follows:

\ \ \ \ \ \ \ \ \ \ {\small$\begin{array}{ccccccccccc}
\Lambda&\stackrel{}{\longrightarrow}&I_2&\stackrel{}{\longrightarrow}&I_2I_1&\stackrel{}{\longrightarrow}&I_2I_1I_2&\stackrel{}{\longrightarrow}&I_2I_1I_2I_1&\stackrel{}{\longrightarrow}&\cdots\\
&\begin{picture}(20,0)\put(2,10){\vector(1,-1){15}}\end{picture}&&\begin{picture}(20,0)\put(2,-5){\vector(1,1){15}}\put(2,10){\vector(1,-1){15}}\end{picture}&&\begin{picture}(20,0)\put(2,-5){\vector(1,1){15}}\put(2,10){\vector(1,-1){15}}\end{picture}&&\begin{picture}(20,0)\put(2,-5){\vector(1,1){15}}\put(2,10){\vector(1,-1){15}}\end{picture}&&\begin{picture}(20,0)\put(2,-5){\vector(1,1){15}}\put(2,10){\vector(1,-1){15}}\end{picture}&\\
&&I_1&\stackrel{}{\longrightarrow}&I_1I_2&\stackrel{}{\longrightarrow}&I_1I_2I_1&\stackrel{}{\longrightarrow}&I_1I_2I_1I_2&\stackrel{}{\longrightarrow}&\cdots
\end{array}$}}

\vskip.5em{\sc Proof }
Using \XEC(2) and \XFA, we have (1). Considering $\Lambda^{\op}$, we have (2). To show (3), we only have to care about arrows, and it is easy.\rule{5pt}{10pt}

%Let $G:=\langle{-1\ 0\choose0\ -1}\rangle\subset\SL_2(K)$ and $\Lambda:=S*G$. Then $\Lambda$ has two primitive orthogonal idempotents $e_1$ and $e_2$.
% and the Loewy series of $\Lambda e_1$ and $\Lambda e_2$ are given by \[J_\Lambda^ne_i/J_\Lambda^{n+1}e_i=S_{i+n}^{n+1}\] for any $n\in\nnn_{\ge0}$, where the index $i+n$ of $S_{i+n}$ is regarded as an element of $\zzz/2\zzz$. 

\vskip,5em
In the rest of this section, assume that $\Lambda$ has at least three simple modules. We next investigate the Loewy series for the factor algebra given by the ideal associated with a set of two simple modules.

\vskip.5em{\bf Lemma \XFJ\ }{\it
(1) The (Loewy) length of $\Lambda/I_i$ is at most two.

(2) Let ${\bf S}:=\{S_i,S_j\}$ with $i\neq j$. Then the Loewy series of $\Lambda/I_{\bf S}$ is as follows:
{\small\[\begin{array}{|c|c|c|c|c|c|c|}\hline
\stackrel{i}{\bullet}\stackrel{}{\mbox{------}}\stackrel{j}{\bullet}&
\stackrel{i}{\bullet}\stackrel{(2\ 1)}{\mbox{------}}\stackrel{j}{\bullet}&
\stackrel{i}{\bullet}\stackrel{(1\ 2)}{\mbox{------}}\stackrel{j}{\bullet}&
\begin{picture}(10,10)\put(7,4){\circle{15}}\end{picture}\stackrel{i}{\bullet}\stackrel{}{\mbox{------}}\stackrel{j}{\bullet}&
\stackrel{i}{\bullet}\stackrel{}{\mbox{------}}\stackrel{j}{\bullet}\begin{picture}(10,10)\put(3,4){\circle{15}}\end{picture}&
\stackrel{i}{\bullet}\stackrel{(3\ 1)}{\mbox{------}}\stackrel{j}{\bullet}&
\stackrel{i}{\bullet}\stackrel{(1\ 3)}{\mbox{------}}\stackrel{j}{\bullet}\\
\hline\def\arraystretch{.5}\left[\begin{array}{c|c}
{\scriptstyle i}&{\scriptstyle j}\\
{\scriptstyle j}&{\scriptstyle i}\end{array}\right]&
\def\arraystretch{.5}\left[\begin{array}{c|c}
{\scriptstyle i}&{\scriptstyle j}\\ 
{\scriptstyle j\ j}&{\scriptstyle i}\\ 
{\scriptstyle i}&{\scriptstyle j}\end{array}\right]&
\def\arraystretch{.5}\left[\begin{array}{c|c}
{\scriptstyle i}&{\scriptstyle j}\\ 
{\scriptstyle j}&{\scriptstyle i\ i}\\ 
{\scriptstyle i}&{\scriptstyle j}\end{array}\right]&
\def\arraystretch{.5}\left[\begin{array}{c|c}
{\scriptstyle i}&{\scriptstyle j}\\ 
{\scriptstyle i\ j}&{\scriptstyle i}\\ 
{\scriptstyle i\ j}&{\scriptstyle i}\\ 
{\scriptstyle i}&{\scriptstyle j}\end{array}\right]&
\def\arraystretch{.5}\left[\begin{array}{c|c}
{\scriptstyle i}&{\scriptstyle j}\\ 
{\scriptstyle j}&{\scriptstyle i\ j}\\ 
{\scriptstyle j}&{\scriptstyle i\ j}\\ 
{\scriptstyle i}&{\scriptstyle j}\end{array}\right]&
\def\arraystretch{.5}\left[\begin{array}{c|c}
{\scriptstyle i}&{\scriptstyle j}\\ 
{\scriptstyle j\ j\ j}&{\scriptstyle i}\\ 
{\scriptstyle i\ i}&{\scriptstyle j\ j}\\ 
{\scriptstyle j\ j\ j}&{\scriptstyle i}\\ 
{\scriptstyle i}&{\scriptstyle j}\end{array}\right]&
\def\arraystretch{.5}\left[\begin{array}{c|c}
{\scriptstyle i}&{\scriptstyle j}\\ 
{\scriptstyle j}&{\scriptstyle i\ i\ i}\\ 
{\scriptstyle i\ i}&{\scriptstyle j\ j}\\ 
{\scriptstyle j}&{\scriptstyle i\ i\ i}\\ 
{\scriptstyle i}&{\scriptstyle j}\end{array}\right]\\
\hline\end{array}\]}}

\vskip.5em{\sc Proof }
By \XDC(1), the category $\pr\Lambda$ of projective $\Lambda$-modules over a 2-CY algebra $\Lambda$ forms a $\tau$-category in the sense of [I1,2]. For any idempotent $e$ of $\Lambda$, the category $\pr(\Lambda/I)$ also forms a $\tau$-category by [I2;1.4]. We give an indication of the proof, referring to [I1,2] for definitions. 

(1) We only have to consider the case when the quiver of $\Lambda$ has a loop at the vertex $i$. Put $\overline{\Lambda}:=\Lambda/I_i$. Since $\pr\overline{\Lambda}$ is a $\tau$-category, we have a minimal projective resolution $\overline{\Lambda}\stackrel{f}{\to}\overline{\Lambda}\stackrel{}{\to}\overline{\Lambda}\to S_i\to0$ such that the map $(f\cdot):\overline{\Lambda}\to\overline{\Lambda}$, which is obtained by applying $\hom_{\overline{\Lambda}}(-,\overline{\Lambda})$, has a simple cokernel. Thus $f$ does not belong to $J_{\overline{\Lambda}}^2$, and the cokernel of $f$ is simple. Thus $\overline{\Lambda}$ has (Loewy) length two.

%By [I1;4.1,7.1], we have a commutative diagram {\scriptsize\[\begin{diag}\Lambda/I_i&\supset&J_{\Lambda/I_i}&\supset&J_{\Lambda/I_i}^2\\ \uparrow&&\uparrow&&\uparrow\\ \Lambda/I_i&\LA{}&\Lambda/I_i&\LA{}&0\\ \uparrow&&\uparrow&&\uparrow\\ 0&\LA{}&\Lambda/I_i&\LA{}&\Lambda/I_i \end{diag}\]} such that each column gives the first two terms of a projective resolution of $J_{\Lambda/I_i}^k$. In particular, $\Lambda/I_i$ has length $2$.

(2) Put $\overline{\Lambda}:=\Lambda/I_{\bf S}$. We explain for the case $\stackrel{i}{\bullet}\stackrel{(3\ 1)}{\mbox{------}}\stackrel{j}{\bullet}$. In this case, the first three terms of the minimal projective resolutions of the simple $\overline{\Lambda}^{\op}$-modules $S_i$ and $S_j$ are given by $\overline{P}_i\to \overline{P}_j^3\to \overline{P}_i\to S_i\to0$ and $\overline{P}_j\to \overline{P}_i\to \overline{P}_j\to S_j\to0$. By [I1;4.1,7.1], we have commutative diagrams
{\scriptsize\[\begin{diag}
\overline{P}_i&\supset&\overline{P}_iJ_{\overline{\Lambda}}&\supset&\overline{P}_iJ_{\overline{\Lambda}}^2&\supset&\overline{P}_iJ_{\overline{\Lambda}}^3&\supset&\overline{P}_iJ_{\overline{\Lambda}}^4&\supset&\overline{P}_iJ_{\overline{\Lambda}}^5\\
\uparrow&&\uparrow&&\uparrow&&\uparrow&&\uparrow&&\uparrow\\
\overline{P}_i&\longleftarrow&\overline{P}_j^3&\longleftarrow&\overline{P}_i^2&\longleftarrow&\overline{P}_j^3&\longleftarrow&\overline{P}_i&\longleftarrow&0\\
\uparrow&*&\uparrow&*&\uparrow&*&\uparrow&*&\uparrow&*&\uparrow\\
0&\longleftarrow&\overline{P}_i&\longleftarrow&\overline{P}_j^3&\longleftarrow&\overline{P}_i^2&\longleftarrow&\overline{P}_j^3&\longleftarrow&\overline{P}_i
\end{diag}\ \ \ \ \ \ \ \ \ \ \begin{diag}
\overline{P}_j&\supset&\overline{P}_jJ_{\overline{\Lambda}}&\supset&\overline{P}_jJ_{\overline{\Lambda}}^2&\supset&\overline{P}_jJ_{\overline{\Lambda}}^3&\supset&\overline{P}_jJ_{\overline{\Lambda}}^4&\supset&\overline{P}_jJ_{\overline{\Lambda}}^5\\
\uparrow&&\uparrow&&\uparrow&&\uparrow&&\uparrow&&\uparrow\\
\overline{P}_j&\longleftarrow&\overline{P}_i&\longleftarrow&\overline{P}_j^2&\longleftarrow&\overline{P}_i&\longleftarrow&\overline{P}_j&\longleftarrow&0\\
\uparrow&*&\uparrow&*&\uparrow&*&\uparrow&*&\uparrow&*&\uparrow\\
0&\longleftarrow&\overline{P}_j&\longleftarrow&\overline{P}_i&\longleftarrow&\overline{P}_j^2&\longleftarrow&\overline{P}_i&\longleftarrow&\overline{P}_j
\end{diag}\]}
called ladders, with the following properties:

(i) Each column gives the first two terms of a projective resolution of $\overline{P}_iJ_{\overline{\Lambda}}^k$ and $\overline{P}_jJ_{\overline{\Lambda}}^k$. 

(ii) The mapping cone of each commutative square $*$ gives the first three terms of a minimal projective resolution of semisimple $\overline{\Lambda}^{\op}$-modules.

By (i), the Loewy series of $\overline{P}_i$ is $(\overline{P}_iJ_{\overline{\Lambda}}^k/\overline{P}_iJ_{\overline{\Lambda}}^{k+1})_{k\ge0}=(S_i,S_j^3,S_i^2,S_j^3,S_i)$, and that of $\overline{P}_j$ is $(\overline{P}_jJ_{\overline{\Lambda}}^k/\overline{P}_jJ_{\overline{\Lambda}}^{k+1})_{k\ge0}=(S_j,S_i,S_j^2,S_i,S_j)$.\rule{5pt}{10pt}

%Since $\overline{\Lambda}$ is selfinjective by \XEF, the co-Loewy series of the $\overline{\Lambda}^{\op}$-module $\overline{\Lambda}$ is the dual of the Loewy series of the $\overline{\Lambda}$-module $\overline{\Lambda}$. Thus the similar calculation as above works, and we obtain the assertion.\rule{5pt}{10pt}

\vskip.5em{\bf Lemma \XFK\ }{\it
For ${\bf S}:=\{S_i,S_j\}$ with $i\neq j$, put $\overline{\Lambda}:=\Lambda/I_{\bf S}$. Let $\overline{P}\in\mod\overline{\Lambda}^{\op}$ be a progenerator. Define $X_k,Y_k\in\mod\overline{\Lambda}^{\op}$ by $X_0=Y_0:=\overline{P}$, $X_{2k+1}:=X_{2k}I_i$, $X_{2k+2}:=X_{2k+1}I_j$, $Y_{2k+1}:=Y_{2k}I_j$ and $Y_{2k+2}:=Y_{2k+1}I_i$ for $k\ge0$. Then we have $X_{m(i,j)}=0=Y_{m(i,j)}$ and $X_k\neq0$, $Y_k\neq0$ for any $k<m(i,j)$.}

%Put $k_n:=i$ if $n$ is even, and $k_n:=j$ otherwise. If $\top M_n$ contains $S_{k_n}$, then $M_{n+1}:=M_nI_{k_n}$. Otherwise, let $M_{n+1}$ be the maximum element among left ideals of $\overline{\Lambda}$ such that $M_n\subset M_{n+1}$ and any composition factor of $M_{n+1}/M_n$ is $S_{k_n}$.
%Then the sequence has period $2m(i,j)$ and forms the diagram {\small\[\begin{diag}M_0=\overline{P}&\supsetneq&M_1&\supsetneq&\cdots&\supsetneq&M_{m(i,j)-1}&\supsetneq&M_{m(i,j)}\\ \parallel&&&&&&&&\parallel\\ M_{2m(i,j)}&\supsetneq&M_{2m(i,j)-1}&\supsetneq&\cdots&\supsetneq&M_{m(i,j)+1}&\supsetneq&0.\end{diag}\]}}

\vskip.5em{\sc Proof }
Obviously we only have to consider the case $\overline{P}=\overline{\Lambda}$. 

(i) We consider the case when the quiver of $\overline{\Lambda}$ has no loop. In this case, we can easily check the assertion by using \XFJ(2). For example, the calculation of $X_k$ and $Y_k$ for the case $\stackrel{i}{\bullet}\stackrel{(2\ 1)}{\mbox{------}}\stackrel{j}{\bullet}$ is as follows: 
{\small\begin{eqnarray*}
X_0=\overline{\Lambda}=\def\arraystretch{.5}\left[\begin{array}{c|c}
{\scriptstyle i}&{\scriptstyle j}\\ 
{\scriptstyle j\ j}&{\scriptstyle i}\\ 
{\scriptstyle i}&{\scriptstyle j}\end{array}\right]\supset
X_1=\def\arraystretch{.5}\left[\begin{array}{c|c}
{\scriptstyle }&{\scriptstyle j}\\ 
{\scriptstyle j\ j}&{\scriptstyle i}\\ 
{\scriptstyle i}&{\scriptstyle j}\end{array}\right]\supset
X_2=\def\arraystretch{.5}\left[\begin{array}{c|c}
{\scriptstyle }&{\scriptstyle }\\ 
{\scriptstyle }&{\scriptstyle i}\\ 
{\scriptstyle i}&{\scriptstyle j}\end{array}\right]\supset
X_3=\def\arraystretch{.5}\left[\begin{array}{c|c}
{\scriptstyle }&{\scriptstyle }\\ 
{\scriptstyle }&{\scriptstyle }\\ 
{\scriptstyle }&{\scriptstyle j}\end{array}\right]\supset
X_4=0\\
Y_0=\overline{\Lambda}=\def\arraystretch{.5}\left[\begin{array}{c|c}
{\scriptstyle i}&{\scriptstyle j}\\ 
{\scriptstyle j\ j}&{\scriptstyle i}\\ 
{\scriptstyle i}&{\scriptstyle j}\end{array}\right]\supset
Y_1=\def\arraystretch{.5}\left[\begin{array}{c|c}
{\scriptstyle i}&{\scriptstyle }\\ 
{\scriptstyle j\ j}&{\scriptstyle i}\\ 
{\scriptstyle i}&{\scriptstyle j}\end{array}\right]\supset
Y_2=\def\arraystretch{.5}\left[\begin{array}{c|c}
{\scriptstyle }&{\scriptstyle }\\ 
{\scriptstyle j\ j}&{\scriptstyle }\\ 
{\scriptstyle i}&{\scriptstyle j}\end{array}\right]\supset
Y_3=\def\arraystretch{.5}\left[\begin{array}{c|c}
{\scriptstyle }&{\scriptstyle }\\ 
{\scriptstyle }&{\scriptstyle }\\ 
{\scriptstyle i}&{\scriptstyle }\end{array}\right]\supset
Y_4=0\end{eqnarray*}}
Thus the assertion follows. One can check the other cases similarly.

%$M_n=\overline{P}_iJ_{\overline{\Lambda}}^n\oplus\overline{P}_jJ_{\overline{\Lambda}}^{n+1}$ holds for any $n\le m$.Moreover,  and $M_n=\overline{P}_iJ_{\overline{\Lambda}}^n\oplus\overline{P}_jJ_{\overline{\Lambda}}^{n+1}$ holds for any $n\le m$, and 

(ii) We consider the case $\begin{picture}(10,10)\put(7,4){\circle{15}}\end{picture}\stackrel{i}{\bullet}\stackrel{}{\mbox{------}}\stackrel{j}{\bullet}$. It follows from \XFJ(1) that $\Lambda/I_i$ has length $2$. Thus the composition factors of $\overline{P}_i/\overline{P}_iI_i$ are two copies of $S_i$. Since $\overline{\Lambda}$ is selfinjective, $\overline{P}_i$ contains a submodule $X$ whose composition factors are two copies of $S_i$. Then $X$ is contained in $\overline{P}_iI_i$, and the composition factors of $\overline{P}_iI_i/X$ are two copies of $S_j$. Thus we can calculate $M_i$ as follows:
{\small\begin{eqnarray*}
X_0=\overline{\Lambda}=\def\arraystretch{.5}\left[\begin{array}{c|c}
{\scriptstyle i}&{\scriptstyle j}\\ 
{\scriptstyle i\ j}&{\scriptstyle i}\\ 
{\scriptstyle i\ j}&{\scriptstyle i}\\ 
{\scriptstyle i}&{\scriptstyle j}\end{array}\right]\supset
X_1=\def\arraystretch{.5}\left[\begin{array}{c|c}
{\scriptstyle }&{\scriptstyle j}\\ 
{\scriptstyle j}&{\scriptstyle i}\\ 
{\scriptstyle i\ j}&{\scriptstyle i}\\ 
{\scriptstyle i}&{\scriptstyle j}\end{array}\right]\supset
X_2=\def\arraystretch{.5}\left[\begin{array}{c|c}
{\scriptstyle }&{\scriptstyle }\\ 
{\scriptstyle }&{\scriptstyle i}\\ 
{\scriptstyle i}&{\scriptstyle i}\\ 
{\scriptstyle i}&{\scriptstyle j}\end{array}\right]\supset
X_3=\def\arraystretch{.5}\left[\begin{array}{c|c}
{\scriptstyle }&{\scriptstyle }\\ 
{\scriptstyle }&{\scriptstyle }\\ 
{\scriptstyle }&{\scriptstyle }\\ 
{\scriptstyle }&{\scriptstyle j}\end{array}\right]\supset
X_4=0\\
Y_0=\overline{\Lambda}=\def\arraystretch{.5}\left[\begin{array}{c|c}
{\scriptstyle i}&{\scriptstyle j}\\ 
{\scriptstyle i\ j}&{\scriptstyle i}\\ 
{\scriptstyle i\ j}&{\scriptstyle i}\\ 
{\scriptstyle i}&{\scriptstyle j}\end{array}\right]\supset
Y_1=\def\arraystretch{.5}\left[\begin{array}{c|c}
{\scriptstyle i}&{\scriptstyle }\\ 
{\scriptstyle i\ j}&{\scriptstyle i}\\ 
{\scriptstyle i\ j}&{\scriptstyle i}\\ 
{\scriptstyle i}&{\scriptstyle j}\end{array}\right]\supset
Y_2=\def\arraystretch{.5}\left[\begin{array}{c|c}
{\scriptstyle }&{\scriptstyle }\\ 
{\scriptstyle j}&{\scriptstyle }\\ 
{\scriptstyle i\ j}&{\scriptstyle }\\ 
{\scriptstyle i}&{\scriptstyle j}\end{array}\right]\supset
Y_3=\def\arraystretch{.5}\left[\begin{array}{c|c}
{\scriptstyle }&{\scriptstyle }\\ 
{\scriptstyle }&{\scriptstyle }\\ 
{\scriptstyle i}&{\scriptstyle }\\ 
{\scriptstyle i}&{\scriptstyle }\end{array}\right]\supset
Y_4=0\end{eqnarray*}}
Thus the assertion follows. We can treat the case $\stackrel{i}{\bullet}\stackrel{}{\mbox{------}}\stackrel{j}{\bullet}\begin{picture}(10,10)\put(3,4){\circle{15}}\end{picture}$ similarly.\rule{5pt}{10pt}

%\vskip.5em{\bf\XFK\ Lemma }{\it Assume $i\neq j$ and there is no edge between $i$ and $j$ in $\Delta$. (1) $T^{(s_is_j)^2}=T$ holds for any $T\in\tilt_1\Lambda$. (2) Put $T_0:=T$, $T_1:=T^{s_i}$, $T_2:=T^{s_is_j}$ and $T_3:=T^{s_is_js_i}$. Then, for some $k\in\zzz/4\zzz$, $\tilt_1\Lambda$ has a subquiver {\small\[\begin{diag}T_k&\RA{i}&T_{k+1}&\RA{j}&T_{k+2}.\\ &\DRA{j}&T_{k+3}&\URA{i} \end{diag}\]}} \vskip-1em{\sc Proof } Put ${\bf S}:=\{S_i,S_j\}$ and define $U$ and $V$ by \XFG. Then neither $S_i$ nor $S_j$ is a composition factor of $\top U_\Lambda$. Since $V\in\tilt_1\Lambda$, we have $S_i\oplus S_j\in\add\top V_\Lambda$ by \XFG(1). Thus we have $S_i\oplus S_j\in\add\top(V/U)_\Lambda$. By \XFG(4), we can put $V/U=S_i^a\oplus S_j^b$ with $ab>0$. Define $W_l$ ($l=1,2$) by $U\subset W_l\subset V$ and \[W_1/U=S_j^b\subset V/U,\ \ \ \ \ W_2/U=S_i^a\subset V/U.\] Since $ab>0$ holds, $\tilt_1\Lambda$ has a subquiver {\small\[\begin{diag} V&\RA{i}&W_1&\RA{j}&U.\\ &\DRA{j}&W_2&\URA{i} \end{diag}\]} \vskip-1em By \XFG(3), $T$ is either $U,W_l$ ($l=1,2$) or $V$. Thus the assertion follows.\rule{5pt}{10pt}

\vskip.5em
The following result is crucial for well-defined action of the affine Weyl and braid groups.

\vskip.5em{\bf Proposition \XFL\ }{\it
Let $i\neq j$ be distinct vertices in $\Delta$,

(1) $T^{(s_is_j)^{m(i,j)}}=T$ holds for any $T\in\tilt_1\Lambda$.

(2) $\tilt_1\Lambda$ has a subquiver
{\small\[\begin{diag}
V&\RA{i}&\bullet&\RA{j}&\bullet&\RA{i}&\cdots\cdots&\RA{i\mbox{\tiny\ or }j}&\bullet&\RA{j\mbox{\tiny\ or }i}&U.\\
&\DRA{j}&\bullet&\RA{i}&\bullet&\RA{j}&\cdots\cdots&\RA{j\mbox{\tiny\ or }i}&\bullet&\URA{i\mbox{\tiny\ or }j}
\end{diag}\]}
with two paths of length $m(i,j)$ such that $T$ is one of these vertices.

(3) We have an equality $I_i\Lotimes_\Lambda I_j\Lotimes_\Lambda I_i\Lotimes_\Lambda I_j\Lotimes_\Lambda\cdots=I_iI_jI_iI_j\cdots=\Lambda(1-e_i-e_j)\Lambda=I_jI_iI_jI_i\cdots=I_j\Lotimes_\Lambda I_i\Lotimes_\Lambda I_j\Lotimes_\Lambda I_i\Lotimes_\Lambda\cdots$, where each derived tensor product and product of ideals contains exactly $m(i,j)$ terms.}

\vskip.5em{\sc Proof }
Put ${\bf S}:=\{S_i,S_j\}$, and consider $U$ and $V$ as defined in \XFG. Then $P:=V/U$ is a progenerator of $\Lambda/I_{\bf S}$ by \XFG(5). Consider the sequences $(X_0,X_1,\cdots,X_{m(i,j)})$ and $(Y_0,Y_1,\cdots,Y_{m(i,j)})$ in \XFK. Then the preimage of these sequences under the natural surjection $V\to P=V/U$ coincides with the sequences $(V,V^{s_i},V^{s_is_j},V^{s_is_js_i},\cdots)$ and $(V^{s_j},V^{s_js_i},V^{s_js_is_j},\cdots)$. Since $X_{m(i,j)}=0=Y_{m(i,j)}$ holds, we have $V^{s_is_js_i\cdots}=U=V^{s_js_is_j\cdots}$ where both $s_is_js_i\cdots$ and $s_js_is_j\cdots$ are products of $m(i,j)$ simple reflections. Thus we have proved (2). Since $T$ belongs to one of these sequences by \XFG(4), we obtain (1). Applying (2) and \XFA(2) to $T:=\Lambda$, we have (3).\rule{5pt}{10pt}

\vskip.5em
We can now prove most of \XFE\ and \XFF.

\vskip.5em{\bf Proof of \XFE(1) and \XFF\ }
We first show \XFE(1). By \XFH\ and \XFL, $W$ acts on $\tilt_1\Lambda$. By \XFB, the action is transitive. For freeness of the action, it is enough to prove \XFF(1) since the geometric representation $\sigma^*:W\to\GL(V^*)$ is injective.

\XFF(2) follows from \XFL(3). We will show \XFF(1). First we show that $[I_i\Lotimes_\Lambda-]$ corresponds to $\sigma^*_i$. If $i\neq j$, then $[I_i\Lotimes_\Lambda P_j]=[I_ie_j]=[P_j]$. Let us calculate $[I_i\Lotimes_\Lambda P_i]=[I_ie_i]$. If the quiver of $\Lambda$ has no loop at the vertex $i$, then we have a minimal projective resolution $0\to P_i\to\bigoplus_{j\neq i}P_j^{k_{i,j}}\stackrel{f}{\to}P_i\to S_i\to 0$ with $\Im f_i=I_ie_i$. Thus we have $[I_ie_i]=[P_i]+\sum_{j=1}^nk_{i,j}[P_j]$. Assume that the quiver of $\Lambda$ has a loop at the vertex $i$. In this case, the composition factors of $P_i/I_ie_i$ are two copies of $S_i$ by \XFJ(1). Since we have a minimal projective resolution $0\to P_i\to P_i\oplus(\bigoplus_{j\neq i}P_j^{k_{i,j}/2})\stackrel{}{\to}P_i\to S_i\to 0$ by our definition of $k_{i,j}$, we have $[I_ie_i]=[P_i]-2[S_i]=[P_i]-2([P_i]-\sum_{j\neq i}k_{i,j}[P_j]/2)=[P_i]+\sum_{j=1}^nk_{i,j}[P_j]$. Thus $[I_i\Lotimes_\Lambda-]$ corresponds to $\sigma^*_i$.

Now take any $w=s_{a_1}\cdots s_{a_k}\in W$. Then $\Lambda^w=I_{a_1}^{\pm1}\Lotimes_\Lambda\cdots\Lotimes_\Lambda I_{a_k}^{\pm1}$ for $I_{a_i}^{-1}:=\Rhom_{\Lambda^{\op}}(I_{a_i},\Lambda)$. Since $[I_{a_i}^{-1}\Lotimes_\Lambda-]=[I_{a_i}\Lotimes_\Lambda-]$ by $(\sigma^*_i)^2=1_{V^*}$, we have that $[\Lambda^w\Lotimes_\Lambda-]=[I_{a_1}\Lotimes_\Lambda-]\cdots[I_{a_k}\Lotimes_\Lambda-]$ corresponds to $\sigma^*_{a_1}\cdots\sigma^*_{a_k}=\sigma^*_w$. Thus \XFF(1) follows.\rule{5pt}{10pt}

\vskip.5em
It remains to prove \XFE(2). For this, we need the following lemma. For simplicity, we write $T\stackrel{i}{\leftrightarrow}U$ if there is an arrow $T\stackrel{i}{\to}U$ or $T\stackrel{i}{\leftarrow}U$ in $\tilt_1\Lambda$. For $w\in W$, we denote by $k(w)$ the length of the shortest path in $\tilt_1\Lambda$ from $\Lambda$ to $\Lambda^w$.

\vskip.5em{\bf Lemma \XFM\ }{\it
Let $w,w^\prime\in W$.

(1) There is an arrow $\Lambda^{w^\prime}\stackrel{i}{\leftrightarrow}\Lambda^w$ in $\tilt_1\Lambda$ if and only if $w=w^\prime s_i$.

(2) There is a subquiver $\Lambda=T_0\stackrel{a_1}{\leftrightarrow}T_1\stackrel{a_2}{\leftrightarrow}\cdots\stackrel{a_k}{\leftrightarrow}T_k=\Lambda^w$ in $\tilt_1\Lambda$ if and only if $w=s_{a_1}\cdots s_{a_k}$ holds.

(3) Any path in $\tilt_1\Lambda$ from $\Lambda$ to $\Lambda^w$ has length $k(w)$.

(4) $\Lambda^w\le\Lambda^{w^\prime}$ in $\tilt_1\Lambda$ if and only if there exists a path in $\tilt_1\Lambda$ from $\Lambda^w$ to $\Lambda^{w^\prime}$.

(5) If there is an arrow $\Lambda^w\stackrel{}{\to}\Lambda^{w^\prime}$ in $\tilt_1\Lambda$, then $k(w^\prime)=k(w)+1$.

(6) $k(w)=l(w)$.

(7) There is a subquiver $\Lambda=T_0\stackrel{a_1}{\to}T_1\stackrel{a_2}{\to}\cdots\stackrel{a_k}{\to}T_k=\Lambda^w$ in $\tilt_1\Lambda$ if and only if $w=s_{a_1}\cdots s_{a_k}$ is a reduced expression of $w$.}

\vskip.5em{\sc Proof }
(1)(2) Immediate from \XFE(1).

(3) We use induction on $k(w)$. If $k(w)=0$, then $w=1$, and there is no non-trivial path in $\tilt_1\Lambda$ from $\Lambda$ to $\Lambda$. Assume that the assertion is true for any $w\in W$ with $k(w)<k$. Fix $w\in W$ with $k(w)=k$. Take a path $\Lambda\stackrel{a_1}{\to}\cdots\stackrel{a_k}{\to}\Lambda^w$ of length $k$, and an arbitrary path $\Lambda\stackrel{b_1}{\to}\cdots\stackrel{b_l}{\to}\Lambda^w$ of length $l$. We will show $l=k$. By \XFL, $\tilt_1\Lambda$ has a subquiver
{\small\[\begin{diag}
\Lambda^{v}&\RA{a_k\mbox{\tiny\ or }b_l}&\cdots\cdots&\RA{a_k}&\Lambda^{wa_kb_l}&\RA{b_l}&\Lambda^{wa_k}&\RA{a_k}&\Lambda^{w}\\
&\DRA{b_l\mbox{\tiny\ or }a_k}&\cdots\cdots&\RA{b_l}&\Lambda^{wb_la_k}&\RA{a_k}&\Lambda^{wb_l}&\URA{b_l}
\end{diag}\]}
consisting of two paths of length $m:=m(a_k,b_l)$. Since there is a path $\Lambda\stackrel{a_1}{\to}\cdots\stackrel{a_{k-1}}{\to}\Lambda^{wa_k}$ of length $k-1$, any path from $\Lambda$ to $\Lambda^{wa_k}$ has length $k-1$ by the inductive assumption. In particular, we have $k(v)=k-m$. Since we have $k(wb_l)\le k(v)+m-1=k-1$ and $k(wb_l)\ge k(w)-1=k-1$, we have $k(wb_l)=k-1$. Since we have a path $\Lambda\stackrel{b_1}{\to}\cdots\stackrel{b_{l-1}}{\to}\Lambda^{wb_l}$ of length $l-1$, we have $l-1=k-1$ by the inductive assumption again. Thus $l=k$.

%{\small\[\begin{diag}\Lambda^{wa_kb_m}&\RA{b_m}&\Lambda^{wa_k}&\RA{a_k}&\Lambda^w\\ &\DRA{a_k}&\Lambda^{wb_m}&\URA{b_m} \end{diag}\ \ \mbox{ or }\ \ \begin{diag}\Lambda^{wa_kb_ma_k}&\RA{a_k}&\Lambda^{wa_kb_m}&\RA{b_m}&\Lambda^{wa_k}&\RA{a_k}&\Lambda^w\\ &\DRA{b_m}&\Lambda^{wb_ma_k}&\RA{a_k}&\Lambda^{wb_m}&\URA{b_m} \end{diag}.\]}

(4) The `if' part is obvious. We show the `only if' part.
By \XEC(1)(ii), we have a path $\cdots\to T_1\to T_0=\Lambda^{w^\prime}$ such that $\Lambda^w\le T_i$ for any $i$. We have $\Lambda^w=T_k$ for some $k$ by (3).

(5) Immediate from (3).

(6) By (2), any path from $\Lambda$ to $\Lambda^w$ gives a presentation of $w$. Thus we have $k(w)\ge l(w)$. Take a presentation $w=s_{a_1}\cdots s_{a_l}$ with $l=l(w)$. By (2) again, there exists a subquiver $\Lambda=T_0\stackrel{a_1}{\leftrightarrow}\cdots\stackrel{a_l}{\leftrightarrow}T_l=\Lambda^w$ in $\tilt_1\Lambda$. This implies $k(w)\le l$ by (5).

(7) Immediate from (2)(5) and (6).\rule{5pt}{10pt}

\vskip.5em
We are now in the position to finish the proof of \XFE.

\vskip.5em{\bf Proof of \XFE(2) }
We first show (ii). By \XFM(1)(5) and (6), there is an arrow $\Lambda^{w^\prime}\stackrel{i}{\to}\Lambda^w$ in $\tilt_1\Lambda$ if and only if $w=w^\prime s_i$ and $l(w^\prime)<l(w)$. This implies that the Hasse quiver of $\tilt_1\Lambda$ and that of the poset $W$ with the right order coincide. It follows from \XFM(4) that (ii) holds. One can show (iii) by a dual argument.

We now show (i) and (iv). Assume $w^\prime\le w$. Take a reduced expression $w=s_{a_1}\cdots s_{a_k}$. Then $\Lambda^w=I_{a_1}\cdots I_{a_k}$ holds by \XFM(7).
Since $w^\prime=s_{a_{i_1}}\cdots s_{a_{i_q}}$ holds for some $1\le i_1<\cdots<i_q\le k$, we have $\Lambda^{w^\prime}\supseteq I_{a_{i_1}}\cdots I_{a_{i_q}}$. Thus $\Lambda^{w^\prime}\supseteq\Lambda^w$ holds.

Conversely, assume $\Lambda^{w^\prime}\supseteq\Lambda^{w}$. Using induction on $\length(\Lambda/\Lambda^w)_\Lambda$, we will show $w^\prime\le w$. Fix $S_i\in\soc(\Lambda/\Lambda^{w})_\Lambda$, so we have $ws_i<w$. If $\Lambda^{w^\prime}\supseteq\Lambda^{ws_i}$, then we have $w^\prime\le ws_i<w$ inductively. If $\Lambda^{w^\prime}\ {\not\supseteq}\ \Lambda^{ws_i}$, then $S_i\in\soc(\Lambda/\Lambda^{w^\prime})_\Lambda$ holds, and we have
{\small\[\begin{diag}
\Lambda^{w^\prime s_i}&\ \ \ \supseteq\ \ \ &\Lambda^{ws_i}\\
\cup&&\cup\\
\Lambda^{w^\prime}&\ \ \ \supseteq\ \ \ &\Lambda^w
\end{diag}\]}
\vskip-1em
Inductively, we have $w^\prime s_i\le ws_i$. Applying [Hu;5.9] to $w^\prime s_i$ and $ws_i$, either $w^\prime\le ws_i$ or $w^\prime\le w$ holds. In any case, we have $w^\prime\le w$. Thus we have shown (i).\rule{5pt}{10pt}

\vskip.5em
The {\it derived Picard group} $\dpic_R(\Lambda)$ of $\Lambda$ was introduced by Yekutieli [Ye2] (see also [RZ][MY]). The elements of $\dpic_R(\Lambda)$ are isoclasses of two-sided tilting complexes $T\in\dd^{\b}(\mod\Lambda\otimes_R\Lambda^{\op})$, and the multiplication of $T$ and $T'$ is given by $T\Lotimes_\Lambda T'$. Then the inverse of $T$ is given by $\Rhom_\Lambda(T,\Lambda)\simeq\Rhom_{\Lambda^{\op}}(T,\Lambda)$. We have a group homomorphism from $\dpic_R(\Lambda)$ to the group $\auteq_R(\dd^{\b}(\mod\Lambda))$ of autoequivalences of $\dd^{\b}(\mod\Lambda)$ defined by $T\mapsto(T\Lotimes_\Lambda-)$. When $\Lambda$ is 2-CY, we have the elements $I_1,\cdots,I_n$ of $\dpic_R(\Lambda)$ which satisfy the braid relations by \XFL(3). Inspired by study in algebraic geometry (e.g. [BO][Bri4][IU]), we have the following natural questions.

\vskip.5em{\bf Questions }
(1) Do $I_1,\cdots,I_n$ together with the shift $[1]$ and the outer automorphism group $\out_R(\Lambda)$ of $\Lambda$ generate $\dpic_R(\Lambda)$?

(2) Is $\dpic_R(\Lambda)$ isomorphic to $(B\timesl\out_R(\Lambda))\times\zzz$ for the affine braid group $B$ and the group $\zzz$ generated by $[1]$?

(3) Is the homomorphism $\dpic_R(\Lambda)\to\auteq_R(\dd^{\b}(\mod\Lambda))$ an isomorphism?

\vskip.5em
Now we consider the case $\Lambda=S*G$ for $S=K[[x,y]]$ for a field $K$ of characteristic 0 and a finite subgroup $G$ of $\SL_2(K)$. In this case, there is a triangle equivalence ({\it McKay correspondence}) between $\dd^{\b}(\mod\Lambda)$ and $\dd^{\b}(\Coh X)$ for a minimal resolution $X$ of the singularity $\Spec S^G$ [KV]. When $G$ is cyclic, Ishii-Uehara [IU] determined generators of the subgroup $\auteq^{\rm FM}(\dd^{\b}(\Coh X))$ consisting of Fourier-Mukai transformations. Consequently, the first question should have a positive answer for this case.

%If this is true, then all algebras which are derived equivalent to $\Lambda$ should be Morita equivalent to $\Lambda$. At least all classical tilting modules are generated by these ideals.
%(2) Does the braid group $B$ (together with the shift $[1]$ and $\aut(\Lambda)$) generate $\aut(\dd^{\b}(\mod\Lambda))$?

\vskip.5em
We end this section by showing that the spherical objects introduced by Seidel-Thomas [ST] give rise to tilting complexes, and hence to autoequivalences of $\dd^{\b}(\mod\Lambda)$. We assume that $\Lambda$ is a projective $R$-module and $\gl\Lambda<\infty$. We call $S\in\dd^{\b}(\mod\Lambda)$ {\it $n$-spherical} ($n>0$) if the following conditions (1) and (2) are satisfied.

(1) $E:=\endm_{\dd(\Mod\Lambda)}(S)$ is a division algebra,

(2) $\dim_E\hom_{\dd(\Mod\Lambda)}(S,S[i])=\left\{\begin{array}{cc}
1&\mbox{if $i=0$ or $n$,}\\
0&\mbox{otherwise.}
\end{array}\right.$

It follows from a result of Keller [K1,2] that condition (2) implies that $S$ comes from an object in $\dd^{\b}(\mod\Lambda\otimes_RE^{\op})$.
%Then, for any $X\in\dd^{\b}(\mod\Lambda)$, both $\bigoplus_{i\in\zzz}\hom_\Lambda(S,X[i])$ and $\bigoplus_{i\in\zzz}\hom_\Lambda(X,S[i])$ are finite dimensional over $E$.
Let us recall the definition due to Seidel and Thomas [ST] of autoequivalences of $\dd^{\b}(\mod\Lambda)$ called {\it twist functors}. We treat here only dual twist functors. To obtain functoriality, we first construct a functor $\ttt_S:\kk^{\b}(\pr\Lambda)\to\kk^{\b}(\mod\Lambda)$. For $P\in\kk^{\b}(\pr\Lambda)$, put
\[\ttt_S(P):=(P\stackrel{e_P}{\longrightarrow}\hom^\bullet_{E^{\op}}(\hom^\bullet_{\Lambda}(P,S),S)).\]
%\ttt_S(P):=(S\otimes^\bullet_E\hom^\bullet_{\Lambda}(S,P)\stackrel{a_P}{\longrightarrow}P),\ \ \ \ \ 
where $e_P$ is the evaluation map, and $\ttt_S(P)$ is defined as mapping cones such that $P$ is in degree 0.
Composing $\ttt_S$ with natural functors $\dd^{\b}(\mod\Lambda)\stackrel{\sim}{\to}\kk^{\b}(\pr\Lambda)$ and $\kk^{\b}(\mod\Lambda)\to\dd^{\b}(\mod\Lambda)$, we obtain an autoequivalence $\ttt_S$ on $\dd^{\b}(\mod\Lambda)$ [ST]. Note that $\ttt_S$ and $\ttt_{S'}$ are isomorphic if $S$ and $S'$ are quasi-isomorphic objects in $\kk^{\b}(\mod\Lambda\otimes_RE^{\op})$.

\vskip.5em{\bf Theorem \XFN\ }{\it
For any $n$-spherical object $S\in\dd^{\b}(\mod\Lambda)$, $T:=\ttt_S(\Lambda)$ is a two-sided tilting complex of $\Lambda$ and there is an isomorphism $\ttt_S\simeq(T\Lotimes_\Lambda-)$ of functors on $\dd^{\b}(\mod\Lambda)$.}

\vskip.5em{\sc Proof }
By definition, we have $T=(\Lambda\stackrel{e_\Lambda}{\longrightarrow}\hom^\bullet_{E^{\op}}(S,S))$.
Since $S$ is a complex of $(\Lambda,E)$-modules, $e_\Lambda$ is a chain homomorphism of complexes of $(\Lambda,\Lambda)$-modules. Thus $T$ is a complex of $(\Lambda,\Lambda)$-modules. Since $\ttt_S$ is an autoequivalence, $T$ is a two-sided tilting complex of $\Lambda$. Now, applying $(-\otimes^\bullet_\Lambda P)$ for any $P\in\kk^{\b}(\pr\Lambda)$, we have
$T\otimes^\bullet_\Lambda P=(P\stackrel{e_\Lambda\otimes P}{\longrightarrow}\hom^\bullet_{E^{\op}}(S,S)\otimes^\bullet_\Lambda P)$. We have natural isomorphisms
{\small\begin{eqnarray*}
&&\hom^\bullet_{E^{\op}}(S,S)\otimes^\bullet_\Lambda P=S\otimes^\bullet_E\hom^\bullet_{E^{\op}}(S,E)\otimes^\bullet_\Lambda P\\
&=&S\otimes^\bullet_E\hom^\bullet_{E^{\op}}(\hom^\bullet_\Lambda(P,S),E)=\hom^\bullet_{E^{\op}}(\hom^\bullet_\Lambda(P,S),S),
\end{eqnarray*}}
and one can easily check that the diagram
{\small\[\begin{diag}
T\otimes^\bullet_\Lambda P:\ \ \ &P&\RA{e_\Lambda\otimes P}&\hom^\bullet_{E^{\op}}(S,S)\otimes^\bullet_\Lambda P\\
&\parallel&&\downarrow\wr&&\\
\ttt_S(P):\ \ \ &P&\RA{e_P}&\hom^\bullet_{E^{\op}}(\hom^\bullet_\Lambda(P,S),S)
\end{diag}\]}
of complexes of $\Lambda$-modules commutes. Thus we have a functorial isomorphism $T\otimes^\bullet_\Lambda P\simeq\ttt_S(P)$ on $\kk^{\b}(\pr\Lambda)$.\rule{5pt}{10pt}

\vskip1.5em
{\bf\XG. 3-Calabi-Yau algebras and cluster algebras }

Cluster algebras (with `no coefficients' and in the skew-symmetric case) are completely determined by a finite quiver with no loops or 2-cycles. We show that for quivers of 3-CY algebras $\Lambda$, tilting theory (with tilting modules of projective dimension at most one, which we assume here) is a nice framework for modelling some of the ingredients in the definition of the corresponding cluster algebra. This motivates a closer investigation of tilting modules over 3-CY algebras, which we have already started in previous sections.

Let $B=(b_{ij})$ be an $n\times n$ skew-symmetric matix with integer entries. The Fomin-Zelevinsky mutation $\mu_k$ ($1\le k\le n$) is defined by $\mu_k(B)= (b_{ij}')$, where
{\small\[b_{ij}'=\left\{\begin{array}{cc}
-b_{ij}&\mbox{if $i=k$ or $j=k$}\\
b_{ij}+\frac{|b_{ik}|b_{kj}+b_{ik}|b_{kj}|}{2}&\mbox{otherwise.}
\end{array}\right.\]}
Then $\mu_k(B)$ is skew-symmetric again and satisfies $\mu_k(\mu_k(B))=B$. We identify $B$ with the quiver $\qq$ with vertices $\{1,2,\cdots,n\}$ and $b_{ij}$ arrows from $i$ to $j$ if $b_{ij}>0$. In this way we have a one-one correspondence between skew-symmetric matix with integer entries and finite quivers with no loops or 2-cycles. Thus for a quiver $\qq$ with no loops or 2-cycles, the Fomin-Zelevinsky mutation $\mu_k(\qq)$, which is again a quiver with no loops or 2-cycles, is defined.

Let $R$ be a 3-dimensional complete local Gorenstein ring with an algebraically closed residue field and $\Lambda$ a basic module-finite $R$-algebra which is 3-CY. The valued quiver of $\Lambda$ can be regarded as a (non-valued) quiver. We first show that when the quiver of $\Lambda$ has no loops or 2-cycles, we can interpret the Fomin-Zelevinsky mutation of the quiver via endomorphism rings $\Gamma$ of the non-projective completions of almost complete projective tilting modules. The situation is especially nice if the quiver of $\Gamma$ also has no loops or 2-cycles, so that the procedure can be repeated. We do not know if this is the case in general. But we give examples of where any successive application of our mutation of $\Lambda$ will give algebras with quivers having no loops or 2-cycles. We also show that all algebras obtained from $\Lambda$ via a sequence of mutations can be constructed from a tilting $\Lambda$-module, even in the case when the quiver of $\Lambda$ has loops and/or 2-cycles (see \XGC(2)). If there are no loops or 2-cycles, it is given by a reflexive tilting $\Lambda$-module (see \XGC(4)). Such a result is of interest for the connection with cluster algebras.

Write $\Lambda =\bigoplus_{i=1}^nP^{(i)}$, and let $Q^{(k)}:=\bigoplus_{i\neq k}P^{(i)}$. Denote by $\check{P}^{(k)}\neq P^{(k)}$ the unique indecomposable $\Lambda$-module such that $\nu_k(\Lambda)=Q^{(k)}\oplus \check{P}^{(k)}$ is a tilting module. We put $\mu_k(\Lambda):=\endm_\Lambda(\nu_k(\Lambda))$, which is a 3-CY algebra again by \XCA(1). Assume that the quiver $\qq_{\Lambda}$ of $\Lambda$, with vertices $1, \cdots, n$, has no loops or 2-cycles. The aim of this section is to use tilting theory to obtain a module theoretical interpretation of the Fomin-Zelevinsky mutation of $\qq_{\Lambda}$ at the vertex $k$. We shall compare $\mu_k(\qq_\Lambda)$ with $\qq_{\mu_k(\Lambda)}$.

Recall that for each $i$ we have a minimal projective resolution 
\[0\longrightarrow P^{(i)} \stackrel{f^{(i)}}\longrightarrow P_2^{(i)}{\longrightarrow}P_1^{(i)}\stackrel{g^{(i)}}\longrightarrow P^{(i)} \longrightarrow S^{(i)}\longrightarrow 0\]
of the simple $\Lambda$-module $S^{(i)}$ (\XDC), and $\check{P}^{(i)}:=\Ker g^{(i)}$. The quiver $\qq_{\Lambda}$ is determined by the maps $g^{(i)}:P^{(i)}_1 \to P^{(i)}$. Since by assumption there are no loops, or equivalently, $S_i$ is 3-spherical, then $P^{(i)}_1$ is in $\add Q^{(i)}$. Further $f^{(i)}:P^{(i)}\to P^{(i)}_2 $ is a minimal left $(\add Q^{(i)})$-approximation, and hence gives rise to the arrows in $\qq_{\Lambda}$ starting at $i$. Since there are no 2-cycles, $P^{(i)}_2$ and $P^{(i)}_1$ have no common indecomposable direct summand. Let $b_{ij}$ be the number of arrows from $i$ to $j$ if there are arrows from $i$ to $j$, and otherwise minus the number of arrows from $j$ to $i$. Then if $b_{ij}>0$, it is the multiplicity of $P^{(j)}$ as a summand of $P_1^{(i)}$, or equivalently, the multiplicity of $P^{(i)}$ in $P_2^{(j)}$. Our first goal is to prove the following, which is one of the main results of this section.

\vskip.5em{\bf Theorem \XGA\ }{\it
Let $\Lambda$ be a basic 3-CY algebra. Assume the quiver $\qq_\Lambda$ of $\Lambda$ has no loops or 2-cycles. Then $\mu_k(\qq_\Lambda)$ is obtained from $\qq_{\mu_k(\Lambda)}$ by removing all 2-cycles.}

\vskip.5em{\sc Proof }
Let $B=(b_{ij})$ be the skew symmetric matrix given by the quiver $\qq_{\Lambda}$ and $B'=(b'_{ij}):=\mu_k(B)$. Put $T:=\nu_k(\Lambda)$ and $\Gamma:=\mu_k(\Lambda)=\endm_\Lambda(T)$. Let $B''=(b_{ij}'')$ be the matrix corresponding to $\qq_{\Gamma}$, after we have removed all 2-cycles. Then we want to show that $b_{ij}'=b_{ij}''$. To avoid confusion we denote the vertex in $\qq_{\Gamma}$ corresponding to $k$ by $\check{k}$, and the same for the matrix $B'$.

In order to compare the quivers $\qq_{\Lambda}$ and $\qq_{\Gamma}$, we need to consider (minimal) projective resolutions of simple $\Gamma$-modules and their relationship to the corresponding resolutions of simple $\Lambda$-modules. We divide up into different cases.

{\bf Case 1 }
Assume that $i=k$. We want to show that $b_{j\check{k}}=-b_{jk}$ and $b_{\check{k}j}= -b_{kj}$. We have exact sequences
{\small\begin{eqnarray*}
&0\to (T,P^{(k)})\to (T,P_2^{(k)})\to (T,\check{P}^{(k)})\to \check{S}^{(k)}\to 0,&\\
&0\to (T,\check{P}^{(k)})\to (T,P_1^{(k)})\to (T,P^{(k)})\to \Ext^1_\Lambda(T, \check{P}^{(k)}).&\end{eqnarray*}}
Since $\Ext^1_\Lambda(T,\check{P}^{(k)})=0$, we have an exact sequence of $\Lambda$-modules
{\small\[0\to (T,\check{P}^{(k)})\to (T,P_1^{(k)})\to (T,P_2^{(k)})\to (T,\check{P}^{(k)})\to \check{S}^{(k)}\to 0.\]}
Since  $P_1^{(k)}$ and $P_2^{(k)}$ are in $\add Q^{(k)}$, and hence in $\add T$, and $(T,P_1^{(k)})$ and $(T,P_2^{(k)})$ have no common indecomposable direct summands, we have a minimal projective resolution of the $\Gamma$-module $\check{S}^{(k)}$. 

We want to show that $\check{S}^{(k)}$ is a simple $\Gamma$-module, or equivalently, that there is no loop in $\qq_{\Gamma}$ at the vertex $\check{k}$. For this we need to show that any non-isomorphic map $a:\check{P}^{(k)}\to \check{P}^{(k)}$ factors through an object in $\add Q^{(k)}$. So consider the commutative diagram
{\small\[\begin{diag}
0&\RA{}&\check{P}^{(k)}&\RA{h}&P_1^{(k)}&\RA{g^{(k)}}&P^{(k)}\\
&& \downarrow^{a}&& \downarrow^{b}&& \downarrow^{c}\\
0&\RA{}&\check{P}^{(k)}&\RA{h}&P_1^{(k)}&\RA{g^{(k)}}&P^{(k)}. 
\end{diag}\]}
Here $b$ and $c$ exist by \XDC(5) and (1). If $c$ is an isomorphism, then $b$ is also an isomorphism since $g^{(k)}$ is right minimal. Hence we get the contradiction that $a$ is an isomorphism. Since $c$ is not an isomorphism, $c$ factors through $g^{(k)}$. It is then easy to see that $a$ factors through $h$, where $P_1^{(k)}$ is in $\add Q^{(k)}$. Hence there is no loop at $k$.

We now compare the projective resolutions of $S^{(k)}$ and $\check{S}^{(k)}$. Since $P_1^{(k)}$ and $P_2^{(k)}$ have no common indecomposable direct summands, the same is the case for $(T,P_1^{(k)})$ and $(T,P_2^{(k)})$. Hence we have $b_{j\check{k}}''=-b_{jk}=b_{jk}'$ and $b_{\check{k}j}''=-b_{kj}=b_{\check{k}j}'$.

{\bf Case 2 }
Assume now that $i\neq k$ and $b_{ik}\ge0$, and consider $b_{ij}$ for $j\neq k$. Then $P_2^{(i)}\in\add Q^{(k)}$ and the multiplicity of $P^{(k)}$ in $P_1^{(i)}$ is $m:=b_{ik}$ because there are no cycles of length 2 in $\qq_{\Lambda}$ by assumption. Decompose $P_1^{(i)}=\overline{P}_1^{(i)}\oplus(P^{(k)})^m$ with $\overline{P}_1^{(i)}\in\add Q^{(k)}$. We can take a commutative diagram
{\small\[\begin{diag}
0&\longrightarrow&P^{(i)}&\RA{g^{(i)}}&P_2^{(i)}&\RA{}&\overline{P}_1^{(i)}\oplus(P^{(k)})^m&\RA{f^{(i)}}&P^{(i)}\\
&&\downarrow&&\downarrow&&\downarrow^{{0\choose1}}\\
0&\longrightarrow&(\check{P}^{(k)})^m&\RA{}&(P_1^{(k)})^m&\RA{(f^{(k)})^m}&(P^{(k)})^m
\end{diag}\]}
of exact sequences. This gives rise to the commutative diagram
{\small\[\begin{diag}
0&\longrightarrow&(T,P^{(i)})&\RA{}&(T,P_2^{(i)})&\RA{}&(T,\overline{P}_1^{(i)}\oplus(P^{(k)})^m)&\RA{}&(T,P^{(i)})&\RA{}&\check{S}^{(i)}&\longrightarrow&0\\
&&\downarrow&&\downarrow&&\downarrow^{{0\choose1}}\\
0&\longrightarrow&(T,(\check{P}^{(k)})^m)&\RA{}&(T,(P_1^{(k)})^m)&\RA{}&(T,(P^{(k)})^m)&\RA{}&0
\end{diag}\]}
of exact sequences. Taking the mapping cone, we get a projective resolution
{\small\[0\to(T,P^{(i)})\to(T,P_2^{(i)})\oplus(T,\check{P}^{(k)})^m\to(T,\overline{P}_1^{(i)})\oplus(T,P_1^{(k)})^m\to(T,P^{(i)})\to\check{S}^{(i)}\to0.\]}

\vskip-1em
We now compare the minimal projective resolutions for $S^{(i)}$ and $\check{S}^{(i)}$, to see the change in quivers when passing from $\qq_\Lambda$ to $\qq_\Gamma$, for arrows starting or ending at $i$. Consider now $b_{ij}^{\prime\prime}$ with $j\neq k$. When passing from $b_{ij}$ to $b_{ij}^{\prime\prime}$, we see that we get something extra if and only if $P^{(j)}\in\add P_1^{(k)}$, which is equivalent to $b_{kj}>0$. Then we get $b_{ij}^{\prime\prime}=b_{ij}+b_{kj}m=b_{ij}+b_{ik}b_{kj}=b_{ij}^\prime$. And if $b_{kj}\le0$, we see that $b_{ij}^{\prime\prime}=b_{ij}$. Hence we see that $b_{ij}^{\prime\prime}=b_{ij}^\prime$ for any $j\neq k$.

{\bf Case 3 }
Assume now that $i\neq k$ and $b_{ik}\le0$. We can show $b_{ij}^{\prime\prime}=b_{ij}^\prime$ for any $j\neq k$ by a dual argument to Case 2.\rule{5pt}{10pt}
% and consider $b_{ij}$ for $j\neq k$. We consider the opposite algebra, which is again 3-CY. For $T=Q^{(k)}\oplus\check{P}^{(k)}$, we consider $T^*:=\hom_\Lambda(T,\Lambda)$, which is of the same type.\rule{5pt}{10pt}
%Then $P_1^{(i)}\in\add Q^{(k)}$ and the multiplicity of $P^{(k)}$ in $P_2^{(i)}$ is $n:=-b_{ik}>0$. Consider $b_{ij}$ for $j\neq k$. We consider the opposite algebra, which is again 3-CY. For $T=\check{Q}^{(k)}\oplus\check{P}^{(k)}$, we consider $T^*:=\hom_\Lambda(T,\Lambda)$, which is of the same type.\rule{5pt}{10pt}

\vskip.5em
In order to continue the process, it is of interest to know if the new quiver $\qq_\Gamma$ also has no 2-cycles. Note that we have seen that it has no loops. We do not know if this is true in general, but we show that if $\Lambda=S*G$ where $S=K[[x,y,z]]$, $G\subset\SL_3(K)$ is given by $G=\langle{\rm diag}(\omega,\omega,\omega)\rangle$ with a primitive third root $\omega$ of $1$, then all iterations of the process in \XGA\ give 3-CY algebras whose quivers have no loops. This is a consequence of the following proposition, where we call a ring $\Lambda$ {\it completely graded} if $\Lambda$ is a direct product(!) $\Lambda=\prod_{i\in\zzz}\Lambda_i$ satisfying $\Lambda_i\Lambda_j\subseteq\Lambda_{i+j}$, and a $\Lambda$-module $M$ {\it completely graded} if $M$ is a direct product $M=\prod_{i\in\zzz}M_i$ satisfying $\Lambda_iM_j\subseteq M_{i+j}$. We denote by $\grmod\Lambda$ the category of completely graded $\Lambda$-modules.

\vskip.5em{\bf Proposition \XGB\ }{\it
Let $\Lambda$ be a 3-CY algebra satisfying the conditions (1) and (2) below.

(1) $\Lambda$ is a completely graded ring $\Lambda=\prod_{j\ge0}\Lambda_j$ with Jacobson radical $J_\Lambda=\prod_{j>0}\Lambda_j$. 

(2) $\Lambda=\bigoplus_{i\in\zzz/3\zzz}P^{(i)}$, and each simple $\Lambda$-module $S^{(i)}:=(P^{(i)})_0$ has a projective resolution
$0\to P^{(i)}(-3)\to P^{(i+1)}(-2)^{b_{i+2}}\to P^{(i+2)}(-1)^{b_{i+1}}\to P^{(i)}\to S^{(i)}\to0$ in $\grmod\Lambda$.

Fix $k\in\zzz/3\zzz$ and put $\Gamma:=\mu_k(\Lambda)$. Then $\Gamma$ satisfies the same conditions as $\Lambda$. Precisely speaking, we introduce a degree on $\check{P}^{(k)}:=\Omega^2(S^{(k)})$ such that the natural inclusion $\check{P}^{(k)}\to(P^{(k+2)})^{b_{k+1}}$ is a morphism in $\grmod\Lambda$. Put $\check{P}^{(k+1)}:=P^{(k+1)}$, $\check{P}^{(k+2)}:=P^{(k+2)}(-1)$ and $T:=\bigoplus_{i\in\zzz/3\zzz}\check{P}^{(i)}\in\grmod\Lambda$. Then the assertions ($\check{1}$) and ($\check{2}$) below hold.

($\check{1}$) $\Gamma$ is a completely graded ring $\Gamma=\prod_{j\ge0}\hom_{\grmod\Lambda}(T,T(j))$ with Jacobson radical $J_{\Gamma}=\prod_{j>0}\hom_{\grmod\Lambda}(T,T(j))$.

($\check{2}$) $\Gamma=\bigoplus_{i\in\zzz/3\zzz}(T,\check{P}^{(i)})$, and each simple $\Gamma$-module $\check{S}^{(i)}:=\hom_\Lambda(T,\check{P}^{(i)})_0$ has projective resolution $0\to(T,\check{P}^{(i)})(-3)\to(T,\check{P}^{(i+2)})(-2)^{b'_{i+1}}\to(T,\check{P}^{(i+1)})(-1)^{b'_{i+2}}\to(T,\check{P}^{(i)})\to\check{S}^{(i)}\to0$ in $\grmod\Lambda$ for $b'_k:=b_{k+1}b_{k+2}-b_k$, $b'_{k+1}:=b_{k+1}$ and $b'_{k+2}:=b_{k+2}$.}

\vskip.5em{\sc Proof }
We can prove the assertion by taking care of degree in the proof of \XGA. 

{\bf Case 1 }
We have exact sequences
{\small\begin{eqnarray*}
&0\to\check{P}^{(k)}(-1)\to P^{(k+2)}(-1)^{b_{k+1}}\to P^{(k)}\to S^{(k)}\to0&\\
&0\to P^{(k)}(-3)\to P^{(k+1)}(-2)^{b_{k+2}}\to\check{P}^{(k)}(-1)\to0&
\end{eqnarray*}}
in $\grmod\Lambda$.
%Then we have a complex {\small\[0\to \check{P}^{(k)}(-3)\to \check{P}^{(k+2)}(-2)^{b_{k+1}}\to \check{P}^{(k+1)}(-1)^{b_{k+2}}\to \check{P}^{(k)}\to0.\]}
Applying $\hom_\Lambda(T,-)$ and the same argument as in Case 1 in the proof of \XGA, we obtain a minimal projective resolution
{\small\[0\to(T,\check{P}^{(k)})(-3)\to(T,\check{P}^{(k+2)})(-2)^{b_{k+1}}\to(T,\check{P}^{(k+1)})(-1)^{b_{k+2}}\to(T,\check{P}^{(k)})\to\check{S}^{(k)}\to0\]}
of $S^{(k)}$ in $\grmod\Gamma$.

{\bf Case 2 }
We consider the minimal projective resolution of the simple $\Gamma$-module $S^{(k+1)}$. We have a commutative diagram
{\small\[\begin{diag}
0&\longrightarrow&P^{(k+1)}(-3)&\RA{}&P^{(k+2)}(-2)^{b_k}&\RA{}&P^{(k)}(-1)^{b_{k+2}}&\RA{}&P^{(k+1)}\\
&&\downarrow&&\downarrow^{f(-2)}&&\parallel\\
0&\longrightarrow&\check{P}^{(k)}(-2)^{b_{k+2}}&\RA{}&P^{(k+2)}(-2)^{b_{k+1}b_{k+2}}&\RA{}&P^{(k)}(-1)^{b_{k+2}}
\end{diag}\]}
of exact sequences. We can choose $f$ to be a morphism in $\grmod\Lambda$. Putting $X:=(P^{(k+2)})^{b_k}$ and $Y:=(P^{(k+2)})^{b_{k+1}b_{k+2}}$ and looking at the degree two part, we have a commutative diagram
{\small\[\begin{diag}
P^{(k+1)}_{-1}=0&\RA{}&X_0&\RA{}&(P^{(k)}_1)^{b_{k+2}}\\
&&\downarrow^{f_0}&&\parallel\\
&&Y_0&\RA{}&(P^{(k)}_1)^{b_{k+2}}
\end{diag}\]}
of exact sequences. Thus $f_0:X_0\to Y_0$ is a monomorphism between semisimple $\Lambda$-modules because $J_{\Lambda}=\prod_{i>0}\Lambda_i$ by our assumption. We can take $g\in\hom_{\grmod\Lambda}(Y,X)$ such that $f_0g_0=1_{X_0}$. Then $fg-1\in\endm_{\grmod\Lambda}(X)$ satisfies $X_0\subset\Ker(fg-1)$. Since $X$ is generated by $X_0$, we obtain $fg-1=0$. Thus $f$ is a split monomorphism 

Consequently, taking a mapping cone as in Case 2 in the proof of \XGA\ and cancelling a trivial direct summand of the complex, we have a complex
\[0\to P^{(k+1)}(-3)\to\check{P}^{(k)}(-2)^{b_{k+2}}\to P^{(k+2)}(-2)^{b_{k+1}b_{k+2}-b_k}\to P^{(k+1)}\]
in $\grmod\Lambda$, which induces a projective resolution
\[0\to(T,\check{P}^{(k+1)})(-3)\to(T,\check{P}^{(k)})(-2)^{b'_{k+2}}\to(T,\check{P}^{(k+2)})(-1)^{b'_k}\to(T,\check{P}^{(k+1)})\to\check{S}^{(k+1)}\to0\]
of $\check{S}^{(k+1)}$ in $\grmod\Gamma$.

{\bf Case 3 }
A dual argument works for $\check{S}^{(k+2)}$.\rule{5pt}{10pt}.

\vskip.5em
Now we consider all iterations of $S*G$ for $G=\langle{\rm diag}(\omega,\omega,\omega)\rangle\subset\SL_3(K)$ with $\omega^3=1$. The Fomin-Zelevinsky mutation of the quiver
\begin{picture}(32,25)
\put(-2,-7){\small$k$}
\put(0,0){\circle*{4}}
\put(32,0){\circle*{4}}
\put(16,24){\circle*{4}}
\put(14,0){\scriptsize $b$}
\put(26,10){\scriptsize $c$}
\put(0,10){\scriptsize $a$}
\put(0,0){\vector(1,0){30}}
\put(30,3){\vector(-2,3){13}}
\put(14,22){\vector(-2,-3){13}}
\end{picture}
at the vertex $k$ is 
\begin{picture}(32,25)
\put(-2,-7){\small $k$}
\put(0,0){\circle*{4}}
\put(32,0){\circle*{4}}
\put(16,24){\circle*{4}}
\put(14,0){\scriptsize $b$}
\put(26,10){\scriptsize $ab-c$}
\put(0,10){\scriptsize $a$}
\put(30,0){\vector(-1,0){30}}
\put(17,22.5){\vector(2,-3){13}}
\put(1,2.5){\vector(2,3){13}}
\end{picture}\ , where $a$, $b$, $c$ and $ab-c$ show the numbers of arrow. The quiver of $S*G$ is the McKay quiver 
\begin{picture}(32,25)
\put(0,0){\circle*{4}}
\put(32,0){\circle*{4}}
\put(16,24){\circle*{4}}
\put(14,0){\scriptsize 3}
\put(26,10){\scriptsize 3}
\put(0,10){\scriptsize 3}
\put(0,0){\vector(1,0){30}}
\put(30,3){\vector(-2,3){13}}
\put(14,22){\vector(-2,-3){13}}
\end{picture}\ of $G$ (see section \XC), and we draw quivers of 3-CY algebras obtained by iterated mutations. It is the picture below, where each quiver has precisely three neighbours. One can check inductively that each triple $(a,b,c)$ satisfies the Markov equation $a^2+b^2+c^2=abc$. All integral solutions of the Markov equation appear in the picture below because it is known that the Fomin-Zelevinky mutation rule $(a,b,c)\mapsto(a,b,ab-c)$ in this case gives all of them.
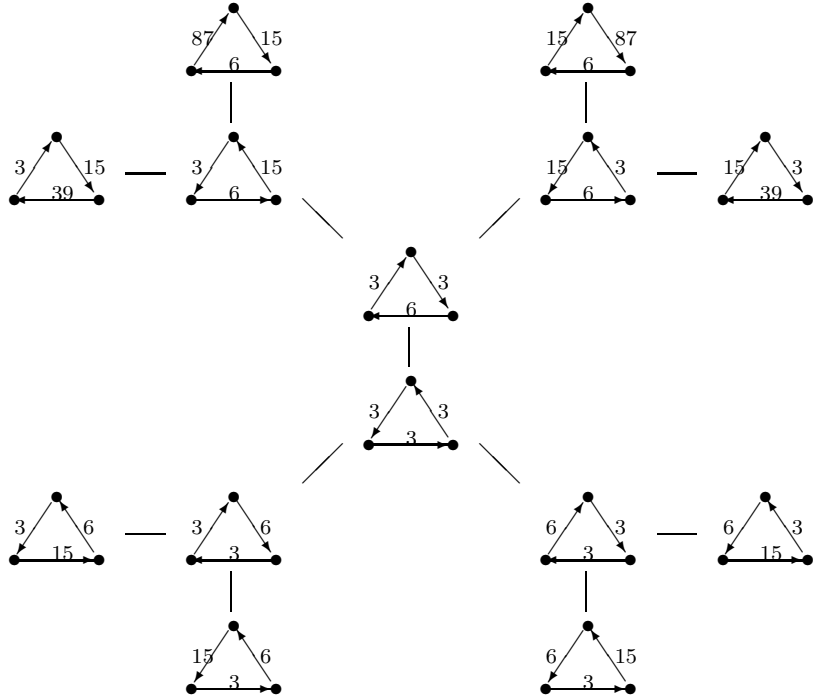
\[\begin{array}{ccccccccc}
&&\begin{picture}(32,25)
%\put(-7,-7){\scriptsize 1}
%\put(32,-7){\scriptsize 2}
%\put(14,28){\scriptsize 3}
\put(0,0){\circle*{4}}
\put(32,0){\circle*{4}}
\put(16,24){\circle*{4}}
\put(14,0){\scriptsize 6}
\put(26,10){\scriptsize 15}
\put(0,10){\scriptsize 87}
\put(30,0){\vector(-1,0){30}}
\put(17,22.5){\vector(2,-3){13}}
\put(1,2.5){\vector(2,3){13}}
\end{picture}
&&&&
\begin{picture}(32,25)
\put(0,0){\circle*{4}}
\put(32,0){\circle*{4}}
\put(16,24){\circle*{4}}
\put(14,0){\scriptsize 6}
\put(26,10){\scriptsize 87}
\put(0,10){\scriptsize 15}
\put(30,0){\vector(-1,0){30}}
\put(17,22.5){\vector(2,-3){13}}
\put(1,2.5){\vector(2,3){13}}
\end{picture}&\\

&&\begin{picture}(2,15)\put(0,0){\line(0,1){15}}\end{picture}&&
&&\begin{picture}(2,15)\put(0,0){\line(0,1){15}}\end{picture}&&\\

\begin{picture}(32,25)
\put(0,0){\circle*{4}}
\put(32,0){\circle*{4}}
\put(16,24){\circle*{4}}
\put(14,0){\scriptsize 39}
\put(26,10){\scriptsize 15}
\put(0,10){\scriptsize 3}
\put(30,0){\vector(-1,0){30}}
\put(17,22.5){\vector(2,-3){13}}
\put(1,2.5){\vector(2,3){13}}
\end{picture}&
\begin{picture}(15,5)\put(0,10){\line(1,0){15}}\end{picture}&
\begin{picture}(32,25)
\put(0,0){\circle*{4}}
\put(32,0){\circle*{4}}
\put(16,24){\circle*{4}}
\put(14,0){\scriptsize 6}
\put(26,10){\scriptsize 15}
\put(0,10){\scriptsize 3}
\put(0,0){\vector(1,0){30}}
\put(30,3){\vector(-2,3){13}}
\put(14,22){\vector(-2,-3){13}}
\end{picture}&&&&
\begin{picture}(32,25)
\put(0,0){\circle*{4}}
\put(32,0){\circle*{4}}
\put(16,24){\circle*{4}}
\put(14,0){\scriptsize 6}
\put(26,10){\scriptsize 3}
\put(0,10){\scriptsize 15}
\put(0,0){\vector(1,0){30}}
\put(30,3){\vector(-2,3){13}}
\put(14,22){\vector(-2,-3){13}}
\end{picture}&
\begin{picture}(15,5)\put(0,10){\line(1,0){15}}\end{picture}&
\begin{picture}(32,25)
\put(0,0){\circle*{4}}
\put(32,0){\circle*{4}}
\put(16,24){\circle*{4}}
\put(14,0){\scriptsize 39}
\put(26,10){\scriptsize 3}
\put(0,10){\scriptsize 15}
\put(30,0){\vector(-1,0){30}}
\put(17,22.5){\vector(2,-3){13}}
\put(1,2.5){\vector(2,3){13}}
\end{picture}\\

&&&\begin{picture}(15,10)\put(15,0){\line(-1,1){15}}\end{picture}&
&\begin{picture}(15,10)\put(0,0){\line(1,1){15}}\end{picture}&&&\\

&&&&\begin{picture}(32,25)
\put(0,0){\circle*{4}}
\put(32,0){\circle*{4}}
\put(16,24){\circle*{4}}
\put(14,0){\scriptsize 6}
\put(26,10){\scriptsize 3}
\put(0,10){\scriptsize 3}
\put(30,0){\vector(-1,0){30}}
\put(17,22.5){\vector(2,-3){13}}
\put(1,2.5){\vector(2,3){13}}
\end{picture}&&&&\\

&&&&\begin{picture}(2,15)\put(0,0){\line(0,1){15}}\end{picture}&&&&\\

&&&&\begin{picture}(32,25)
\put(0,0){\circle*{4}}
\put(32,0){\circle*{4}}
\put(16,24){\circle*{4}}
\put(14,0){\scriptsize 3}
\put(26,10){\scriptsize 3}
\put(0,10){\scriptsize 3}
\put(0,0){\vector(1,0){30}}
\put(30,3){\vector(-2,3){13}}
\put(14,22){\vector(-2,-3){13}}
\end{picture}&&&&\\

&&&\begin{picture}(15,10)\put(0,0){\line(1,1){15}}\end{picture}&
&\begin{picture}(15,10)\put(15,0){\line(-1,1){15}}\end{picture}&&&\\

\begin{picture}(32,25)
\put(0,0){\circle*{4}}
\put(32,0){\circle*{4}}
\put(16,24){\circle*{4}}
\put(14,0){\scriptsize 15}
\put(26,10){\scriptsize 6}
\put(0,10){\scriptsize 3}
\put(0,0){\vector(1,0){30}}
\put(30,3){\vector(-2,3){13}}
\put(14,22){\vector(-2,-3){13}}
\end{picture}&
\begin{picture}(15,5)\put(0,10){\line(1,0){15}}\end{picture}&
\begin{picture}(32,25)
\put(0,0){\circle*{4}}
\put(32,0){\circle*{4}}
\put(16,24){\circle*{4}}
\put(14,0){\scriptsize 3}
\put(26,10){\scriptsize 6}
\put(0,10){\scriptsize 3}
\put(30,0){\vector(-1,0){30}}
\put(17,22.5){\vector(2,-3){13}}
\put(1,2.5){\vector(2,3){13}}
\end{picture}
&&&&\begin{picture}(32,25)
\put(0,0){\circle*{4}}
\put(32,0){\circle*{4}}
\put(16,24){\circle*{4}}
\put(14,0){\scriptsize 3}
\put(26,10){\scriptsize 3}
\put(0,10){\scriptsize 6}
\put(30,0){\vector(-1,0){30}}
\put(17,22.5){\vector(2,-3){13}}
\put(1,2.5){\vector(2,3){13}}
\end{picture}&
\begin{picture}(15,5)\put(0,10){\line(1,0){15}}\end{picture}&
\begin{picture}(32,25)
\put(0,0){\circle*{4}}
\put(32,0){\circle*{4}}
\put(16,24){\circle*{4}}
\put(14,0){\scriptsize 15}
\put(26,10){\scriptsize 3}
\put(0,10){\scriptsize 6}
\put(0,0){\vector(1,0){30}}
\put(30,3){\vector(-2,3){13}}
\put(14,22){\vector(-2,-3){13}}
\end{picture}\\

&&\begin{picture}(2,15)\put(0,0){\line(0,1){15}}\end{picture}&&
&&\begin{picture}(2,15)\put(0,0){\line(0,1){15}}\end{picture}&&\\

&&\begin{picture}(32,25)
\put(0,0){\circle*{4}}
\put(32,0){\circle*{4}}
\put(16,24){\circle*{4}}
\put(14,0){\scriptsize 3}
\put(26,10){\scriptsize 6}
\put(0,10){\scriptsize 15}
\put(0,0){\vector(1,0){30}}
\put(30,3){\vector(-2,3){13}}
\put(14,22){\vector(-2,-3){13}}
\end{picture}&&&&
\begin{picture}(32,25)
\put(0,0){\circle*{4}}
\put(32,0){\circle*{4}}
\put(16,24){\circle*{4}}
\put(14,0){\scriptsize 3}
\put(26,10){\scriptsize 15}
\put(0,10){\scriptsize 6}
\put(0,0){\vector(1,0){30}}
\put(30,3){\vector(-2,3){13}}
\put(14,22){\vector(-2,-3){13}}
\end{picture}&&
\end{array}\]

\vskip.5em
The above picture gives the Hasse graph of $\tilt_1(S*G)$ (see \XGC(1) below). We note that a similar picture appeared in the classification of exceptional vector bundles over $\ppp^2$ due to Gorodentsev-Rudakov [GR][Rud]. It also appeared in recent work of Bridgeland [Bri2,3] on t-structures on the derived category of the total space of the canonical line bundle $\oo_{\ppp^2}(-3)$ on $\ppp^2$.

Now we consider other 3-CY algebras $S*G$ for $G=\langle{\rm diag}(\zeta,\zeta^2,\zeta^2)\rangle\subset\SL_3(K)$ with $\zeta^5=1$. We draw a few quivers of 3-CY algebras obtained by iterated mutations. Unfortunately we do not know whether they coincide with mutations of algebras.
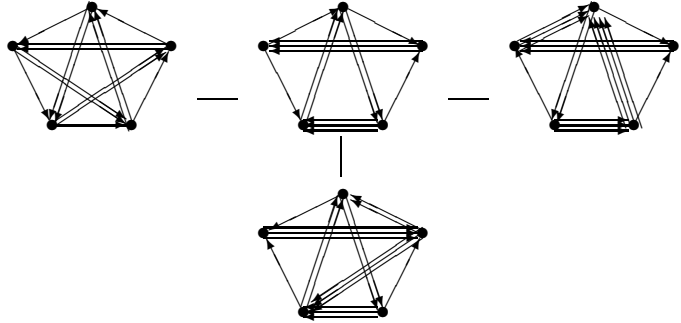
\[\begin{array}{ccccc}
\begin{picture}(60,47)
\put(15,0){\circle*{4}}
\put(45,0){\circle*{4}}
\put(60,30){\circle*{4}}
\put(30,45){\circle*{4}}
\put(0,30){\circle*{4}}
\put(15,0){\vector(1,0){28}}
\put(45,0){\vector(1,2){14}}
\put(60,30){\vector(-2,1){28}}
\put(30,45){\vector(-2,-1){28}}
\put(0,30){\vector(1,-2){14}}
\put(15,-1){\vector(3,2){43}}
\put(15,1){\vector(3,2){43}}
\put(45,2){\vector(-1,3){14}}
\put(44,-2){\vector(-1,3){15}}
\put(60,29){\vector(-1,0){58}}
\put(60,31){\vector(-1,0){58}}
\put(30,43){\vector(-1,-3){14}}
\put(29,47){\vector(-1,-3){15}}
\put(0,29){\vector(3,-2){43}}
\put(0,31){\vector(3,-2){43}}
\end{picture}&
\begin{picture}(15,5)\put(0,10){\line(1,0){15}}\end{picture}&
\begin{picture}(60,47)
\put(15,0){\circle*{4}}
\put(45,0){\circle*{4}}
\put(60,30){\circle*{4}}
\put(30,45){\circle*{4}}
\put(0,30){\circle*{4}}
\put(43,-2){\vector(-1,0){28}}
\put(43,0){\vector(-1,0){28}}
\put(43,2){\vector(-1,0){28}}
\put(45,0){\vector(1,2){14}}
\put(32,44){\vector(2,-1){27}}
\put(2,31){\vector(2,1){27}}
\put(0,30){\vector(1,-2){14}}
\put(31,44){\vector(1,-3){14}}
\put(29,43){\vector(1,-3){14}}
\put(60,28){\vector(-1,0){58}}
\put(60,30){\vector(-1,0){58}}
\put(60,32){\vector(-1,0){58}}
\put(16,1){\vector(1,3){14}}
\put(14,2){\vector(1,3){14}}
\end{picture}&
\begin{picture}(15,5)\put(0,10){\line(1,0){15}}\end{picture}&
\begin{picture}(60,47)
\put(15,0){\circle*{4}}
\put(45,0){\circle*{4}}
\put(60,30){\circle*{4}}
\put(30,45){\circle*{4}}
\put(0,30){\circle*{4}}
\put(15,-2){\vector(1,0){28}}
\put(15,0){\vector(1,0){28}}
\put(15,2){\vector(1,0){28}}
\put(45,0){\vector(1,2){14}}
\put(32,44){\vector(2,-1){27}}
\put(1,32){\vector(2,1){27}}
\put(1,30){\vector(2,1){27}}
\put(1,28){\vector(2,1){27}}
\put(14,2){\vector(-1,2){14}}

\put(42,-1){\vector(-1,3){14}}
\put(44,-1){\vector(-1,3){14}}
\put(46,-1){\vector(-1,3){14}}
\put(48,-1){\vector(-1,3){14}}
\put(60,28){\vector(-1,0){58}}
\put(60,30){\vector(-1,0){58}}
\put(60,32){\vector(-1,0){58}}
\put(30,43){\vector(-1,-3){14}}
\put(29,47){\vector(-1,-3){15}}
\end{picture}\\
&&\begin{picture}(2,15)\put(0,0){\line(0,1){15}}\end{picture}&&\\
&&\begin{picture}(60,47)
\put(15,0){\circle*{4}}
\put(45,0){\circle*{4}}
\put(60,30){\circle*{4}}
\put(30,45){\circle*{4}}
\put(0,30){\circle*{4}}
\put(43,-2){\vector(-1,0){28}}
\put(43,0){\vector(-1,0){28}}
\put(43,2){\vector(-1,0){28}}
\put(45,0){\vector(1,2){14}}
\put(60,29){\vector(-2,1){27}}
\put(60,31){\vector(-2,1){27}}
\put(30,45){\vector(-2,-1){28}}
\put(14,2){\vector(-1,2){13}}
\put(60,28){\vector(-3,-2){42}}
\put(60,30){\vector(-3,-2){42}}
\put(60,32){\vector(-3,-2){42}}
\put(0,28){\vector(1,0){58}}
\put(0,30){\vector(1,0){58}}
\put(0,32){\vector(1,0){58}}
\put(31,44){\vector(1,-3){14}}
\put(29,43){\vector(1,-3){14}}
\put(16,1){\vector(1,3){14}}
\put(14,2){\vector(1,3){14}}
\end{picture}&&
\end{array}\]

%\begin{picture}(60,47)\put(15,0){\circle*{4}}\put(45,0){\circle*{4}}\put(60,30){\circle*{4}}\put(30,45){\circle*{4}}\put(0,30){\circle*{4}}\put(44,0){\vector(-1,0){28}}\put(45,0){\vector(1,2){14}}\put(60,32){\vector(-2,1){28}}\put(60,30){\vector(-2,1){28}}\put(60,28){\vector(-2,1){28}}\put(30,45){\vector(-2,-1){28}}\put(14,2){\vector(-1,2){13}}\put(58,28){\vector(-3,-2){42}}\put(58,30){\vector(-3,-2){43}}\put(16,1){\vector(1,3){14}}\put(14,2){\vector(1,3){14}}\put(0,29){\vector(3,-2){43}}\put(0,31){\vector(3,-2){43}}\put(0,33){\vector(3,-2){43}}\end{picture}

Motivated by the connection with cluster algebras, we would like to view the tilting modules as analogs to clusters, and hence we want to show that any 3-CY algebra $\Gamma$ obtained via a sequence of mutations from a 3-CY algebra $\Lambda$ can be obtained directly as the endomorphism ring of a tilting $\Lambda$-module. We already know from section \XE\ that if $T=\bigoplus_{i=1}^nT_i$ is a basic tilting $\Lambda$-module, with the $T_i$ indecomposable, then there is for each $k=1,\cdots,n$ a unique indecomposable $\Lambda$-module $\check{T}_k$ with $\check{T}_k\ {\not\simeq}\ T_k$ such that $\nu_k(T):=(\bigoplus_{i\neq k}T_i)\oplus\check{T}_k$ is a tilting module, and we have also given an explicit description of $\nu_k(T)$.

We now point out the relationship between mutation of algebras and of tilting modules in part (1) of the next result. Part (3) is analogous to a basic property of the Fomin-Zelevinsky mutation. We know that all algebras obtained from a 3-CY algebra $\Lambda$ by successive applications of mutations are derived equivalent, and hence by [Ri1] can be obtained from $\Lambda$ as an endomorphism ring of a tilting complex. But we have a better result (2) saying that if $\Gamma$ is obtained from $\Lambda$ by a sequence of our special tilting modules, then it can be obtained directly with one tilting module. A more general version of part (4) (including an alternative approach) is given in the next section \XHK(3).

\vskip.5em{\bf Proposition \XGC\ }{\it
Let $\Lambda$ be basic 3-CY.

(1) $\mu_k(\endm_\Lambda(T))=\endm_\Lambda(\nu_k(T))$ for any tilting $\Lambda$-module $T$.

(2) $\mu_{k_m}\circ\cdots\circ\mu_{k_1}(\Lambda)=\endm_\Lambda(\nu_{k_m}\circ\cdots\circ\nu_{k_1}(\Lambda))$ for any $k_1,\cdots,k_m$.

(3) $\mu_k(\mu_k(\Lambda)))=\Lambda$.

(4) Assume in (2) that the quivers of $\Lambda_i:=\mu_{k_i}\circ\cdots\circ\mu_{k_1}(\Lambda)$ have no loops for any $i$ ($0\le i<m$). Then $\nu_{k_m}\circ\cdots\circ\nu_{k_1}(\Lambda)$ is reflexive and isomorphic to $(T_1\otimes_{\Lambda_1}T_2\otimes_{\Lambda_2}\cdots\otimes_{\Lambda_{m-1}}T_m)^{**}$ for $T_i:=\nu_{k_i}(\Lambda_{i-1})$.}

\vskip.5em{\sc Proof }
(1) follows from \XEF(6). Using (1) repeatedly, we have (2). We obtain (3) by $\mu_k(\mu_k(\Lambda))=\endm_\Lambda(\nu_k(\nu_k(\Lambda))=\endm_\Lambda(\Lambda)=\Lambda$. We obtain (4) by using \XEF(5) repeatedly.\rule{5pt}{10pt}

\vskip.5em
In addition to considering, for some fixed 3-CY algebra $\Lambda$, the tilting modules to be the analogs of the clusters and the indecomposable partial tilting modules $M$ as the analogs of the cluster variables,
we also have an interpretation of the exchange multiplication rule. Let $T=\bigoplus_{i=1}^nT_i$ be a tilting module, and assume that $T$ is connected with the cluster $\underline{x}=\{x_1,\cdots,x_n\}$ in such a way that the $T_i$ are associated with the cluster variables $x_i$. (If we start by fixing a correspondence $\underline{u}=\{u_1,\cdots,u_n\}\mapsto\Lambda=\bigoplus_{i=1}^nQ_i$,
we reach $\Lambda^\prime=\endm_\Lambda(T)$ by a sequence of mutations. We choose $\underline{x}=\{x_1\cdots,x_n\}$ to be the cluster obtained by the same sequence of mutations applied to $\underline{u}=\{u_1,\cdots,u_n\}$. This is the procedure used in [BMR2] in the context of cluster categories.
For $\overline{T}=T/T_k$, we have the minimal right and left $\add\overline{T}$-approximations $B\to T_k$ and $T_k\to B^\prime$, with $B=\bigoplus_{i\neq k}T_i^{r_i}$ and $B^\prime=\bigoplus_{i\neq k}T_i^{s_i}$, where $r_i\ge0$ and $s_i\ge0$. We then have $T_k\cdot\check{T}_k=\prod T_i^{r_i}+\prod T_i^{s_i}$.

As we asked in section \XE, it would be interesting to know if every tilting module $T$ can be obtained from the tilting module $\Lambda$ over a 3-CY algebra $\Lambda$ by a finite sequence of mutations, that is, if the Hasse quiver of $\tilt_1\Lambda$ is connected. Another interesting problem is whether there is a one-one correspondence between the cluster variables and the indecomposable partial tilting $\Lambda$-modules, inducing a one-one correspondence between clusters and tilting modules.

In addition to the work on cluster categories serving as a model, there is also a connection with the modelling of some class of cluster algebras `with coefficients' by special modules over preprojective algebras, from [GLS]. They consider maximal rigid modules $M$ (see section \XH\ for definition) over a preprojective algebra $\Lambda$ of a Dynkin diagram. For the case of $\Lambda$ being of finite representation type, the stable category $\underline{\mod}\Lambda$ is actually equivalent to a cluster category [BMRRT], and the maximal rigid modules are closely related to the cluster-tilting objects, as well as coinciding with the maximal $1$-orthogonal modules of [I3,4]. Here $\endm_\Lambda(M)$ has global dimension $3$, and the tilting theory over this algebra is relevant. These algebras are sort of a degenerate version of 3-CY algebras. The relationship is similar to  the relationship between the invariant ring $S^G$ and the skew group ring $S*G$ where $S=K[[x_1,\cdots,x_n]]$ and $G$ is a finite subgroup of $\SL_3(k)$.

\vskip1.5em{\bf\XH. Non-commutative crepant resolutions }

In this section we improve results on tilting modules for 3-CY algebras from section \XG. We show that if we can pass from $\Lambda$ to $\Lambda'$ by taking the endomorphism algebra of a tilting module, and the same way from $\Lambda'$ to $\Lambda''$, then we can also pass directly from $\Lambda$ to $\Lambda''$ in this way. For this, constructing new tilting modules from old ones via homomorphism spaces and tensor products is crucial. Actually, we work in a more general context, investigating the non-commutative crepant resolutions (NCCR) of Van den Bergh, extending his definition to non-commutative algebras. For 3-CY algebras there turns out to be a close relationship to reflexive tilting modules. A main result is the solution of a conjecture of Van den Bergh on derived equivalence of NCCR for 3-dimensional algebras.
%When an algebra $\Lambda'$ is obtained from an algebra $\Lambda$ as the endomorphism algebra of a tilting module, and $\Lambda''$ is obtained from $\Lambda'$ the same way, then $\Lambda$ and $\Lambda''$ are derived equivalent. But in general $\Lambda''$ is not isomorphic to the endomorphism algebra of a tilting $\Lambda$-module. We can however prove such a result for 3-CY algebras, and this is useful for the connection with cluster algebras. Acturally, we work in a more general context, involving the non-commutative crepant resolutions (NCCR) of Van den Bergh, and prove a conjecture of Van den Bergh on 3-dimensional module-finite algebras. We investigate the relationship between NCCR and tilting modules for 3-CY and 3-CY$^-$ algebras.

Throughout this section, let $R$ be a normal Gorenstein domain with $\dim R=d$ and $\Lambda$ a module-finite $R$-algebra such that the structure morphism $R\to\Lambda$ is injective.
Generalizing the following definition of Van den Bergh [Va1,2], we say that $M$ gives a {\it non-commutative crepant resolution} (NCCR for short) $\Gamma:=\endm_\Lambda(M)$ of $\Lambda$ if 

(1) $M\in\ref\Lambda$ is a height one generator (\XBD) of $\Lambda$, and

(2) $\Gamma_{\dn{p}}$ is an $R_{\dn{p}}$-order (in the sense of section \XB) with $\gl\Gamma_{\dn{p}}=\height\dn{p}$ for any $\dn{p}\in\Spec R$.

Obviously, one can replace $\Spec R$ by $\Max R$ in condition (2). If $\Gamma$ is $d$-CY, then condition (2) is satisfied by \XCB.

In [Va1,2], Van den Bergh gave a non-commutative analogue of a conjecture of Bondal-Orlov [BO]: {\it All NCCR of a normal Gorenstein domain $\Lambda$ are derived equivalent}. In fact, he proved this conjecture for 3-dimensional terminal singularities $\Lambda$. In this section, using a method from [I4], we show that {\it his conjecture is true for arbitrary 3-dimensional module-finite algebras}.
In the original definition in [Va1,2], the height one generator condition is not assumed. If $\Lambda$ is a normal domain, then any non-zero reflexive $\Lambda$-module is a height one generator. Thus our definition coincides with the original one in this case. For the non-commutative situation, if we drop the height one generator condition, then $\Lambda:=\def\arraystretch{.5}\left(\begin{array}{cc}{\scriptstyle R}&{\scriptstyle R}\\ {\scriptstyle xR}&{\scriptstyle R}\end{array}\right)$ with $R:=K[[x,y,z]]$ has NCCR $\Lambda=\endm_\Lambda(\Lambda)$ and $R=\endm_\Lambda(\def\arraystretch{.5}\left(\begin{array}{c}{\scriptstyle R}\\ {\scriptstyle R}\end{array}\right))$, and $\Lambda$ and $R$ are not derived equivalent. Thus the height one generator condition seems to be appropriate for our considerations.

For NCCR it would be interesting to know for which $R$-algebras $\Lambda$ they exist, what the relationship is with $\Lambda$ when they exist, and what the connection is between different NCCR of the same algebra $\Lambda$. Also, it would be nice to describe the $\Lambda$-modules giving rise to some NCCR. Some such questions will be investigated along the way.

We start with some easy properties of NCCR, which will be useful later. In particular, reflexive equivalences play a crucial role.

\vskip.5em{\bf Proposition \XHA\ }{\it
We have the following for a module-finite $R$-algebra $\Lambda$.

(1) Assume that $M\in\ref\Lambda$ gives a NCCR of $\Lambda$.

\strut\kern1em
(i) $M_{\dn{p}}\in\ref\Lambda_{\dn{p}}$ gives a NCCR of $\Lambda_{\dn{p}}$ for any $\dn{p}\in\Spec R$.

\strut\kern1em
(ii) $\Lambda_{\dn{p}}$ is an $R_{\dn{p}}$-order with $\gl\Lambda_{\dn{p}}=1$ for any $\dn{p}\in\Spec R$ with $\height\dn{p}=1$.

\strut\kern1em
(iii) $M$ is a height one progenerator of $\Lambda$, and we have a reflexive equivalence $\fff:=\hom_\Lambda(M,-):\ref\Lambda\to\ref\endm_\Lambda(M)$.

(2) Let $\fff:\ref\Lambda\to\ref\Gamma$ be any reflexive equivalence. Then $N\in\ref\Lambda$ gives a NCCR of $\Lambda$ if and only if $\fff(N)\in\ref\Gamma$ gives a NCCR of $\Gamma$.}

\vskip.5em{\sc Proof }
(1)(i) Obvious from the definition of NCCR.
%Obviously $M_{\dn{p}}$ is a height one generator. Since $\Gamma:=\endm_\Lambda(M)$ is an $R$-order with $\gl\Gamma=d$, $\Gamma_{\dn{p}}$ is an $R_{\dn{p}}$-order with $\gl\Gamma_{\dn{p}}=\height\dn{p}$, and $\Gamma_{\dn{p}}=\endm_{\Lambda_{\dn{p}}}(M_{\dn{p}})$.

(ii) Assume that $M$ gives a NCCR $\Gamma=\endm_\Lambda(M)$. Since $M_{\dn{p}}$ is a generator for $\Lambda_{\dn{p}}$, there exists an idempotent $e$ in $\ma_n(\Gamma_{\dn{p}})$ for some $n$ such that $\Lambda_{\dn{p}}$ is Morita equivalent to $e\ma_n(\Gamma_{\dn{p}})e$. Since $\ma_n(\Gamma_{\dn{p}})$ is an $R_{\dn{p}}$-order with $\gl\ma_n(\Gamma_{\dn{p}})=1$, we have that $\Lambda_{\dn{p}}$ is an $R_{\dn{p}}$-order with $\gl\Lambda_{\dn{p}}=1$ by taking completion and using \XCL(3).

(iii) Fix any $\dn{p}\in\Spec R$ with $\height\dn{p}=1$. By (ii), we have $\gl\Lambda_{\dn{p}}=1$. Since $M_{\dn{p}}$ is a torsionfree $R_{\dn{p}}$-module, $M_{\dn{p}}$ is a projective $\Lambda_{\dn{p}}$-module. Thus the assertion follows by \XBD.

(2) Since $\endm_\Lambda(N)=\endm_\Gamma(\fff(N))$ holds, we only have to show that $N$ is a height one generator of $\Lambda$ if and only if $\fff(N)$ is a height one generator of $\Gamma$. Fix $\dn{p}\in\Spec R$ with $\height\dn{p}=1$. Then $\fff$ induces an equivalence $\fff_{\dn{p}}:\ref\Lambda_{\dn{p}}\to\ref\Gamma_{\dn{p}}$ by \XBD(2). Since $\gl\Lambda_{\dn{p}}=\gl\Gamma_{\dn{p}}=1$ by (1)(ii), $N_{\dn{p}}$ is a generator if and only if $\ref\Lambda_{\dn{p}}=\add N_{\dn{p}}$ if and only if $\ref\Gamma_{\dn{p}}=\add\fff_{\dn{p}}(N_{\dn{p}})$ if and only if $\fff_{\dn{p}}(N_{\dn{p}})$ is a generator. Thus the assertion follows.\rule{5pt}{10pt}

\vskip.5em
We note here the following necessary (and sufficient for $d\le2$) conditions for $\Lambda$ to have some NCCR even though we do not use it in the rest of this paper.

\vskip.5em{\bf Proposition \XHB\ }{\it
If $\Lambda$ has a NCCR, then (1) and (2) hold. The converse holds if $d\le2$.

(1) $\Lambda_{\dn{p}}$ is an $R_{\dn{p}}$-order with $\gl\Lambda_{\dn{p}}=\height\dn{p}$ for any $\dn{p}\in\Spec R$ with $\height\dn{p}\le1$.

(2) There exists $M\in\ref\Lambda$ such that $\ref\Lambda_{\dn{p}}=\add M_{\dn{p}}$ for any $\dn{p}\in\Spec R$ with $\height\dn{p}=2$. In particular, $\Lambda_{\dn{p}}$ is representation-finite.}

\vskip.5em{\sc Proof }
(1) follows from \XHA(1)(ii). Since $\Gamma_{\dn{p}}$ is a NCCR of $\Lambda_{\dn{p}}$ by \XHA(1)(i), we have a reflexive equivalence $\hom_{\Lambda_{\dn{p}}}(M_{\dn{p}},-):\ref\Lambda_{\dn{p}}\to\ref\Gamma_{\dn{p}}$ by \XHA(1)(iii). Since $\gl\Gamma_{\dn{p}}=2$, $\ref\Gamma_{\dn{p}}$ consists of projective $\Gamma_{\dn{p}}$-modules. Thus $\ref\Lambda_{\dn{p}}=\add M_{\dn{p}}$.

We now show the converse if $d\le2$. If $d\le 1$, then $\Lambda$ itself gives a NCCR. If $d=2$, then we take $M$ in (2). Since $\ref\Lambda_{\dn{p}}$ is closed under kernels, we have that $\gl\endm_{\Lambda_{\dn{p}}}(M_{\dn{p}})=2$ [A1]. Thus $M$ gives a NCCR of $\Lambda$.\rule{5pt}{10pt}

\vskip.5em
For algebras which are $d$-CY$^-$ or $d$-CY, we have some nice properties.

\vskip.5em{\bf Proposition \XHC\ }{\it
(1) If a $n$-CY$^-$ algebra $\Lambda$ has a NCCR, then $\height\dn{p}=n$ and $\Lambda_{\dn{p}}$ is a symmetric $R_{\dn{p}}$-order for any $\dn{p}\in\Max R$.

(2) Any NCCR of a $d$-CY$^-$ algebra $\Lambda$ is $d$-CY.

(3) Any reflexive tilting module over a $d$-CY algebra $\Lambda$ gives a NCCR.}

\vskip.5em{\sc Proof }
(1) Since $R$ is domain and $\ref\Lambda$ is non-empty, we have $0\in\Supp{}_R\Lambda$. Thus the structure morphism $R\to\Lambda$ is injective. Now the assertion follows from \XCA(3) and \XCB.

(2) Let $\Gamma$ be a NCCR of $\Lambda$. Fix $\dn{p}\in\Max R$. Since $\Gamma_{\dn{p}}$ is reflexive equivalent to a symmetric $R_{\dn{p}}$-order $\Lambda_{\dn{p}}$ by (1) and \XHA(1)(iii), it is a symmetric $R_{\dn{p}}$-order with $\gl\Gamma_{\dn{p}}=\height\dn{p}$ by \XBD. Thus $\Gamma$ is $d$-CY.

(3) Let $T$ be a reflexive tilting $\Lambda$-module with $\Gamma:=\endm_\Lambda(T)$. Then $\Gamma$ is also $d$-CY by \XCA(1). Thus $\Gamma_{\dn{p}}$ is an $R_{\dn{p}}$-order with $\gl\Gamma_{\dn{p}}=\height\dn{p}$ for any $\dn{p}\in\Max R$ by \XCB. We only have to show that $T$ is a height one generator of $\Lambda$. For any $\dn{p}\in\Spec R$ with $\height\dn{p}=1$, we have $\depth{}_{R_{\dn{p}}}T_{\dn{p}}=1$ and $\gl\Lambda_{\dn{p}}=1$. Thus $T_{\dn{p}}$ is a projective tilting $\Lambda_{\dn{p}}$-module by \XBC\ and \XBG. Hence it is a generator.\rule{5pt}{10pt}

%We have an exact sequence $0\to\Gamma_{\dn{p}}\to T_0{}_{\dn{p}}\to T_1{}_{\dn{p}}\to0$, which splits by \XHF(1). Denote by $\ggg:\ref\Gamma\to\ref\Lambda$ the inverse of the equivalence $\fff:\ref\Lambda\to\ref\Gamma$ in \XHA. Then $M_{\dn{p}}=\ggg(\Gamma_{\dn{p}})\in\add\ggg(T_{\dn{p}})=\add N_{\dn{p}}$.\rule{5pt}{10pt}

\vskip.5em
In the rest of this section, we assume $d\le3$. 
First we consider the relationship between reflexive tilting modules and tilting modules of projective dimension at most one over a 3-CY algebra $\Lambda$.
%By \XBC(1), any reflexive tilting $\Lambda$-modules has projective dimension at most one. The following proposition gives a sufficient condition for the converse. In general, the converse does not hold. 

\vskip.5em{\bf Proposition \XHD\ }{\it
(1) Assume that $\dim R\le 3$ and $\Lambda_{\dn{p}}$ is an $R_{\dn{p}}$-order with $\gl\Lambda_{\dn{p}}=\height\dn{p}$ for any $\dn{p}\in\Spec R$.

\strut\kern1em
(i) Any $M\in\ref\Lambda$ satisfies $\pd{}_\Lambda M\le1$. 

\strut\kern1em
(ii) Any reflexive tilting $\Lambda$-module has projective dimension at most one. Conversely, if $\Lambda_{\dn{p}}$ is Morita equivalent to a local ring for any $\dn{p}\in\Spec R$ with $\height\dn{p}\le2$, then any tilting $\Lambda$-module of projective dimension at most one is reflexive.

(2) We use the notation in \XCN. Then $R:=S^G$ is an isolated singularity if and only if $G$ acts freely on $K^d\backslash\{0\}$. In this case, reflexive tilting modules over $\Lambda:=S*G$ are exactly tilting $\Lambda$-module of projective dimension at most one.}

\vskip.5em{\sc Proof }
(1)(i) For any $\dn{p}\in\Spec R$, we have $\depth{}_{R_{\dn{p}}}M_{\dn{p}}\ge\min\{2,\height{\dn{p}}\}$ and $\gl\Lambda_{\dn{p}}=\height\dn{p}$. Thus $\pd{}_{\Lambda_{\dn{p}}}M_{\dn{p}}\le1$ by \XBC, and we have $\pd{}_\Lambda M\le1$.

(ii) Let $T$ be a tilting $\Lambda$-module with $\pd{}_\Lambda T\le1$. For any $\dn{p}\in\Spec R$, it follows from \XBG\ that $T_{\dn{p}}$ is a tilting $\Lambda_{\dn{p}}$-module. If $\height\dn{p}=3$, then we have $\depth{}_{R_{\dn{p}}}T_{\dn{p}}\ge2$ by \XBC. Now assume $\height\dn{p}\le2$. Since $\Lambda_{\dn{p}}$ is Morita equivalent to a local ring, any tilting $\Lambda_{\dn{p}}$-module is projective. Thus we have $\depth{}_{R_{\dn{p}}}T_{\dn{p}}\ge\height\dn{p}$. Since $T$ satisfies the S$_2$ condition, it is reflexive.

(2) For the first assertion, we refer to [IY;8.2]. Now we show the second assertion. Fix any $\dn{p}\in\Spec R\backslash\Max R$. Since $S_{\dn{p}}$ is a CM module over a regular local ring $R_{\dn{p}}$, it is a free $R_{\dn{p}}$-module. Thus the assertion follows from (1)(ii) since $(S*G)_{\dn{p}}=\endm_R(S)_{\dn{p}}=\endm_{R_{\dn{p}}}(S_{\dn{p}})$ (e.g. proof of \XCN) is Morita equivalent to $R_{\dn{p}}$.\rule{5pt}{10pt}

\vskip.5em
Now we put $\Lambda:=S*G$ for $G=\langle{\rm diag}(1,-1,-1)\rangle$. Then $S^G$ is not an isolated singularity, and $\Lambda$ is isomorphic to a complete tensor product $K[[x_1]]\widehat{\otimes}_K(K[[x_2,x_3]]*H)$ for $H:=\langle{\rm diag}(-1,-1)\rangle$. For any $T\in\tilt_1(K[[x_2,x_3]]*H)$, we have $K[[x_1]]\widehat{\otimes}_KT\in\tilt_1\Lambda$. This is not reflexive if $T$ is not projective. Thus $\Lambda$ has non-reflexive tilting modules of projective dimension at most one.

Recall that a module-finite $R$-algebra $\Lambda$ is called an {\it isolated singularity} if $\gl\Lambda_{\dn{p}}=\height\dn{p}$ for any $\dn{p}\in\Spec R\backslash\Max R$. 
To prove derived equivalence of different NCCR, we need an easy lemma on depth, and a relationship between depth and vanishing of $\ext^1$ for isolated singularities.

\vskip.5em{\bf Lemma \XHE\ }{\it
Let $R$ be local and $0\to X_t\stackrel{f_t}{\to}X_{t-1}\stackrel{f_{t-1}}{\to}\cdots\stackrel{f_3}{\to}X_2\stackrel{f_2}{\to}X_1\stackrel{f_1}{\to}X_0\to0$
an exact sequence with $X_0\in\flmod\Lambda$ for some $t\ge0$. If $\depth X_i\ge i$ for any $i>0$, then $X_0=0$.}

\vskip.5em{\sc Proof }
Put $Y_i:=\Im f_{i}$ for $1\le i\le t$. Inductively, we will show $\depth Y_i\ge i$. This is true for $i=t$. Now we assume $\depth Y_{i+1}\ge i+1$ and consider the exact sequence
$0\to Y_{i+1}\to X_i\to Y_i\to0$
with $\depth Y_{i+1}\ge i+1$ and $\depth X_i\ge i$. Applying $\hom_\Lambda(\Lambda/J_\Lambda,-)$, we see that $\depth Y_i\ge i$. In particular, $X_0=Y_1\in\flmod\Lambda$ satisfies $\depth X_0\ge0$, i.e. $\hom_\Lambda(\Lambda/J_\Lambda,X_0)=0$, and hence $X_0=0$.\rule{5pt}{10pt}

%For example, any $R$-algebra with $\gl\Lambda=d=\dim R$ (e.g. $d$-CY algebra) is an isolated singularity.

\vskip.5em{\bf Lemma \XHF\ }{\it
Assume that $R$ is local with $\dim R=3$ and $\Lambda$ is an isolated singularity. Let $M,N\in\ref\Lambda$.

(1) For any $\dn{p}\in\Spec R\backslash\Max R$, we have that $M_{\dn{p}}$ is a projective $\Lambda_{\dn{p}}$-module.

(2) $\ext^i_\Lambda(M,X)$ ($i>0$) has finite length for any $X\in\mod\Lambda$.

(3) If $\depth\hom_\Lambda(M,N)\ge3$, then $\ext^1_\Lambda(M,N)=0$.}

\vskip.5em{\sc Proof }
(1) follows from the S$_2$ condition $\depth{}_{R_{\dn{p}}}M_{\dn{p}}\ge\min\{2,\height\dn{p}\}$ and \XBC. (2) follows from (1). Now we show (3). Consider the exact sequence $0\to\Omega M\to P\to M\to0$ where $P$ is projective. Then
$0\to\hom_\Lambda(M,N)\to\hom_\Lambda(P,N)\to\hom_\Lambda(\Omega M,N)\stackrel{}{\to}\ext^1_\Lambda(M,N)\to0$
is an exact sequence with $\depth\hom_\Lambda(M,N)\ge3$, $\depth\hom_\Lambda(P,N)\ge2$ and $\depth\hom_\Lambda(\Omega M,N)\ge2$ (e.g. \XBD(1)). Since $\ext^1_\Lambda(M,N)$ has finite length by (2), it follows that $\ext^1_\Lambda(M,N)=0$ by \XHE.\rule{5pt}{10pt}

%Then this is obtained from a complex \[{\bf X}:0\to X_1\stackrel{g}{\to}X_0\stackrel{f}{\to}M_2\to0\] with $X_i\in\add M_1$ by applying $\fff:=\hom_\Lambda(M_1,-)$. Since $\depth\Ker g\ge2$ and $\hom_\Lambda(M_1,\Ker g)=0$ holds, we obtain $\Ker g=0$. (2) We have an exact sequence $0\to X_1\stackrel{g^\prime}\to\Ker f\to C_1\to0$ with a right $(\add M_2)$-approximation $g^\prime$. By , $C_1$ has finite length. Since both $X_1$ and $\Ker f$ are reflexive, we obtain $C_1=0$ by \XHE. Thus we have an exact sequence \[0\to X_1\stackrel{g}{\to}X_0\stackrel{f}{\to}M_2\to C_0\to0\] with $C_0\in\flmod\Lambda$. By $\depth\Gamma_2\ge3$ and \XHF(3), $\ext^1_\Lambda(M_2,M_2)=0$ holds. Thus \[\hom_\Lambda({\bf X},M_2):0\to\endm_\Lambda(M_2)\stackrel{f\cdot}{\to}\hom_\Lambda(X_0,M_2)\stackrel{g\cdot}{\to}\hom_\Lambda(X_1,M_2)\to0\] is exact by \XHE(2). (3) Since we have a commutative diagram \[\begin{diag}\hom_\Lambda({\bf X},M_2):&0&\to&\endm_\Lambda(M_2)&\stackrel{f\cdot}{\longrightarrow}&\hom_\Lambda(X_0,M_2)&\stackrel{g\cdot}{\longrightarrow}&\hom_\Lambda(X_1,M_2)&\to&0\\ &&&\downarrow^{\fff_{M_2,M_2}}&&\wr\downarrow^{\fff_{X_0,M_2}}&&\wr\downarrow^{\fff_{X_1,M_2}}\\ \hom_{\Gamma_1}({\bf P},U):&0&\to&\endm_{\Gamma_1}(U)&\longrightarrow&\hom_{\Gamma_1}(P_0,U)&\longrightarrow&\hom_{\Gamma_1}(P_1,U)&\to&\ext^1_{\Gamma_1}(U,U)&\to&0 \end{diag}\] of exact sequences, we obtain $\ext^1_{\Gamma_1}(U,U)=0$ and $\endm_{\Gamma_1}(U)=\Gamma_2$.\rule{5pt}{10pt}

\vskip.5em
Now we can prove the main result in this section, where we say that a $\Gamma$-module $N$ is {\it rigid} if $\ext^1_\Gamma(N,N)=0$. It generalizes results in [I4], where $\Gamma$ is assumed to be an order which is an isolated singularity and $M$ is assumed to be a CM $\Gamma$-module.

\vskip.5em{\bf Theorem \XHG\ }{\it
Let $R$ be a normal Gorenstein domain with $\dim R\le3$ and $\Lambda$ a module-finite algebra. For $M_i\in\ref\Lambda$, put $\Gamma_i:=\endm_\Lambda(M_i)$ and $U:=\hom_\Lambda(M_1,M_2)$.

(1) If $\Gamma_1$ is a NCCR of $\Lambda$ and $\Gamma_2$ is an order, then $U$ is a reflexive rigid $\Gamma_1$-module with $\endm_{\Gamma_1}(U)=\Gamma_2$.

(2) If $\Gamma_i$ ($i=1,2$) is a NCCR of $\Lambda$, then $U$ is a reflexive tilting $\Gamma_1$-module with $\endm_{\Gamma_1}(U)=\Gamma_2$.}

\vskip.5em{\sc Proof }
(1) Since we have a reflexive equivalence $\hom_\Lambda(M_1,-):\ref\Lambda\to\ref\Gamma_1$ by \XHA(1)(iii), we have $U\in\ref\Gamma_1$ and $\endm_{\Gamma_1}(U)=\Gamma_2$. Fix $\dn{p}\in\Spec R$. If $\height\dn{p}\le2$, then $U_{\dn{p}}\in\ref(\Gamma_1)_{\dn{p}}$ and $\gl(\Gamma_1)_{\dn{p}}=\height\dn{p}$ imply that $U_{\dn{p}}$ is a projective $(\Gamma_1)_{\dn{p}}$-module. In particular, we have $\ext^1_{(\Gamma_1)_{\dn{p}}}(U_{\dn{p}},U_{\dn{p}})=0$. Now we assume $\height\dn{p}=3$. Since $\depth\endm_{(\Gamma_1)_{\dn{p}}}(U_{\dn{p}})=3$ and $(\Gamma_1)_{\dn{p}}$ is an isolated singularity, we have $\ext^1_{(\Gamma_1)_{\dn{p}}}(U_{\dn{p}},U_{\dn{p}})=0$ by \XHF(3). Thus we have $\ext^1_{\Gamma_1}(U,U)=0$.

%Let \[0\to P_1\to P_0\to U\to0\] be a projective resolution. Applying $\hom_{\Gamma_1}(-,U)$, we have an exact sequence \[0\to\endm_{\Gamma_1}(U)\to\hom_{\Gamma_1}(P_0,U)\to\hom_{\Gamma_1}(P_1,U)\to\ext^1_{\Gamma_1}(U,U)\to0.\] Since we have a reflexive equivalence $\hom_\Lambda(M_1,-):\ref\Lambda\to\ref\Gamma_1$ by \XHA(3), we have $\endm_{\Gamma_1}(U)=\Gamma_2$. Thus we have $\depth\endm_{\Gamma_1}(U)=3$ and $\depth\hom_{\Gamma_1}(P_i,U)\ge2$. Since $\ext^1_{\Gamma_1}(U,U)$ has finite length by \XHF(2), we have $\ext^1_{\Gamma_1}(U,U)=0$ by \XHE.

(2) Since $U=\hom_{\Lambda^{\op}}(M_2^*,M_1^*)$, we have $\endm_{\Gamma_2^{\op}}(U)=\Gamma_1$ and $\ext^1_{\Gamma_2^{\op}}(U,U)=0$ by (1). We have $\pd{}_{\Gamma_2^{\op}}U\le1$ by \XHD(1)(i). Let $0\to Q_1\to Q_0\to U\to0$ be a projective resolution of the $\Gamma_2^{\op}$-module $U$. Applying $\hom_{\Gamma_2^{\op}}(-,U)$, we have an exact sequence $0\to\Gamma_1\to\hom_{\Gamma_2^{\op}}(Q_0,U)\to\hom_{\Gamma_2^{\op}}(Q_1,U)\to0$ of $\Gamma_1$-modules. Since $\hom_{\Gamma_2^{\op}}(Q_i,U)\in\add{}_{\Gamma_1}U$ holds, $U$ is a tilting $\Gamma_1$-module.\rule{5pt}{10pt}

\vskip.5em
As a direct consequence, we have the following desired result.

\vskip.5em{\bf Corollary \XHH\ }{\it
Let $R$ be a normal Gorenstein domain with $\dim R\le3$ and $\Lambda$ a module-finite algebra.

(1) Then all NCCR of $\Lambda$ are derived equivalent.

(2) Assume that $R$ is complete local. If $M_1$ and $M_2$ are $\Lambda$-modules giving NCCR, then $M_1$ and $M_2$ have the same number of non-isomorphic indecomposable summands.}

\vskip.5em{\sc Proof }
(1) follows directly from \XHG(2). For (2), we use that the Grothendieck group of $\endm_\Lambda(M_i)$ has a basis consisting of isoclasses of indecomposable projective modules, and that derived equivalences preserve Grothendieck groups.\rule{5pt}{10pt}

\vskip.5em
Note that \XHG\ also gives a correspondence between $\Lambda$-modules giving rise to NCCR and a subset of the reflexive tilting $\Gamma_1$-modules. For 3-CY$^-$ algebras we have the following improvement of this.

\vskip.5em{\bf Theorem \XHI\ }{\it
Let $\Lambda$ be 3-CY$^-$ and $M$ a $\Lambda$-module giving a NCCR $\Gamma:=\endm_\Lambda(M)$ and $\fff:=\hom_{\Lambda}(M,-):\ref\Lambda\to\ref\Gamma$ the induced reflexive equivalence.

(1) $\fff$ gives a one-one correspondence between $\Lambda$-modules giving NCCR and reflexive tilting $\Gamma$-modules.

(2) $\fff$ gives a one-one correspondence between reflexive $\Lambda$-modules whose endomorphism rings are orders and reflexive rigid $\Gamma$-modules.

(3) For $N\in\ref\Lambda$, $\endm_\Lambda(N)$ is an order if and only if $N$ is a direct summand of some $\Lambda$-module giving a NCCR.

(4) $\Lambda$ has a generator giving a NCCR.}

\vskip.5em{\sc Proof }
(1) Since $\fff$ is full and faithful, we only have to show that $\fff$ gives a surjective map, in view of \XHG(2). Take any reflexive tilting $\Gamma$-module $T$. Choose $N\in\ref\Lambda$ such that $\fff(N)=T$. By \XHC(2), $\Gamma$ is 3-CY. By \XHC(3), $T$ gives a NCCR of $\Gamma$. By \XHA(2), $N$ gives a NCCR of $\Lambda$.

(2) Since any reflexive rigid module is a direct summand of a reflexive tilting module by the Bongartz completion \XBH, the assertion follows by \XHG(1) and (1), and using that if $\endm_\Lambda(M)$ is an order and $N$ is a direct summand of $M$, then $\endm_\Lambda(N)$ is an order.

(3) The `if' part is obvious, and the `only if' part follows by (1)(2) and the Bongartz completion \XBH.

(4) Since $\endm_\Lambda(\Lambda)=\Lambda$ is an order, the assertion follows by (3).\rule{5pt}{10pt}

%\vskip.5em{\bf Corollary }{\it Let $\Lambda$ be a 3-CY$^-$ algebra with a NCCR. Then the following conditions for $N\in\ref\Lambda$ are equivalent. (1) $\endm_\Lambda(N)$ is an order. (2) $N$ is a direct summand of some $\Lambda$-module giving a NCCR.} \vskip.5em{\sc Proof } (1)$\Rightarrow$(2) is obvious. We will show (2)$\Rightarrow$(1). Let $\Gamma:=\endm_\Lambda(M)$ be a NCCR and $\fff:=\hom_\Lambda(M,-):\ref\Lambda\to\ref\Gamma$ a reflexive equivalence. Then $\fff(N)$ is a partial tilting $\Gamma$-module by \XHE. Using the Bongartz completion in \XBH, $\fff(N)$ is a direct summand of some classical tilting module $T$. By \XHI, there exists $L\in\ref\Lambda$ giving a NCCR such that $T=\fff(L)$. Since $\fff$ is an equivalence, $N$ is a direct summand of $L$.\rule{5pt}{10pt}

\vskip.5em
Choosing $M:=\Lambda$ in \XHI\ for a 3-CY algebra, we have the following remarkable relationship between tilting modules and NCCR.

\vskip.5em{\bf Corollary \XHJ\ }{\it
Let $\Lambda$ be 3-CY.

(1) $\Lambda$-modules giving NCCR are exactly reflexive tilting $\Lambda$-modules.

(2) Reflexive $\Lambda$-modules whose endomorphism rings are orders are exactly reflexive rigid $\Lambda$-modules.

(3) Let $T$ be a reflexive tilting $\Lambda$-module and $\Gamma:=\endm_\Lambda(T)$. Then the reflexive euivalence $\hom_\Lambda(T,-):\ref\Lambda\to\ref\Gamma$ gives a one-one correspondence between reflexive tilting $\Lambda$-modules and reflexive tilting $\Gamma$-modules.}

\vskip.5em
The following corollary shows that reflexive tilting modules over 3-CY algebras are closed under taking tensor products and homomorphisms. This is a quite peculiar property of 3-CY algebras. Especially (3) below gives another explanation of \XGC(4).

\vskip.5em{\bf Corollary \XHK\ }{\it
Let $\Lambda$ be 3-CY. 

(1) For any reflexive tilting $\Lambda$-modules $T_1$ and $T_2$ and
$\Gamma_i:=\endm_\Lambda(T_i)$, we have that $U:=\hom_\Lambda(T_1,T_2)$ is a reflexive tilting $\Gamma_1$-module with $\endm_{\Gamma_1}(U)=\Gamma_2$.

(2) $(-)^*\simeq\hom_\Lambda(-,\Lambda^*)$ gives a one-one correspondence between right reflexive tilting modules and left reflexive tilting modules.

(3) For any reflexive tilting $\Lambda^{\op}$-module $T_1$ and $\Lambda$-module $T_2$ and $\Gamma_i:=\endm_\Lambda(T_i)$, it follows that $U:=(T_1\otimes_\Lambda T_2)^{**}$ is a reflexive tilting $\Gamma_1$-module with $\endm_{\Gamma_1}(U)=\Gamma_2$.}

\vskip.5em{\sc Proof }
(1) Immediate from \XHJ(3).

(2) $\Lambda^*$ is a reflexive tilting $\Lambda$-module by \XHJ(1). We only have to put $T_2:=\Lambda^*$ in (1).

(3) Since $(T_1\otimes_\Lambda T_2)^{**}=\hom_\Lambda(T_2,T_1^*)^*$ holds, the assertion follows from (1) and (2).\rule{5pt}{10pt}

\vskip.5em
Van den Bergh raised the following question in [Va2;4.4]: {\it If $\Lambda$ has a NCCR, then does there exist a CM $\Lambda$-module giving a NCCR?} We give a positive answer for isolated singularities.

\vskip.5em{\bf Proposition \XHL\ }{\it
For any 3-CY$^-$ algebra $\Lambda$ which is an isolated singularity, then $\Lambda$ has a NCCR if and only if $\Lambda$ has a CM generator giving a NCCR.}

\vskip.5em{\sc Proof }
We only have to show the `only if' part. By \XHI(4), $\Lambda$ has a generator $M$ giving a NCCR. Since $\Lambda$ is an isolated singularity, we have $\ext^1_\Lambda(M,M)=0$ by \XHF(3), in particular, $\ext^1_\Lambda(M,\Lambda)=0$ holds. Taking localization and applying \XCD(5)(ii), we have $M\in\cm\Lambda$ since $M$ is reflexive.\rule{5pt}{10pt}

\vskip.5em
It was shown in [I4] that modules giving NCCR is closed related to maximal $1$-orthogonal modules introduced in [I3].
We call $M\in\cm\Lambda$ a {\it maximal $1$-orthogonal $\Lambda$-module} if
{\small\[\add M=\{X\in\cm\Lambda\ |\ \ext^1_\Lambda(M,X)=0\}=\{X\in\cm\Lambda\ |\ \ext^1_\Lambda(X,M)=0\}.\]}
%For example, $S$ is a maximal $1$-orthogonal $S^G$-module [I3,4].
If $\Lambda$ is an order, then any maximal $1$-orthogonal $\Lambda$-module $M$ satisfies $\Lambda\oplus\Lambda^*\in\add M$. Using [I4;5.2.1] and \XHL, we immediately obtain the following result for complete regular local $R$. Later we shall show in \XHR\ that this is valid for arbitrary $R$.

\vskip.5em{\bf Corollary \XHM\ }{\it
Let $R$ be complete regular local and $\Lambda$ a 3-CY$^-$ algebra which is an isolated singularity. Then $\Lambda$ has a NCCR if and only if $\Lambda$ has a maximal $1$-orthogonal module.}

\vskip.5em{\sc Proof }
Since $\Lambda$ is 3-CY$^-$, we have $\Lambda^*\simeq\Lambda$ as $\Lambda$-modules by \XCB. Now assume that $M\in\cm\Lambda$ is a generator.
It was shown in [I4;5.2.1] that $M$ gives a NCCR of $\Lambda$ if and only if $M$ is maximal $1$-orthogonal. Thus the assertion follows from \XHL.\rule{5pt}{10pt}

\vskip.5em
We now investigate the relationship between rigid modules and NCCR for 3-CY$^-$ algebras, in particular for those which are isolated singularities. The following will be useful.

\vskip.5em{\bf Lemma \XHN\ }{\it
Let $M\in\ref\Lambda$ be a generator, $\Gamma:=\endm_\Lambda(M)$ and $\fff:=\hom_\Lambda(M,-):\ref\Lambda\to\ref\Gamma$ a reflexive equivalence (\XBD(2)). Then we have a functorial monomorphism
$\ext^1_\Gamma(\fff(X),\fff(Y))\subseteq\ext^1_\Lambda(X,Y)$
for any $X,Y\in\ref\Lambda$.}

\vskip.5em{\sc Proof }
Since $M$ is a generator, there exists an exact sequence $0\to X_1\to M_0\stackrel{f}{\to}X\to0$ with a right $(\add M)$-approximation $f$. Then $0\to\fff(X_1)\to\fff(M_0)\to\fff(X)\to0$ is an exact sequence in $\mod\Gamma$ where $\fff(M_0)$ is projective. We have the following exact commutative diagram:
{\small\[\begin{diag}
\hom_\Gamma(\fff(M_0),\fff(Y))&\longrightarrow&\hom_\Gamma(\fff(X_1),\fff(Y))&\longrightarrow&\ext^1_\Gamma(\fff(X),\fff(Y))&\longrightarrow&0\\
\parallel&&\parallel\\
\hom_\Lambda(M_0,Y)&\longrightarrow&\hom_\Lambda(X_1,Y)&\longrightarrow&\ext^1_\Lambda(X,Y)
\end{diag}\]}
\vskip-1em
Thus we have a monomorphism $\ext^1_\Gamma(\fff(X),\fff(Y))\subseteq\ext^1_\Lambda(X,Y)$.\rule{5pt}{10pt}

\vskip.5em{\bf Theorem \XHO\ }{\it
Let $\Lambda$ be a 3-CY$^-$ algebra with a NCCR and $N\in\ref\Lambda$. Then (1)$\Rightarrow$(2)$\Leftrightarrow$(3) holds. If $\Lambda$ is an isolated singularity, then (1)--(3) are equivalent.

(1) $N$ is rigid.

(2) $\endm_\Lambda(N)$ is an order.

(3) $N$ is a direct summand of some $\Lambda$-module giving a NCCR.}

\vskip.5em{\sc Proof }
(2)$\Leftrightarrow$(3) follows by \XHI(3). If $\Lambda$ is an isolated singularity, then we can show (2)$\Rightarrow$(1) by taking localization and applying \XHF(3).

(1)$\Rightarrow$(2) Assume that $N$ is rigid. By \XHI(4), $\Lambda$ has a generator $M$ giving a NCCR. Put $\Gamma:=\endm_\Lambda(M)$ and let $\fff:=\hom_\Lambda(M,-):\ref\Lambda\to\ref\Gamma$ be a reflexive equivalence. It follows from \XHN\ that $\ext^1_\Gamma(\fff(N),\fff(N))=0$. Thus $\fff(N)$ is a reflexive rigid $\Gamma$-module. By \XHI(2), $\endm_\Lambda(N)$ is an order.\rule{5pt}{10pt}

\vskip.5em
We illustrate with the following example.

\vskip.5em{\bf Theorem \XHP\ }{\it
Let $K$ be a field of characteristic zero, $G$ a finite subgroup of $\SL_3(K)$, $S:=K[[x,y,z]]$ and $g$ the number of irreducible representations of $G$.

(1) For any $S^G$-module $M$ giving a NCCR, the number of non-isomorphic indecomposable direct summands of $M$ is exactly $g$.

(2) Any reflexive rigid $S^G$-module $N$ is a direct summand of an $S^G$-module giving a NCCR. Thus the number of non-isomorphic indecomposable direct summands of $N$ is at most $g$.}

\vskip.5em{\sc Proof }
$S$ gives a NCCR $S*G$ of a 3-CY$^-$-algebra $S^G$, which has
exactly $g$ non-isomorphic indecomposable direct summands [I4].
Thus (1) follows by \XHH(2). Now (2) follows by (1) and \XHO(1)$\Rightarrow$(3).\rule{5pt}{10pt}

\vskip.5em
We say that a reflexive rigid $\Lambda$-module $N$ {\it maximal rigid} if $L\in\ref\Lambda$ and $\ext^1_\Lambda(N\oplus L,N\oplus L)=0$ imply $L$ is in $\add N$ (c.f. [GLS]).

\vskip.5em{\bf Corollary \XHQ\ }{\it
Let $\Lambda$ be a 3-CY$^-$ algebra which is an isolated singularity and has a NCCR.

(1) $M$ in $\ref\Lambda$ gives a NCCR if and only if it is maximal rigid.

(2) Any reflexive equivalence $\ref\Lambda\to\ref\Gamma$ gives a one-one correspondence between rigid $\Lambda$-modules and rigid $\Gamma$-modules.}

\vskip.5em{\sc Proof }
(1) This follows directly from the equivalence of \XHO(1) and (3).

(2) By \XHO(1) and (2), rigidity depends only on the endomorphism ring.\rule{5pt}{10pt}

\vskip.5em
We have the following generalization of \XHM\ for arbitrary $R$.

\vskip.5em{\bf Theorem \XHR\ }{\it
Let $\Lambda$ a 3-CY$^-$ algebra which is an isolated singularity. 

(1) CM $\Lambda$-modules giving NCCR are exactly maximal $1$-orthogonal $\Lambda$-modules.

(2) $\Lambda$ has a NCCR if and only if $\Lambda$ has a maximal $1$-orthogonal module.}

\vskip.5em{\sc Proof }
By \XHL, we only have to show (1).

(i) Assume that $M\in\cm\Lambda$ is maximal $1$-orthogonal.
Put $\Gamma:=\endm_\Lambda(M)$ and $\fff:=\hom_\Lambda(M,-)$.
%Take $L\in\ref\Lambda$ such that $M\oplus L$ is rigid.
%Since $M$ is a generator, we have $L\in\cm\Lambda$.
%Thus $\Ext^1_\Lambda(M,L)=0$ implies $L\in\add M$. This shows that $M$ is maximal rigid. By \XHQ(1), $M$ gives a NCCR of $\Lambda$.
Take an exact sequence $0\to\Omega M\to P\to M\to0$ where $P$ is projective.
Applying $\hom_\Lambda(-,M)$, we have an exact sequence $0\to\Gamma\to\hom_\Lambda(P,M)\to\hom_\Lambda(\Omega M,M)\to\Ext^1_\Lambda(M,M)=0$.
Since $\hom_\Lambda(P,M)\in\cm R$ and $\hom_\Lambda(\Omega M,M)\in\ref R$, we have that $\Gamma$ is an $R$-order.

Since any $\dn{p}\in\Max R$ satisfies $\height\dn{p}=3$ by \XCB(1), we only have to show $\gl\Gamma\le 3$.
For any $X\in\mod\Gamma$, take a projective resolution $Q_1\stackrel{f}{\to}Q_0\to X\to0$ where $Q_i$ is projective.
By Yoneda's Lemma on $\add M$, there exists $g\in\hom_\Lambda(M_1,M_0)$ such that $M_i\in\add M$ and $f=\fff(g)$.
Put $Y:=\Ker g$ and take an exact sequence $0\to Z\to M_2\stackrel{h}{\to}Y\to0$ with a right $(\add M)$-approximation $h$.
Then $Z\in\cm\Lambda$ holds.
Applying $\fff$, we have an exact sequence $\fff(M_2)\stackrel{\fff(h)}{\longrightarrow}\fff(Y)\to\Ext^1_\Lambda(M,Z)\to\Ext^1_\Lambda(M,M_2)$. Since $\fff(h)$ is surjective and $M$ is rigid, we have $\Ext^1_\Lambda(M,Z)=0$. Thus $Z\in\add M$.
Consequently, $\pd{}_\Gamma X\le 3$ holds since we have a projective resolution
$0\to\fff(Z)\to\fff(M_2)\to Q_1\stackrel{f}{\to}Q_0\to X\to0$.

(ii) Assume that $M\in\cm\Lambda$ gives a NCCR $\Gamma:=\endm_\Lambda(M)$.
Put $\fff:=\hom_\Lambda(M,-)$. By \XHQ(1), $M$ is maximal rigid. Since $\Lambda$ is 3-CY$^-$ and $M\in\cm\Lambda$, we have that $M\oplus\Lambda$ and $M\oplus\Lambda^*$ are rigid. In particular, we have $\Lambda\oplus\Lambda^*\in\add M$.

Take $X\in\cm\Lambda$ such that $\Ext^1_\Lambda(X,M)=0$. Take a projective resolution $P_1\to P_0\to X^*\to0$ of a $\Lambda^{\op}$-module $X^*$.
Applying $(-)^*$, we have an exact sequence $0\to X\to P_0^*\to P_1^*$ with $P_i^*\in\add\Lambda^*\subset\add M$.
Applying $\fff$, we have an exact sequence $0\to\fff(X)\to\fff(P_0^*)\to\fff(P_1^*)$ with projective $\Gamma$-modules $\fff(P_i^*)$.
Since $\gl\Gamma=3$, we have $\pd{}_\Gamma\fff(X)\le1$.
Thus we can take a projective resolution $0\to Q_1\to Q_0\to\fff(X)\to0$.
By Yoneda's Lemma on $\add M$, there exists a complex $0\to M_1\to M_0\to X\to0$ with $M_i\in\add M$ such that $0\to\fff(M_1)\to\fff(M_0)\to\fff(X)\to0$ is isomorphic to the above projective resolution. Since $M$ is a generator, we have that $0\to M_1\to M_0\to X\to0$ is exact. Since $\Ext^1_\Lambda(X,M)=0$ by our assumption, we have $X\in\add M$.

Take $Y\in\cm\Lambda$ such that $\Ext^1_\Lambda(M,Y)=0$. 
Since $M^*\in\cm\Lambda^{\op}$ gives a NCCR and $\Ext^1_{\Lambda^{\op}}(Y^*,M^*)=0$, we have $Y^*\in\add M^*$ and $Y\in\add M$.
Consequently, $M$ is maximal $1$-orthogonal.\rule{5pt}{10pt}

\vskip.5em
%We now relate maximal rigid modules to the maximal $1$-orthogonal modules. 
While a maximal $1$-orthogonal module is maximal rigid, the converse does not hold in general even if it is CM. For example, the simple singularity $\Lambda:=K[[x_1,x_2,x_3,x_4]]/(x_1^{2n+1}+x_2^2+x_3^2+x_4^2)$ does not have non-projective rigid CM modules [Yo]. Thus $\Lambda$ is a maximal rigid $\Lambda$-module, which is not maximal $1$-orthogonal. It seems to be difficult to know when the converse holds. Now we show that some kind of converse holds for 3-CY$^-$ algebras.

%We call $M\in\ref\Lambda$ a {\it maximal $1$-orthogonal $\Lambda$-module} if {\small\[\add M=\{X\in\ref\Lambda\ |\ \ext^1_\Lambda(M,X)=0\}=\{X\in\ref\Lambda\ |\ \ext^1_\Lambda(X,M)=0\}.\]} For example, $S$ is a maximal $1$-orthogonal $S^G$-module [I4]. By definition, any maximal $1$-orthogonal $\Lambda$-module $M$ satisfies $\Lambda\in\add M$. Thus $M$ must be in $\cm\Lambda$.
%For $M\in\ref\Lambda$, put $\Gamma:=\endm_\Lambda(M)$ and \[\cm M:=\{X\in\ref\Lambda\ |\ \hom_\Lambda(M,X)\in\cm\Gamma\}.\] Of course, $\cm\Lambda$ coincides with the usual defined one. We have an equivalence\[\fff:=\hom_\Lambda(M,-):\cm M\to\cm\Gamma.\]
%\vskip.5em{\bf\XHS\ Proposition }{\it Any maximal $1$-orthogonal object is maximal rigid.}\vskip.5em{\sc Proof } Assume $M\oplus N$ is rigid for $N\in\ref\Lambda$. Then $\ext^1_\Lambda(N,\Lambda)=0$ implies $N\in\cm\Lambda$. Since $\ext^1_\Lambda(M,N)=0$ holds, we have $N\in\add M$ by the maximal $1$-orthogonality of $M$.\rule{5pt}{10pt}

\vskip.5em{\bf Proposition \XHS\ }{\it
Let $\Lambda$ be a 3-CY$^-$ algebra which is an isolated singularity and has a NCCR. Then $M\in\ref\Lambda$ is maximal rigid if and only if there exists a reflexive equivalence $\fff:\ref\Lambda\to\ref\Gamma$ such that $\fff(M)$ is a maximal $1$-orthogonal $\Gamma$-module.}

\vskip.5em{\sc Proof }
Let $M$ be maximal rigid. Then $M$ gives a NCCR $\Gamma:=\endm_\Lambda(M)$ by \XHQ(1). Then $\fff:=\hom_\Lambda(M,-):\ref\Lambda\to\ref\Gamma$ is a reflexive equivalence. Since $\Gamma$ is an order with $\gl\Gamma=3$, we have $\cm\Gamma=\add\Gamma$ by taking localization and applying \XBC. Thus $\fff(M)=\Gamma$ is a maximal $1$-orthogonal $\Gamma$-module. The other implication follows from \XHQ(2) since a maximal $1$-orthogonal module is maximal rigid.\rule{5pt}{10pt}

{\footnotesize
\begin{center}
{\bf References}
\end{center}

[AV] M. Artin, J.-L. Verdier: Reflexive modules over rational double points. Math. Ann. 270 (1985), no. 1, 79--82.

[A1] M. Auslander: Representation dimension of Artin algebras. Lecture notes, Queen Mary College, London, 1971.

[A2] M. Auslander: Functors and morphisms determined by objects. Representation theory of algebras (Proc. Conf., Temple Univ., Philadelphia, Pa., 1976), pp. 1--244. Lecture Notes in Pure Appl. Math., Vol. 37, Dekker, New York, 1978. 

[A3] M. Auslander: Rational singularities and almost split sequences. Trans. Amer. Math. Soc. 293 (1986), no. 2, 511--531. 

[AB] M. Auslander, M. Bridger: Stable module theory. Memoirs of the American Mathematical Society, No. 94, American Mathematical Society, Providence, R.I. 1969 146 pp.

[APR] M. Auslander, M. I. Platzeck, I. Reiten: Coxeter functors without diagrams.  Trans. Amer. Math. Soc. 250 (1979), 1--46.

[AR1] M. Auslander, I. Reiten: Applications of contravariantly finite subcategories. Adv. Math. 86 (1991), no. 1, 111--152.

[AR2] M. Auslander, I. Reiten: $D$Tr-periodic modules and functors. Representation theory of algebras (Cocoyoc, 1994), 39--50, CMS Conf. Proc., 18, Amer. Math. Soc., Providence, RI, 1996.

[AS] M. Auslander, S. Smalo: Almost split sequences in subcategories.  J. Algebra  69  (1981), no. 2, 426--454.

[BT] R. Berger, R. Taillefer: Poincare-Birkhoff-Witt Deformations of Calabi-Yau Algebras, \newline arXiv:math.RT/0610112.

[BB] A. Bjorner, F. Brenti: Combinatorics of Coxeter groups. Graduate Texts in Mathematics, 231. Springer, New York, 2005.

[Boc] R. Bocklandt: Graded Calabi Yau Algebras of dimension 3, arXiv:math.RA/0603558.

[BO] A. Bondal, D. Orlov: Semiorthogonal decomposition for algebraic varieties, arXiv:alg-geom/9506012.

[Bon1] K. Bongartz: Tilted algebras. Representations of algebras (Puebla, 1980), pp. 26--38, Lecture Notes in Math., 903, Springer, Berlin-New York, 1981. 

[Bon2] K. Bongartz: Algebras and quadratic forms. J. London Math. Soc. (2) 28 (1983), no. 3, 461--469.

[Bra] A. Braun, in preparation.

[BGS] K.A. Brown, I.G. Gordon, C.H. Stroppel: Cherednik, Hecke and quantum algebras as free Frobenius and Calabi-Yau extensions, arXiv:math.RT/0607170.

[Bri1] T. Bridgeland: Flops and derived categories. Invent. Math. 147 (2002), no. 3, 613--632.

[Bri2] T. Bridgeland: t-structures on some local Calabi-Yau varieties. J. Algebra 289 (2005), no. 2, 453--483.

[Bri3] T. Bridgeland: Stability conditions on a non-compact Calabi-Yau threefold. Comm. Math. Phys. 266 (2006), no. 3, 715--733.

[Bri4] T. Bridgeland: Stability conditions and Kleinian singularities, arXiv:math.AG/0508257.

[BKR] T. Bridgeland, A. King, M. Reid: The McKay correspondence as an equivalence of derived categories. J. Amer. Math. Soc. 14 (2001), no. 3, 535--554.

[BH] W. Bruns, J. Herzog: Cohen-Macaulay rings. Cambridge Studies in Advanced Mathematics, 39. Cambridge University Press, Cambridge, 1993.

[BMR1] A. Buan, R. Marsh, I. Reiten: Cluster-tilted algebras. Trans. Amer. Math. Soc. 359 (2007), no. 1, 323--332.

[BMR2] A. Buan, R. Marsh, I. Reiten: Cluster mutation via quiver representations, to appear in Comm. Math. Helv., arXiv:math.RT/0412077.

[BMRRT] A. Buan, R. Marsh, M. Reineke, I. Reiten, G. Todorov: Tilting theory and cluster combinatorics. Adv. Math. 204 (2006), no. 2, 572--618.

[CCS] P. Caldero, F. Chapoton, R. Schiffler: Quivers with relations arising from clusters ($A\sb n$ case). Trans. Amer. Math. Soc. 358 (2006), no. 3, 1347--1364.

%[CHU] F. Coelho, D. Happel, L. Unger: Complements to partial tilting modules. J. Algebra 170 (1994), no. 1, 184--205. 

[CRe] C. W. Curtis, I. Reiner: Methods of representation theory. Vol. I. With applications to finite groups and orders. Pure and Applied Mathematics. A Wiley-Interscience Publication. John Wiley \& Sons, Inc., New York, 1981.

[CRo] J. Chuang, R. Rouquier, in preparation.

[EG] E. G. Evans, P. Griffith: Syzygies. London Mathematical Society Lecture Note Series, 106. Cambridge University Press, Cambridge, 1985.

[F] H. Foxby: Bounded complexes of flat modules. J. Pure Appl. Algebra 15 (1979), no. 2, 149--172.

[FZ1] S. Fomin, A. Zelevinsky: Cluster algebras. I. Foundations. J. Amer. Math. Soc. 15 (2002), no. 2, 497--529.

[FZ2] S. Fomin, A. Zelevinsky: Cluster algebras. II. Finite type classification. Invent. Math. 154 (2003), no. 1, 63--121.

[GLS] C. Geiss, B. Leclerc, J. Schr\"oer: Rigid modules over preprojective algebras. Invent. Math. 165 (2006), no. 3, 589--632.

[G] V. Ginzburg, Calabi-Yau algebras, arXiv:math.AG/0612139.

[GN1] S. Goto, K. Nishida: Finite modules of finite injective dimension over a Noetherian algebra. J. London Math. Soc. (2) 63 (2001), no. 2, 319--335. 

[GN2] S. Goto, K. Nishida: Towards a theory of Bass numbers with application to Gorenstein algebras. Colloq. Math. 91 (2002), no. 2, 191--253. 

[GR] A. L. Gorodentsev, A. N. Rudakov: Exceptional vector bundles on projective spaces. Duke Math. J. 54 (1987), no. 1, 115--130.

[GV] G. Gonzalez-Sprinberg, J.-L. Verdier: Construction geometrique de la correspondance de McKay. Ann. Sci. Ecole Norm. Sup. (4) 16 (1983), no. 3, 409--449 (1984). 

[Hae] J. Haefner: On Gorenstein, Frobenius and symmetric orders. Nova J. Algebra Geom. 2 (1993), no. 4, 315--359.

[Har1] R. Hartshorne: Residues and duality. Lecture notes of a seminar on the work of A. Grothendieck, given at Harvard 1963/64. With an appendix by P. Deligne. Lecture Notes in Mathematics, No. 20 Springer-Verlag, Berlin-New York 1966.

[Har2] R. Hartshorne: Algebraic geometry. Graduate Texts in Mathematics, No. 52. Springer-Verlag, New York-Heidelberg, 1977.

[Hap] D. Happel: Triangulated categories in the representation theory of finite-dimensional algebras. London Mathematical Society Lecture Note Series, 119. Cambridge University Press, Cambridge, 1988.

[HPR1] D. Happel, U. Preiser, C. M. Ringel: Vinberg's characterization of Dynkin diagrams using subadditive functions with application to $D{\rm Tr}$-periodic modules. Representation theory, II (Proc. Second Internat. Conf., Carleton Univ., Ottawa, Ont., 1979), pp. 280--294, Lecture Notes in Math., 832, Springer, Berlin, 1980.

[HPR2] D. Happel, U. Preiser, C. M. Ringel: Binary polyhedral groups and Euclidean diagrams. Manuscripta Math. 31 (1980), no. 1-3, 317--329.

[HU1] D. Happel, L. Unger: Almost complete tilting modules. Proc. Amer. Math. Soc.  107  (1989),  no. 3, 603--610. 

%[HU1] D. Happel, L. Unger: Partial tilting modules and covariantly finite subcategories. Comm. Algebra 22 (1994), no. 5, 1723--1727.

[HU2] D. Happel, L. Unger; On a partial order of tilting modules. Algebr. Represent. Theory 8 (2005), no. 2, 147--156.

[He] J. Herzog: Ringe mit nur endlich vielen Isomorphieklassen von maximalen, unzerlegbaren Cohen-Macaulay-Moduln. Math. Ann. 233 (1978), no. 1, 21--34.

[Hu]  J. E. Humphreys: Reflection groups and Coxeter groups. Cambridge Studies in Advanced Mathematics, 29. Cambridge University Press, Cambridge, 1990.

[IU] A. Ishii, H. Uehara: Autoequivalences of derived categories on the minimal resolutions of $A\sb n$-singularities on surfaces. J. Differential Geom. 71 (2005), no. 3, 385--435.

[I1] O. Iyama: $\tau$-categories I: Ladders, Algebr. Represent. Theory 8 (2005), no. 3, 297--321.

[I2] O. Iyama: $\tau$-categories II: Nakayama pairs and Rejective subcategories, Algebr. Represent. Theory 8 (2005), no. 4, 449--477.

[I3] O. Iyama: Higher dimensional Auslander-Reiten theory on maximal orthogonal subcategories, to appear in Adv. Math., arXiv:math.RT/0407052.

[I4] O. Iyama: Auslander correspondence, to appear in Adv. Math., arXiv:math.RT/0411631.

[I5] O. Iyama: Maximal orthogonal subcategories of triangulated categories satisfying Serre duality, Mathematisches Forschungsinstitut Oberwolfach Report, no. 6 (2005), 353--355.

[IY] O. Iyama, Y. Yoshino: Mutations in triangulated categories and rigid Cohen-Macaulay modules, arXiv:math.RT/0607736.

[K1] B. Keller: Deriving DG categories. Ann. Sci. Ecole Norm. Sup. (4) 27 (1994),  no. 1, 63--102.

[K2] B. Keller: On the construction of triangle equivalences. Derived equivalences for group rings, 155--176, Lecture Notes in Math., 1685, Springer, Berlin, 1998.

[K3] B. Keller: On triangulated orbit categories. Doc. Math. 10 (2005), 551--581.

[KR] B. Keller, I. Reiten: Cluster-tilted algebras are Gorenstein and stably Calabi-Yau, to appear in Adv. Math., arXiv:math.RT/0512471.

[KV] M. Kapranov, E. Vasserot: Kleinian singularities, derived categories and Hall algebras. Math. Ann. 316 (2000), no. 3, 565--576.

[MRZ] R. Marsh, M. Reineke, A. Zelevinsky: Generalized associahedra via quiver representations. Trans. Amer. Math. Soc. 355 (2003), no. 10, 4171--4186.

[Ma] H. Matsumura: Commutative ring theory. Cambridge Studies in Advanced Mathematics, 8. Cambridge University Press, Cambridge, 1986.

[Mc] J. McKay: Graphs, singularities, and finite groups. The Santa Cruz Conference on Finite Groups (Univ. California, Santa Cruz, Calif., 1979), pp. 183--186,  Proc. Sympos. Pure Math., 37, Amer. Math. Soc., Providence, R.I., 1980. 

[MY] J. Miyachi, A. Yekutieli: Derived Picard groups of finite-dimensional hereditary algebras.  Compositio Math.  129  (2001),  no. 3, 341--368.

[Ra] M. Ramras: Maximal orders over regular local rings of dimension two. Trans. Amer. Math. Soc. 142 1969 457--479.

[Re] I. Reiner: Maximal orders. London Mathematical Society Monographs, No. 5. Academic Press, London-New York, 1975.

[RV1] I. Reiten, M. Van den Bergh: Two-dimensional tame and maximal orders of finite representation type. Mem. Amer. Math. Soc. 80 (1989).

[RV2] I. Reiten, M. Van den Bergh:  Noetherian hereditary abelian categories satisfying Serre duality. J. Amer. Math. Soc.  15  (2002), no. 2, 295--366.

[RR] I. Reiten, C. Riedtmann: Skew group algebras in the representation theory of Artin algebras. J. Algebra 92 (1985), no. 1, 224--282.

[Ri1] J. Rickard: Morita theory for derived categories. J. London Math. Soc. (2) 39 (1989), no. 3, 436--456.

[Ri2] J. Rickard, private communication.

[RS] C. Riedtmann, A. Schofield: On a simplicial complex associated with tilting modules. Comment. Math. Helv. 66 (1991), no. 1, 70--78.

[Rud] A. N. Rudakov: Markov numbers and exceptional bundles on $P\sp 2$. Math. USSR-Izv. 32 (1989), no. 1, 99--112.

[Rum1] W. Rump: Non-commutative regular rings. J. Algebra 243 (2001), no. 2, 385--408.

[Rum2] W. Rump: Non-commutative Cohen-Macaulay rings.  J. Algebra  236  (2001),  no. 2, 522--548.

[RZ] R. Rouquier, A. Zimmermann: Picard groups for derived module categories. Proc. London Math. Soc. (3)  87  (2003),  no. 1, 197--225.

[ST] P. Seidel, R. Thomas: Braid group actions on derived categories of coherent sheaves. Duke Math. J. 108 (2001), no. 1, 37--108.

[U] L. Unger: Schur modules over wild, finite-dimensional path algebras with three simple modules. J. Pure Appl. Algebra 64 (1990), no. 2, 205--222.

[Va1] M. Van den Bergh: Three-dimensional flops and noncommutative rings. Duke Math. J. 122 (2004), no. 3, 423--455.

[Va2] M. Van den Bergh: Non-commutative crepant resolutions. The legacy of Niels Henrik Abel, 749--770, Springer, Berlin, 2004.

[Va3] M. Van den Bergh: Introduction to super potentials, Mathematisches Forschungsinstitut Oberwolfach Report, no. 6 (2005), 396--398.

[Ve] J.-L. Verdier: Des categories derivees des categories abeliennes. Asterisque No. 239 (1996). 

[Yo] Y. Yoshino: Cohen-Macaulay modules over Cohen-Macaulay rings. London Mathematical Society Lecture Note Series, 146. Cambridge University Press, Cambridge, 1990. 

[Ye1] A. Yekutieli: Dualizing complexes over noncommutative graded algebras. J. Algebra 153 (1992), no. 1, 41--84.

[Ye2] A. Yekutieli: Dualizing complexes, Morita equivalence and the derived Picard group of a ring. J. London Math. Soc. (2) 60 (1999), no. 3, 723--746.}

\vskip.5em{\footnotesize
{\sc Graduate School of Mathematics, Nagoya University
Chikusa-ku, Nagoya, 464-8602
Japan}

{\it iyama@math.nagoya-u.ac.jp}

\vskip.5em
{\sc Institutt for matematiske fag
NTNU
7491 Trondheim, Norway}

{\it idunr@math.ntnu.no}
}
\end{document}